\documentclass[a4paper,english]{smfart} %Right-numbered equations

\usepackage[T1]{fontenc} 
\usepackage[utf8]{inputenc}
\usepackage[english]{babel}
\usepackage{newpxtext} %Palatino
	\usepackage{eulervm} %Euler style

\usepackage{amssymb,textcomp,mathrsfs,braket,commath,relsize}

\usepackage{stmaryrd}  %Avoid clash with extpfeil
\expandafter\def\csname opt@stmaryrd.sty\endcsname
{only,shortleftarrow,shortrightarrow}

\usepackage{bm} %To use many alphabets
	\SetSymbolFont{stmry}{bold}{U}{stmry}{m}{n} %Avoids meaningless warning for bold math symbols

\makeatletter %Print 2020AMSsubjclass (instead of 2010)
\@namedef{subjclassname@2020}{
	\textup{2020} Mathematics Subject Classification}
\makeatother

\usepackage{mathtools}
    \mathtoolsset{showonlyrefs} %Numbering only for referred equations 
	
\usepackage{smfthm} 
	\NumberTheoremsAs{subsubsection}\SwapTheoremNumbers

\usepackage{etoolbox} 
    \newcommand{\addQEDstyle}[2]{\AtBeginEnvironment{#1}{\pushQED{\qed}\renewcommand{\qedsymbol}{#2}}
    \AtEndEnvironment{#1}{\popQED}} %Symbol at the end of environment: call as \addQEDstyle{environmentname}{symbolname}
    \addQEDstyle{rema}{$\triangle$}
    \addQEDstyle{remark*}{$\triangle$} %Triangles at the end of remarks
	\addQEDstyle{exem}{$\triangle$}
    \addQEDstyle{exem*}{$\triangle$} %Triangles at the end of examples
	
\usepackage{tikz}
	\usetikzlibrary{graphs,shapes.geometric,decorations.pathreplacing,calligraphy} %For graphs, pentagons,curly braces
	\tikzset{
		nodes={circle,draw},
		every loop/.style={}, %removes arrow head from loops
	}
	
\usepackage{float}

\AtBeginDocument{\def\MR#1{}} %No MR number in bibliography
\apptocmd{\sloppy}{\hbadness 10000\relax}{}{} %For badboxes in .bbl (uses etoolbox)

\makeatletter %To print 2020AMSsubjclass
	\@namedef{subjclassname@2020}{
		\textup{2020} Mathematics Subject Classification}
\makeatother

\usepackage{mymacros}

\usepackage[marginpar]{todo}

\usepackage[colorlinks,psdextra,pdfencoding=auto,backref=page]{hyperref}
	\hypersetup{colorlinks,linkcolor={red!50!black},filecolor={green!50!black},urlcolor={blue!80!black},citecolor={blue!50!black}}

\begin{document}
	\frontmatter 
	
	\title[(Quasi) root systems and hyperplane arrangements]{A colourful classification of (quasi) root systems and hyperplane arrangements} %[header]{main title}

	\author[G.~Rembado]{Gabriele Rembado} 

	\address[G.~Rembado]{Hausdorff Centre for Mathematics, University of Bonn, 60 Endenicher Allee, D-53115 Bonn (Germany)}
		\email{gabriele.rembado@hcm.uni-bonn.de}
		
% 	\subjclass{17B22,52C35}
	
	
	\thanks{The author is supported by the Deutsche Forschungsgemeinschaft (DFG, German Research Foundation) under Germany Excellence Strategy - GZ 2047/1, Projekt-ID 390685813.}
	
	\begin{abstract}
		We introduce a class of graphs with coloured edges to encode subsystems of the classical root systems, which in particular classify them up to equivalence. 
		We further use the graphs to describe root-kernel intersections, as well as restrictions of root (sub)systems on such intersections, generalising the regular part of a Cartan subalgebra.
		We also consider a slight variation to encode the hyperplane arrangements only, showing there is a unique noncrystallographic arrangement that arises.
		Finally, a variation of the main definition leads to elementary classifications of closed and Levi root subsystems.
	\end{abstract}

	{\let\newpage\relax\maketitle} %No page break after title

	%Tablo of content
	\setcounter{tocdepth}{1}  %Stop at sections
		\tableofcontents

	\mainmatter 
	
	\section*{Introduction}

	Let $G$ be a connected complex reductive Lie group, and $T \sse G$ a maximal torus (such as the group of invertible diagonal matrices inside $\GL_n(\mb C)$).
	Denote $\mf g = \Lie(G)$ the Lie algebra of $G$, and let $\mf t = \Lie(T) \sse \mf g$ be the associated Cartan subalgebra.
	
	Previous papers~\cite{doucot_rembado_tamiozzo_2022_local_wild_mapping_class_groups_and_cabled_braids,doucot_rembado_2022_topology_of_irregular_isomonodromy_times_on_a_fixed_pointed_curve} introduced ``wild'' mapping class groups of type $\mf g$, associated with the root system $\Phi_{\mf g} = \Phi(\mf g, \mf t) \sse \mf t^{\dual}$ and some further data: the missing piece comes from moduli spaces arising in complex algebraic geometry, parametrising isomorphism/gauge classes of meromorphic connections on principal $G$-bundles over a Riemann surface, i.e. the \emph{de Rham} spaces of the wild nonabelian Hodge correspondence~\cite{sabbah_1999_harmonic_metrics_and_connections_with_irregular_singularities,biquard_boalch_2004_wild_nonabelian_hodge_theory_on_curves}.
	See~\cite{doucot_rembado_tamiozzo_2022_local_wild_mapping_class_groups_and_cabled_braids,doucot_rembado_2022_topology_of_irregular_isomonodromy_times_on_a_fixed_pointed_curve,boalch_doucot_rembado_2022_twisted_local_wild_mapping_class_groups} for details and references, in particular the relation with the wild character varieties (the \emph{Betti} spaces) and the admissible deformations of wild Riemann surfaces~\cite{boalch_2014_geometry_and_braiding_of_stokes_data_fission_and_wild_character_varieties}---as well as the reviews~\cite{boalch_2012_hyperkaehler_manifolds_and_nonbelian_hodge_theory_of_irregular_curves,boalch_2017_wild_character_varieties_meromorphic_hitchin_systems_and_dynkin_graphs}.
	The papers ~\cite{rembado_2018_quantisation_of_moduli_spaces_and_connections,rembado_2019_simply_laced_quantum_connections_generalising_kz,rembado_2020_symmetries_of_the_simply_laced_quantum_connections_and_quantisation_of_quiver_varieties,felder_rembado_2023_singular_modules_for_affine_lie_algebras,andersen_malusa_rembado_2022_genus_one_complex_quantum_chern_simons_theory,andersen_malusa_rembado_2021_sp1_symmetric_hyperkaehler_quantisation} instead considered the ``quantum'' side of this story---the quantisation of the de Rham spaces.
	
	\vspace{5pt}
	
	The upshot are spaces $\bm B \sse \mf t^p$ for an integer $p \geq 1$, coming with a product decomposition
	\begin{equation}
		\label{eq:factors_deformation_space}
		\bm B = \prod_{i = 1}^p \bm B_i \, , \qquad \bm B_i \sse \mf t \, ,
	\end{equation}
	and the pure local ``wild'' mapping class group is then the fundamental group 
	\begin{equation*}
		\Gamma = \pi_1(\bm B) \simeq \prod_i \pi_1(\bm B_i) \, .
	\end{equation*}
	Now each factor in~\eqref{eq:factors_deformation_space} is a ``restricted'' root-hyperplane complement, in the following sense.
	If $\Phi' \sse \Phi$ is an inclusion of root systems, we consider the flat 
	\begin{equation}
		\label{eq:root_kernel}
		U \ceqq \Ker(\Phi') = \bigcap_{\Phi'} \Ker(\alpha) \sse \mf t \, ,
	\end{equation}
	of the root-hyperplane arrangement of $\Phi \sse \mf t^{\dual}$, and then the complement of the remaining hyperplanes:
	\begin{equation}
		\label{eq:restricted_complement}
		\bm B(\Phi',\Phi) \ceqq U \cap \bigcup_{\Phi \sm \Phi'} \bigl( \mf t \sm \Ker(\alpha) \bigr) \sse U \, .
	\end{equation}
	When nonempty,~\eqref{eq:restricted_complement} is the complement of the hyperplane arrangement (in $U$) obtained upon restriction of the additional roots, i.e.
	\begin{equation}
		\label{eq:restricted_complement_2}
		\bm B(\Phi',\Phi) = U \, \bigsm \, \bigcup_{\Phi \sm \Phi'} \Ker \bigl( \eval[1]{\alpha}_U \bigr) \, .
	\end{equation}
	Then there exists an increasing sequence 
	\begin{equation}
		\label{eq:root_system_sequence}
		\Phi_1 \sse \dm \sse \Phi_{p+1} \ceqq \Phi_{\mf g} \, ,
	\end{equation}
	of root (sub)systems, such that $\bm B_i = \bm B_i(\Phi_i,\Phi_{i+1})$ for $i \in \set{1,\dc,p}$.
	
	The generic case corresponds to the inclusion $\vn \sse \Phi_{\mf g}$ (for $p = 1$), which brings about the regular part
	\begin{equation*}
		\bm B(\vn,\Phi_{\mf g}) = \mf t_{\reg} = \mf t \, \bigsm \, \bigcup_{\Phi_{\mf g}} \Ker(\alpha) \sse \mf t \, ,
	\end{equation*}
	of the Cartan subalgebra.
	Hence this example of local wild mapping class group, whose role in 2d gauge theory was first understood in~\cite{boalch_2002_g_bundles_isomonodromy_and_quantum_weyl_groups}, leads to pure (generalised/Artin--Tits) braid groups of type $\mf g$~\cite{brieskorn_1971_die_fundamentalgruppe_des_raumes_der_regulaeren_orbits_einer_endlichen_komplexen_spiegelungsgruppe,deligne_1972_les_immeubles_des_groupes_de_tresses_generalises,brieskorn_saito_1972_artin_gruppen_und_coxeter_gruppen}.
	
	\vspace{5pt}
	
	The problem we consider here is thus to classify the hyperplane arrangements of the ``restricted'' system
	\begin{equation}
		\label{eq:restricted_system}
		\eval[1]{\Phi}_U \ceqq \Set{ \eval[1]{\alpha}_U | \alpha \in \Phi \sm \Phi'} \sse U^{\dual} 
	\end{equation}
	of linear functionals.
	
	\vspace{5pt}
	
	Note~\eqref{eq:restricted_system} is a symmetric subset, but it need \emph{not} be a root system; and further even the resulting hyperplane arrangement is not that of a root system in general---i.e. it is not \emph{crystallographic}.
	Moreover it is \emph{not} true that the resulting Coxeter/reflection group is obtained by restricting (to $U$) the elements of the Weyl group preserving~\eqref{eq:root_kernel}---even in type $A$.
	
	In view of the aforementioned obstructions, an explicit description of~\eqref{eq:restricted_system} in principle relies on a study of:
	\begin{enumerate}
		\item root subsystems $\Phi \sse \Phi_{\mf g}$;
		
		\item for any inclusion $\Phi' \sse \Phi$ of two such, the set of equivalence classes of functionals $\alpha,\beta \in \Phi$ with respect to the relation
		\begin{equation*}
			\alpha \sim \beta \qquad \text{if} \qquad U \sse \Ker(\alpha - \beta) \, .
		\end{equation*}
	\end{enumerate}
	This latter is also the quotient set $\Phi \bs U^{\perp} \sse \mf t^{\dual} \bs U^{\perp}$, within the quotient vector space (of $\mf t^{\dual}$) modulo the annihilator 
	\begin{equation*}
		U^{\perp} = \Set{ \lambda \in \mf t^{\dual} | \eval[1]{\lambda}_U = 0 } \sse \mf t^{\dual} \, .
	\end{equation*}
	But in our situation $U^{\perp} = \mb C \Phi' \sse \mf t^{\dual}$, so we can equivalently consider the quotient set 
	\begin{equation*}
		\Phi \bs \mb C \Phi' \sse \mb C \Phi \bs \mb C \Phi' \, ,
	\end{equation*}
	which is thus naturally in bijection with~\eqref{eq:restricted_system}: we refer to either as the \emph{quotient} of $\Phi$, modulo $\Phi'$.\fn{
		E.g. if $\Phi' = \set{\pm \gamma}$ for some $\gamma \in \Phi$ then the quotient is simply the set of $\gamma$-strings in $\Phi$.}
	
	While a study of such quotients and hyperplane arrangements can be carried in general, the geometric setup of~\cite{doucot_rembado_tamiozzo_2022_local_wild_mapping_class_groups_and_cabled_braids,doucot_rembado_2022_topology_of_irregular_isomonodromy_times_on_a_fixed_pointed_curve} only brings about a particular type of root subsystems: namely the root systems of iterated centraliser of semisimple elements, viz. (reductive) Levi factors of parabolic Lie subalgebras of $\mf g$.
	These are thus the \emph{Levi} (root) subsystems, leading to the ``fission'' of~\cite{boalch_2007_quasi_hamiltonian_geometry_of_meromorphic_connections,boalch_2009_through_the_analytic_halo_fission_via_irregular_singularities}---in which case the complement~\eqref{eq:restricted_complement_2} is nonempty.
	
	\vspace{5pt}
	
	In this paper instead we consider the more abstract problem of studying~\eqref{eq:restricted_system} for \emph{all} root subsystems of a simple Lie algebra $\mf g$.
	This yields the general semisimple case by taking direct products, and adding an Abelian factor (in the reductive case) does not modify the root system---only the ambient Cartan subalgebra.
	More precisely, we propose an elementary classification of quotients of the classical root subsystems, using certain graphs with coloured edges which retain more information than the Dynkin diagrams, and which do \emph{not} rely on a choice of base of simple roots.

	\vspace{5pt}
	
	Note~\cite{oshima_2006_a_classification_of_subsystems_of_a_root_system} also deals with the classification of irreducible root (sub)systems: in particular the statement of Cor.~\ref{cor:root_system_classification} is coherent with the tables in \S~10 of op.~cit; but this will just come as a corollary of the main construction, and we further get to \emph{non}crystallographic systems/arrangements by considering quotients.
	Also we should point out our quotients are different from those of~\cite{howlett_1980_normalizers_of_parabolic_subgroups_of_reflection_groups}, but the normalisers of parabolic subgroups considered there should be relevant to the study of the reflection groups of~\eqref{eq:restricted_system}.
	
	Finally, using an (invariant) scalar product $( \cdot \mid \cdot) \cl \mf t \otimes \mf t \to \mb C$, the vanishing locus of a subset $\Phi' \sse \mf t^{\dual}$ corresponds to the orthogonal subspace of its image under the vector space isomorphism $( \cdot \mid \cdot )^{\sharp} \cl \mf t^{\dual} \to \mf t$.
	Then restricting a linear functional $\lambda \in \mf t^{\dual}$ to $U = \Ker(\Phi') \sse \mf t$ is the same as taking the orthogonal projection of $( \cdot | \cdot )^{\sharp}(\lambda) \in \mf t$ onto $U$: in this viewpoint one is thus led to study projections of (co)roots onto intersections of other (co)roots' kernels, cf.~\cite{dijols_2019_projections_of_root_systems}.
	The viewpoint we take here instead does \emph{not} rely on an invariant product.
	
	\subsection*{Main results and layout}
	
	In \S~\ref{sec:type_A} we review the special linear case, which is the prototype for the other classical types.
	
	In \S~\ref{sec:crystallographs} we introduce \emph{bichromatic} graphs $(\mc G,c)$, i.e. graphs endowed with a red/green colouring $c \cl \mc G_1 \to \set{R,G}$ of each edge, and describe the main correspondence which maps them bijectively to symmetric subsets of the nonreduced root system $BC_n$---if $\mc G$ has $n \geq 1$ nodes, cf. \S~\ref{sec:appendix}.
	In particular there are ``classical'' graphs corresponding to root systems of type $A$, $B$, $C$, $D$ and $BC$, depicted in Ex.~\ref{ex:classical_crystallographs}.
	
	In \S~\ref{sec:main_correspondence} we introduce a special class of bichromatic graphs, which correspond to (crystallographic) root subsystems of $BC_n$ under the main correspondence above: they are thus called \emph{crystallo}graphs (see Def.~\ref{def:crystallographs}).
	
	\begin{enonce*}{Theorem (cf. Thm.~\ref{thm:crystallographs_classify_root_systems})}
		
		There is a canonical inclusion-preserving bijection between crystallographs on $n \geq 1$ nodes and root subsystems of $BC_n$.
	\end{enonce*}
	
	In \S~\ref{sec:crystallograph_classification} we classify crystallographs directly from their defining properties, proving the following.
	
	\begin{enonce*}{Theorem (cf. Thm.~\ref{thm:crystallograph_classification}, Lem.~\ref{lem:bipartite_graph_is_type_A} and Cor.~\ref{cor:root_system_classification})}
		
		All crystallographs are disjoint unions of ``classical'' ones, up to acting on root (sub)systems via the Weyl group; in particular no exceptional types arise from root subsystems of classical simple Lie algebras.
	\end{enonce*}

	In \S~\ref{sec:kernels} we use crystallographs to explicitly determine the vanishing locus of any root subsystem.
	
% 	\begin{enonce*}{Theorem (cf. Prop.~\ref{prop:kernels})}
% 		
% 		Let $\Phi \sse BC_n$ be a root subsystem, and $\mc G^{\Phi}$ the associated crystallograph; then $\Ker(\Phi)$ has a natural basis consisting of the set of type-$A$ connected components of $\mc G^{\Phi}$---including trivial ones.
% 	\end{enonce*}
	
	In \S~\ref{sec:quotients} we define ``quotients'' graphs $\mc G \bs \mc G'$ of nested crystallographs $\mc G' \sse \mc G$, with the explicit algorithm of Def.~\ref{def:quotient_graph}.
	We study the obstruction for such a quotient to be a crystallograph, and accordingly introduce a weakened notion of \emph{quasi}-crystallographs in Def.~\ref{def:quasi_crystallographs}
	The main statement about these two new classes of bichromatic graphs is the following.
	
	\begin{enonce*}{Theorem (cf. Thm.~\ref{thm:quotient_graphs_classify_restricted_arrangements}, Thm.~\ref{thm:quasi_crystallographs_classify_quotients}, and Prop.~\ref{prop:quasi_crystallographs_are_quotients})}
		
		Let $\Phi' \sse \Phi \sse BC_n$ be nested root systems, and $\mc G' = \mc G^{\Phi'}$, $\mc G = \mc G^{\Phi}$ the associated crystallographs.
		Then the restricted system~\eqref{eq:restricted_system} is associated with the quotient graph $\mc G \bs \mc G'$ under the main correspondence, and it is a quasi-crystallograph; and conversely all quasi-crystallographs arise from quotients of root systems.
	\end{enonce*}

	Further we completely classify this extended class of graphs, thereby explicitly describing all possible quotients of root systems.
	
	\begin{enonce*}{Theorem (cf. Cor.~\ref{cor:quasi_crystallographs_classification} and Rk.~\ref{rk:exotic_arrangement})}
		
		Quasi-crystallographs are disjoint unions of ``classical'' graphs, and two other ``exotic'' components---which do \emph{not} generically correspond to root systems.
	\end{enonce*}

	In \S~\ref{sec:arrangements} we modify the main definition to introduce \emph{projective} crystallographs, in Def.~\ref{def:projective_crystallographs}.
	They retain less information than crystallographs, since are used to (completely) encode the data of subsets of the root-hyperplane arrangement of type $B_n/C_n$; further a natural ``projectification'' operation $\mc G \mapsto \mb P(\mc G)$ produces such a graph from any crystallograph, and corresponds to taking the canonical projection $\pi \cl \mf t^{\dual} \sm \set{0} \to \mb P \bigl( \mf t^{\dual} \bigr)$ of root systems---as their kernels are dilation-independent.
	Then we adapt the previous arguments to prove the following.
	
	\begin{enonce*}{Theorem (cf. Cor.~\ref{cor:projective_crystallograph_classification} and Thm.~\ref{thm:projective_crystallograph_classificy_hyperplane_arrangements})}
		
		Projective crystallographs on $n \geq 1$ nodes are in natural inclusion-preserving bijection with root-hyperplane sub-arrangements of $B_n/C_n/BC_n$, and are all obtained by projectifying crystallographs.
		Moreover they are disjoint unions of ``classical'' projective crystallographs.
	\end{enonce*}

	Then we consider the compatibility of quotients and projectifications, to get to the desired classification of quotient/restricted hyperplane arrangements.
	
	\begin{enonce*}{Theorem (Cf. Lem.~\ref{lem:compatibility_quotients_projectifications} and Thm.~\ref{thm:projectified_quotients_classify_restricted_arrangements})}
		Let $\Phi' \sse \Phi \sse BC_n$ be nested root systems, and $\mc G' = \mb P (\mc G^{\Phi'})$, $\mc G = \mb P(\mc G^{\Phi})$ the associated projective crystallographs.
		Then the quotient $\mc G \bs \mc G'$ is the projectification of the quasi-crystallograph associated to the restricted system~\eqref{eq:restricted_system}, and thus controls the quotient/restricted hyperplane arrangement.
		
		Moreover all projectifications of quasi-crystallographs are disjoint unions of the ``classical'' ones, and of a \emph{single} exotic component.
	\end{enonce*}

	The latter exotic hyperplane arrangement is as follows (cf.~\cite{doucot_rembado_tamiozzo_2022_local_wild_mapping_class_groups_and_cabled_braids}).
	Let $r,s \geq 0$ be integers, and consider the standard complex coordinates $\bm z = (z_1,\dc,z_{r+s})$ on $\mb C^{r+s}$.
	Then the exotic arrangement contains the hyperplanes 
	\begin{equation*}
		H^{\pm}_{ij} = \Set{ \bm z \in \mb C^{r+s} | z_i \pm z_j = 0 } \sse \mb C^{r+s} \, , \qquad i \neq j \in \set{1,\dc,r+s} \, ,
	\end{equation*}
	viz. the root-hyperplanes of type $D_{r+s}$; and the hyperplanes 
	\begin{equation*}
		H_i = \Set{ \bm z \in \mb C^{r+s} | z_i = 0 } \sse \mb C^{r+s} \, , \qquad i \in \set{1,\dc,r} \, ,
	\end{equation*}
	viz. the root hyperplanes of type $B_r/C_r$---in the standard embedding $\mb C^r \hra \mb C^r \times \mb C^s$. 
	
	This arrangement thus plays a distinguished role within the theory of irreducible root systems.
	
	Finally, in \S\S~\ref{sec:closed_subystems} and~\ref{sec:levi_subsystems} we use yet another variation of the main correspondence to classify closed subsystems of roots, and further Levi subsystems.
		
	\vspace{5pt}
	
	Some standard notions, notations and conventions, used throughout the paper, are summarised in \S~\ref{sec:appendix}.
	The end of examples and remarks is signaled by a ``$\triangle$''.
	
	\section{Prototype: type \texorpdfstring{$A$}{A}}
	\label{sec:type_A}
	
	\subsection{}
	
	The main idea can be seen in type $A$: let $\mf g = \mf{sl}_n(\mb C)$, with root system $A_{n-1}$.
	
	If $\Phi \sse A_{n-1}$ is a root subsystem, for each $i \in \ul n \ceqq \set{1,\dc,n}$ we can consider the subset
	\begin{equation}
		\label{eq:type_A_partition_sets}
		I_i = \set{i} \cup \Set{ j \in \ul n | \alpha^-_{ij} \in \Phi } \sse \ul n \, ,
	\end{equation}
	where $\alpha^-_{ij} \cl V \to \mb C$ are the type-$A$ roots for $V = \mb C^n$ (cf. \S~\ref{sec:appendix}).
	The fact that $\Phi$ is closed under the root-hyperplane reflections of its elements implies the subsets~\eqref{eq:type_A_partition_sets} yield a partition of $\ul n$: 
	\begin{equation}
		\label{eq:type_A_partition}
		\ul n = \coprod_{i \in J} I_i \, , \qquad J \ceqq \Set{ \min (I_i) | i \in \ul n } \sse \ul n \, .
	\end{equation}
	In turn there is a root system decomposition
	\begin{equation}
		\label{eq:type_A_decomposition}
		\Phi \simeq \bigoplus_J A_{I_i} \, , \qquad A_{I_i} \ceqq \Set{ \alpha^-_{jk} | j,k \in I_i } \sse \Phi \, .
	\end{equation}
	
	\subsection{}
	
	This situation can be naturally encoded into a graph $\mc G = \mc G^{\Phi}$ on the set of nodes $\mc G_0 = \ul n$, by putting an (unoriented) straight edge $\set{j,k}$ if and only if $\alpha_{jk} \in \Phi$.
	Then~\eqref{eq:type_A_decomposition} is equivalent to the fact that $\mc G$ decomposes into a disjoint union 
	\begin{equation*}
		\mc G = \coprod_J \mc G(i) \, ,
	\end{equation*}
	of (sub)graphs $\mc G(i) \sse \mc G$ for $i \in J$, each of which is a complete graph on the set of nodes $I_i \sse \ul n$---without loop edges.
	Importantly \emph{all} root subsystems are thus parametrised by such graphs.
	
	Moreover one can use these latter to describe quotients, as follows.
	An inclusion $\Phi' \sse \Phi$ corresponds to a subgraph $\mc G'= \mc G^{\Phi'}$ of $\mc G$, on the same set of nodes, and the system of restricted roots~\eqref{eq:restricted_system} is determined by the missing edges.
	More precisely there may be pairs of connected components $\mc G'(i),\mc G'(j) \sse \mc G'$, for $i,j \in J$, such that there are nodes $i_0 \in \mc G'(i)_0$ and $j_0 \in \mc G'(j)_0$ which are adjacent in $\mc G$---i.e. $\set{i_0,j_0} \in \mc G_1 \sm \mc G'_1$; in the next picture the new straight edge (lying in the larger graph only) is depicted by a dashed line:
	\begin{figure}[H]
		\centering
		\begin{tikzpicture}
			\node (a) at (0,0) {};
			\node (b) [label=below:$i_0$] at (2,0) {};
			\node (c) at (1,sqrt 3) {};
			\node (d) [draw=none] at (1,-.5) {$\mc G'(i)$};
			\node (e) [label=below:$j_0$] at (4,0) {};
			\node (f) at (6,0) {};
			\node (g) at (6,2) {};
			\node (h) at (4,2) {};
			\node (i) [draw=none] at (5,-.5) {$\mc G'(j)$};
			\graph{
				(a) -- (b) -- (c) -- (a),
				(e) -- (f) -- (g) -- (h) -- (e),
				(e) -- (g),
				(f) -- (h),
				(b) --[dashed] (e)
			};
		\end{tikzpicture}
		\caption{Example: two connected components of $\mc G'$ are linked in $\mc G$.}
	\end{figure}
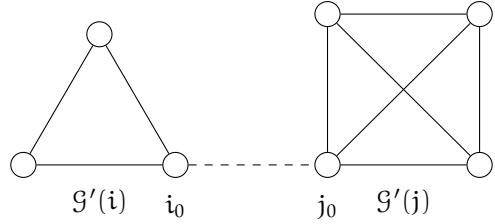

	Further the fact that $\Phi$ is a root (sub)system then implies all pairs of nodes in two such components are adjacent in the larger graph, i.e. there exists a connected component of $\mc G$ containing both $\mc G'(i)$ and $\mc G'(j)$ as subgraphs.
	This is because the Weyl group element corresponding to $\sigma^-_{ij} \in A_{n-1}$ maps $\alpha^-_{jk} \mapsto \alpha^-_{ki}$, for any triple $\set{i,j,k} \sse \ul n$ of distinct indices (see below more generally).	
	
	Then we can construct a ``quotient'' graph, denoted $\mc G \bs \mc G'$, with the following algorithm: its set of nodes is that of connected components of $\mc G'$, viz. the set $J' \sse \ul n$ associated with $\Phi'$ as in~\eqref{eq:type_A_partition}, and two components are adjacent if and only if they are ``fused'' in $\mc G$.
	
	It turns out $\mc G \bs \mc G'$ is still a disjoint union of complete graphs (see below more generally), so it encodes a type-$A$ root system inside $\mb C^{\abs{J'}}$; and further there is a canonical vector space isomorphism $\Ker(\Phi') \simeq \mb C^{\abs{J'}}$, identifying $A_{n-1}$ with the root system of $\mf{gl}(V)$---vanishing on the centre $\mb C \Id_V \sse \mf{gl}(V)$.
	
	Hence this graph-theoretic description yields in particular a proof of the following (known):	
	\begin{prop}[Cf.~\cite{doucot_rembado_tamiozzo_2022_local_wild_mapping_class_groups_and_cabled_braids}]
		\, \\[-10pt]
		\begin{enumerate}
			\item If $\Phi \sse A_{n-1}$ is a root subsystem, then it is isomorphic to a direct sum of irreducible type-$A$ root systems.
			
			\item If $\Phi' \sse \Phi \sse A_{n-1}$ is an inclusion of root systems, the quotient $\Phi \bs \Phi'$ is isomorphic to a root system of type $A$, and in particular $A_{n-1} \bs \Phi' \simeq A_{\abs{J}' - 1}$.
		\end{enumerate}
	\end{prop}
	
	\subsubsection{}
	
	A proof of a stronger statement, involving all classical simple Lie algebras, is given below: in brief a more general class of graphs (extending $\mc G' \sse \mc G$ and $\mc G \bs \mc G'$ from this section) leads to an explicit description of \emph{all} root subsystems, and of \emph{all} linear functionals that arise upon the restriction~\eqref{eq:restricted_system}.
	The main idea is to add a different type of straight edges to encode the new roots that arise in type $D$, and further loop edges to encode the short/long roots in type $B/C$.
	Finally, a variation thereof can be introduced to only encode the resulting hyperplane arrangements (cf. \S~\ref{sec:arrangements}), and to refine the classification to closed/Levi subsystems only (cf. \S\S~\ref{sec:closed_subystems} and~\ref{sec:levi_subsystems})---which is a void restriction in type $A$.
	
	\begin{rema}
		In particular in type $A$ the ``restricted'' root-hyperplane complement~\eqref{eq:restricted_complement_2} is identified with a product of spaces of (ordered) configurations of points in the complex affine line, whose fundamental groups are pure braid groups.
		If $\Phi = A_{n-1}$ there is a single component, and the fundamental group is the pure braid group on $\abs{J'}$ strands.
		See~\cite{doucot_rembado_tamiozzo_2022_local_wild_mapping_class_groups_and_cabled_braids} about the general case.
	\end{rema}	
	
	\section{Bichromatic graphs and symmetric sets of roots}
	\label{sec:crystallographs}
	
	\subsection{}
	
	Let $\set{R,G}$ be the set of colours ``red'' and ``green''.
	
	\begin{defi}
		A \emph{bichromatic} graph is a graph $\mc G = (\mc G_0,\mc G_1,m)$ equipped with a colour function $c \cl \mc G_1 \to \set{R,G}$, assigning colours to each edge, such that $c(e,m_e) \neq c(e',m_e')$ if $e = e' \sse \mc G_0$.
		
		A bichromatic \emph{subgraph} of $(\mc G,c)$ is a subgraph $\mc G' \sse \mc G$ equipped with the restricted colour function $c' = \eval[1]{c}_{\mc G'_1}$.\fn{
			At times the colour function will be (abusively) omitted.}
	\end{defi}
	
	\subsection{}
	
	The colouring condition is that any two edges incident at the same node(s) have different colours, so there are at most two such edges.
	We thus say a bichromatic graph is \emph{complete} if it has exactly two edges incident at each pair of (possibly coinciding) nodes.

	If the colour function is constant we speak of ``red/green graphs'' according to whether $c(\mc G_1) = \set R$ or $c(\mc G_1) = \set G$---this is now just an ordinary ``monochromatic'' graph.
	We denote $\mc G_R, \mc G_G \sse \mc G$ the red/green subgraphs of the bichromatic graph $\mc G$ obtained by keeping its red/green edges only, i.e. the graphs defined by
	\begin{equation*}
		(\mc G_{\bullet})_0 \ceqq \mc G_0 \, , \qquad (\mc G_{\bullet})_1 \ceqq c^{-1}(\bullet) \sse \mc G_1 \, , \qquad \bullet \in \set{R,G} \, .
	\end{equation*}
	Clearly all the information about $(\mc G,c)$ is contained in the pair $(\mc G_R,\mc G_G)$.
	
	Recall an ordinary graph is \emph{simply-laced} if it has no repeated edges and no loop edges.
	
	\begin{defi}
		A bichromatic graph $(\mc G,c)$ is \emph{simply-laced} if $(\mc G_R,\mc G_G)$ is a pair of simply-laced (monochromatic) graphs.
	\end{defi}
	
	\subsubsection{}
	
	This means $\mc G$ has no loop edges, and we denote $\ol{\mc G} \sse \mc G$ the simply-laced bichromatic graph obtained from $\mc G$ by removing its loop edges: it is the maximal simply-laced bichromatic subgraph of $\mc G$.
	In this terminology a complete bichromatic graph is thus \emph{not} simply-laced, so a \emph{simply-laced complete} bichromatic graph is rather a bichromatic graph with every possible straight edge---but no loop edges.
	
	\subsection{Correspondence with symmetric systems}
	
	Now for an integer $n \geq 1$ we use the incidence of a bichromatic graph $(\mc G,c)$ on the set of nodes $\mc G_0 = \ul n$ to store information about certain subsets of $BC_n$, in the following (main) correspondence.
	(Recall from \S~\ref{sec:appendix} that we consider the vector space $V = \mb C^n$ with canonical basis $(e_i)_{i \in \ul n}$).
	
	If $e = \set{i,j} \in \mc G_1$ is a straight edge, then $\Phi^{\mc G,c}$ contains:
	\begin{enumerate}
		\item the roots $\pm e^-_{ij}$, if $c(e) = R$;
		
		\item the roots $\pm e^+_{ij}$, if $c(e) = G$.
	\end{enumerate}
	
	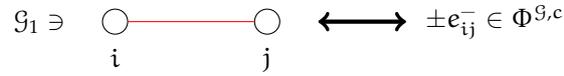
\begin{figure}[H]
		\centering
		\begin{tikzpicture}
			\node (a) [label=below:$i$] at (0,0) {};
			\node (b) [label=below:$j$] at (2,0) {};
			\node (c) [draw=none] at (2.5,0) {};
			\node (d) [draw=none] at (4,0) {};
			\node (e) [draw=none] at (5,0) {$\pm e^-_{ij} \in \Phi^{\mc G,c}$};
			\node (f) [draw=none] at (-1,0) {$\mc G_1 \ni$};
			\graph{
				(a) --[red] (b),
				(c) --[to-to,ultra thick] (d)
			};
		\end{tikzpicture}
		\caption{Main correspondence, straight edges (I).}
	\end{figure}

	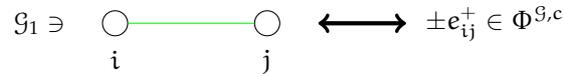
\begin{figure}[H]
		\centering
		\begin{tikzpicture}
			\node (a) [label=below:$i$] at (0,0) {};
			\node (b) [label=below:$j$] at (2,0) {};
			\node (c) [draw=none] at (2.5,0) {};
			\node (d) [draw=none] at (4,0) {};
			\node (e) [draw=none] at (5,0) {$\pm e^+_{ij} \in \Phi^{\mc G,c}$};
			\node (f) [draw=none] at (-1,0) {$\mc G_1 \ni$};
			\graph{
				(a) --[green] (b),
				(c) --[to-to,ultra thick] (d)
			};
		\end{tikzpicture}
		\caption{Main correspondence, straight edges (II).}
	\end{figure}
	
	Analogously if $l = \set{k} \in \mc G_1$ is a loop edge of $\mc G$ then $\Phi^{\mc G,c}$ contains:
	\begin{enumerate}
		\item the roots $\pm e_k$, if $c(l) = R$;
		
		\item the roots $\pm 2e_k$, if $c(l) = G$.
	\end{enumerate}
	
	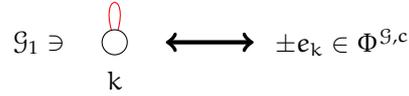
\begin{figure}[H]
		\centering
		\begin{tikzpicture}
			\node (a) [label=below:$k$] at (0,0) {};
			\node (c) [draw=none] at (.5,0) {};
			\node (d) [draw=none] at (2,0) {};
			\node (e) [draw=none] at (3,0) {$\pm e_k \in \Phi^{\mc G,c}$};
			\node (f) [draw=none] at (-1,0) {$\mc G_1 \ni$};
			\graph{
				(c) --[to-to,ultra thick] (d)
			};
			\path
			 (a) edge [loop above,red] (a);
		\end{tikzpicture}
		\caption{Main correspondence, loop edges (I).}
		\label{fig:main_correspondence_loop_edges_I}
	\end{figure}
	
	\begin{figure}[H]
		\centering
		\begin{tikzpicture}
			\node (a) [label=below:$k$] at (0,0) {};
			\node (c) [draw=none] at (.5,0) {};
			\node (d) [draw=none] at (2,0) {};
			\node (e) [draw=none] at (3,0) {$\pm 2e_k \in \Phi^{\mc G,c}$};
			\node (f) [draw=none] at (-1,0) {$\mc G_1 \ni$};
			\graph{
				(c) --[to-to,ultra thick] (d)
			};
			\path
			 (a) edge [loop above,green] (a);
		\end{tikzpicture}
		\caption{Main correspondence, loop edges (II).}
		\label{fig:main_correspondence_loop_edges_II}
	\end{figure}
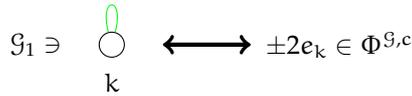

	By construction $\Phi^{\mc G,c}$ is closed under the involution $\alpha \mapsto -\alpha$ of $BC_n$---it is a \emph{symmetric} subset.
	Conversely any symmetric subset $\Phi \sse BC_n$ is associated with a unique bichromatic graph $(\mc G^{\Phi},c)$ on $n$ nodes by inverting the above prescription.

	By construction these correspondences are mutually inverse order-preserving bijections between bichromatic graphs on $n \geq 1$ nodes and symmetric subsets of the root system $BC_n$.

	\begin{exem}[Classical bichromatic graphs]	
		\label{ex:classical_crystallographs}
		
		The classical (irreducible) root systems yield the following bichromatic graphs: a simply-laced complete red graph for type $A$, a simply-laced complete bichromatic graph for type $D$, a simply-laced complete bichromatic graph with red loop edges (resp. green loop edges) glued at all nodes for type $B$ (resp. type $C$), and finally a complete bichromatic graph for type $BC$.
		See e.g. below for $n = 4$ nodes:
		
		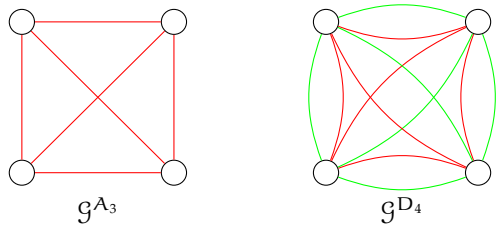
\begin{figure}[H]
			\centering
			\begin{tikzpicture}
				\node (a) at (0,0) {};
				\node (b) at (2,0) {};
				\node (c) at (2,2) {};
				\node (d) at (0,2) {};
				\node(e) [draw=none] at (1,-.5) {$\mc G^{A_3}$};
				\node (f) at (4,0) {};
				\node (g) at (6,0) {};
				\node (h) at (6,2) {};
				\node (i) at (4,2) {};
				\node (j) [draw=none] at (5,-.5) {$\mc G^{D_4}$};
				\graph{
					(a) --[red] (b) --[red] (c) --[red] (d) --[red] (a),
					(a) --[red] (c),
					(b) --[red] (d),
					(f) --[red,bend left=20] (g) -- [red,bend left=20] (h) --[red,bend left=20] (i) --[red,bend left=20] (f),
					(f) --[green,bend right=20] (g) -- [green,bend right=20] (h) --[green,bend right=20] (i) --[green,bend right=20] (f),
					(f) --[red,bend left=20] (h),
					(f) --[green,bend right=20] (h),
					(g) --[red,bend left=20] (i),
					(g) --[green,bend right=20] (i)
				};
			\end{tikzpicture}
			\caption{Examples of simply-laced classical graphs.}
		\end{figure}
		
		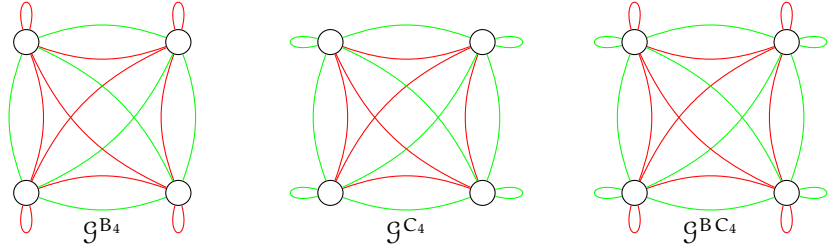
\begin{figure}[H]
			\centering
			\begin{tikzpicture}
				\node (a) at (0,0) {};
				\node (b) at (2,0) {};
				\node (c) at (2,2) {};
				\node (d) at (0,2) {};
				\node(e) [draw=none] at (1,-.5) {$\mc G^{B_4}$};
				\node (f) at (4,0) {};
				\node (g) at (6,0) {};
				\node (h) at (6,2) {};
				\node (i) at (4,2) {};
				\node (j) [draw=none] at (5,-.5) {$\mc G^{C_4}$};
				\node (k) at (8,0) {};
				\node (l) at (10,0) {};
				\node (m) at (10,2) {};
				\node (n) at (8,2) {};
				\node (o) [draw=none] at (9,-.5) {$\mc G^{BC_4}$};
				\graph{
					(a) --[red,bend left=20] (b) --[red,bend left=20] (c) --[red,bend left=20] (d) --[red,bend left=20] (a),
					(a) --[green,bend right=20] (b) -- [green,bend right=20] (c) --[green,bend right=20] (d) --[green,bend right=20] (a),
					(a) --[red,bend left=20] (c),
					(a) --[green,bend right=20] (c),
					(b) --[red,bend left=20] (d),
					(b) --[green,bend right=20] (d),
					(f) --[red,bend left=20] (g) -- [red,bend left=20] (h) --[red,bend left=20] (i) --[red,bend left=20] (f),
					(f) --[green,bend right=20] (g) -- [green,bend right=20] (h) --[green,bend right=20] (i) --[green,bend right=20] (f),
					(f) --[red,bend left=20] (h),
					(f) --[green,bend right=20] (h),
					(g) --[red,bend left=20] (i),
					(g) --[green,bend right=20] (i),
					
					(k) --[red,bend left=20] (l) -- [red,bend left=20] (m) --[red,bend left=20] (n) --[red,bend left=20] (k),
					(k) --[green,bend right=20] (l) -- [green,bend right=20] (m) --[green,bend right=20] (n) --[green,bend right=20] (k),
					(k) --[red,bend left=20] (m),
					(k) --[green,bend right=20] (m),
					(l) --[red,bend left=20] (n),
					(l) --[green,bend right=20] (n)
				};
				\path
					(a) edge [red,loop below] (a)
					(b) edge [red,loop below] (b)
					(c) edge [red,loop above] (c)
					(d) edge [red,loop above] (d)
					(f) edge [green,loop left] (f)
					(g) edge [green,loop right] (g)
					(h) edge [green,loop right] (h)
					(i) edge [green,loop left] (i)
					(k) edge [red,loop below] (k)
					(l) edge [red,loop below] (l)
					(m) edge [red,loop above] (m)
					(n) edge [red,loop above] (n)
					(k) edge [green,loop left] (k)
					(l) edge [green,loop right] (l)
					(m) edge [green,loop right] (m)
					(n) edge [green,loop left] (n);
			\end{tikzpicture}
			\caption{Examples of non-simply-laced classical graphs. \qedhere}
		\end{figure}
	\end{exem}
	
	\begin{rema}[Dual reading]
		In this construction $\Phi^{\mc G}$ is naturally a subset of $V = \mb C^n$, but one may equivalently work with subsets of $V^{\dual}$ using the dual basis $e_i^{\dual}$---defined by $\braket{ e_i^{\dual},e_j } = \delta_{ij}$ in the canonical pairing $\braket{ \cdot,\cdot } \cl V^{\dual} \otimes V \to \mb C$.
		In this case we consider the linear functionals 
		\begin{equation*}
			\alpha^{\pm}_{ij} = e^{\dual}_i \pm e^{\dual}_j \, , \quad \alpha_i = e^{\dual}_i \in V^{\dual} \, ,
		\end{equation*}
		and again get a bijection $\mc G \mapsto \Phi^{\mc G} \sse V^{\dual}$ between bichromatic graphs and symmetric subsets of the dual/inverse root system of $BC_n$: note $e_i^{\dual} = ( \cdot \mid \cdot )^{\flat}(e_i) \in V^{\dual}$ in the isomorphism $( \cdot \mid \cdot)^{\flat} \cl V \to V^{\dual}$ provided by the (canonical) invariant scalar product.
		The notation is then closer to the prototype in \S~\ref{sec:type_A}, identifying $V = \mb C^n$ with the Cartan subalgebra of $\mf{gl}_n(\mb C)$---so the root system lies in the dual.
	\end{rema}
	
	\section{Crystallographs and root subsystems}
	\label{sec:main_correspondence}
	
	\subsection{}
	
	The crucial fact is that one can explicitly describe the class of bichromatic graphs that map to root subsystems, rather than arbitrary symmetric subsets of roots.
	This is achieved by encoding the transformations of roots under mutual root-hyperplane reflections, in the following (main) definition.

	\begin{defi}
		\label{def:crystallographs}
		
		A bichromatic graph $(\mc G,c)$ is a \emph{crystallograph}\fn{
			As in ``\emph{crystallo}graphic (bichromatic) \emph{graph}''.}
		if the following conditions hold:
		\begin{enumerate}
			\item if a node is incident to two distinct straight edges, then $\mc G$ contains a straight edge closing the triangle; the third side is red if and only if the first two have the same colour---else it is green;
			
			\begin{figure}[H]
				\centering
				\begin{tikzpicture}
					\node (a) at (0,0) {};
					\node (b) at (2,0) {};
					\node (c) at (1,sqrt 3) {};
					\node (d) [draw=none] at (2,sqrt 3/2) {};
					\node (e) [draw=none] at (4,sqrt 3/2) {};
					\node (f) at (4,0) {};
					\node (g) at (6,0) {};
					\node (h) at (5,sqrt 3) {};
					\graph{
						(a) --[red] (b) --[red] (c),
						(d) ->[ultra thick] (e),
						(f) --[red] (g) --[red] (h) --[red] (f)
					};
				\end{tikzpicture}
				\caption{Straight edge closure condition for crystallographs (I).}
				\label{fig:straight_edge_closure_condition_crystallographs_I}
			\end{figure}
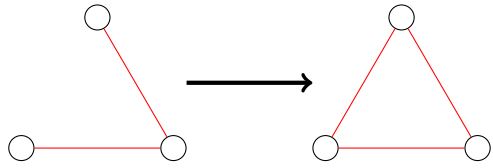
			
			\begin{figure}[H]
				\centering
				\begin{tikzpicture}
					\node (a) at (0,0) {};
					\node (b) at (2,0) {};
					\node (c) at (1,sqrt 3) {};
					\node (d) [draw=none] at (2,sqrt 3/2) {};
					\node (e) [draw=none] at (4,sqrt 3/2) {};
					\node (f) at (4,0) {};
					\node (g) at (6,0) {};
					\node (h) at (5,sqrt 3) {};
					\graph{
						(a) --[red] (b) --[green] (c),
						(d) ->[ultra thick] (e),
						(f) --[red] (g) --[green] (h) --[green] (f)
					};
				\end{tikzpicture}
				\caption{Straight edge closure condition for crystallographs (II).}
				\label{fig:straight_edge_closure_condition_crystallographs_II}
			\end{figure}
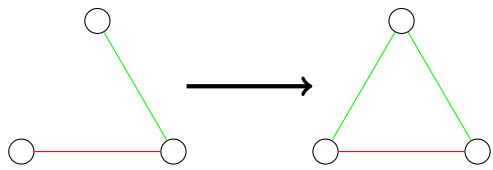
			
			\begin{figure}[H]
				\centering
				\begin{tikzpicture}
					\node (a) at (0,0) {};
					\node (b) at (2,0) {};
					\node (c) at (1,sqrt 3) {};
					\node (d) [draw=none] at (2,sqrt 3/2) {};
					\node (e) [draw=none] at (4,sqrt 3/2) {};
					\node (f) at (4,0) {};
					\node (g) at (6,0) {};
					\node (h) at (5,sqrt 3) {};
					\graph{
						(a) --[green] (b) --[green] (c),
						(d) ->[ultra thick] (e),
						(f) --[green] (g) --[green] (h) --[red] (f)
					};
				\end{tikzpicture}
				\caption{Straight edge closure condition for crystallographs (III).}
				\label{fig:straight_edge_closure_condition_crystallographs_III}
			\end{figure}
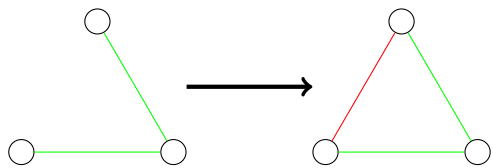
				
			\item if a node is incident to a loop edge and a straight edge, then $\mc G$ contains the straight edge of opposite colour, as well as the loop edge of the same colour at the opposite end.
			
			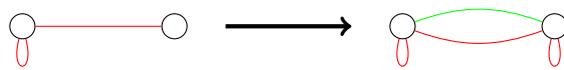
\begin{figure}[H]
				\centering
				\begin{tikzpicture}
					\node (a) at (0,0) {};
					\node (b) at (2,0) {};
					\node (c) [draw=none] at (2.5,0) {};
					\node (d) [draw=none] at (4.5,0) {};
					\node (e) at (5,0) {};
					\node (f) at (7,0) {};
					\graph{
						(a) --[red] (b),
						(c) ->[ultra thick] (d),
						(e) --[red,bend right=20] (f),
						(e) --[green,bend left=20] (f)
					};
					\path
						(a) edge [red,loop below] (a)
						(e) edge [red,loop below] (e)
						(f) edge [red,loop below] (f);
				\end{tikzpicture}
				\caption{Loop edge closure condition for crystallographs (I).} 
				\label{fig:red_loop_edge_closure_crystallograph_I}
			\end{figure}

			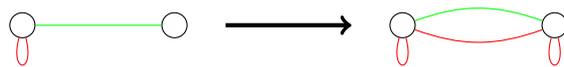
\begin{figure}[H]
				\centering
				\begin{tikzpicture}
					\node (a) at (0,0) {};
					\node (b) at (2,0) {};
					\node (c) [draw=none] at (2.5,0) {};
					\node (d) [draw=none] at (4.5,0) {};
					\node (e) at (5,0) {};
					\node (f) at (7,0) {};
					\graph{
						(a) --[green] (b),
						(c) ->[ultra thick] (d),
						(e) --[red,bend right=20] (f),
						(e) --[green,bend left=20] (f)
					};
					\path
						(a) edge [red,loop below] (a)
						(e) edge [red,loop below] (e)
						(f) edge [red,loop below] (f);
				\end{tikzpicture}
				\caption{Loop edge closure condition for crystallographs (II).} 
				\label{fig:red_loop_edge_closure_crystallograph_II}
			\end{figure}
			
			\begin{figure}[H]
				\centering
				\begin{tikzpicture}
					\node (a) at (0,0) {};
					\node (b) at (2,0) {};
					\node (c) [draw=none] at (2.5,0) {};
					\node (d) [draw=none] at (4.5,0) {};
					\node (e) at (5,0) {};
					\node (f) at (7,0) {};
					\graph{
						(a) --[red] (b),
						(c) ->[ultra thick] (d),
						(e) --[red,bend right=20] (f),
						(e) --[green,bend left=20] (f)
					};
					\path
						(a) edge [green,loop below] (a)
						(e) edge [green,loop below] (e)
						(f) edge [green,loop below] (f);
				\end{tikzpicture}
				\caption{Loop edge closure condition for crystallographs (III).} 
				\label{fig:green_loop_edge_closure_crystallograph_I}
			\end{figure}
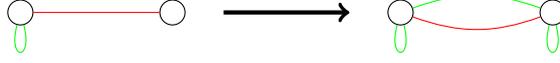
			
			\begin{figure}[H]
				\centering
				\begin{tikzpicture}
					\node (a) at (0,0) {};
					\node (b) at (2,0) {};
					\node (c) [draw=none] at (2.5,0) {};
					\node (d) [draw=none] at (4.5,0) {};
					\node (e) at (5,0) {};
					\node (f) at (7,0) {};
					\graph{
						(a) --[green] (b),
						(c) ->[ultra thick] (d),
						(e) --[red,bend right=20] (f),
						(e) --[green,bend left=20] (f)
					};
					\path
						(a) edge [green,loop below] (a)
						(e) edge [green,loop below] (e)
						(f) edge [green,loop below] (f);
				\end{tikzpicture}
				\caption{Loop edge closure condition for crystallographs (IV).} 
				\label{fig:green_loop_edge_closure_crystallograph_II}
			\end{figure}
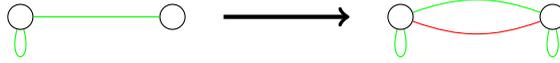
		\end{enumerate}
	\end{defi}
	
	\subsubsection{}
	
	The fact that the classical root systems of Ex.~\ref{ex:classical_crystallographs} correspond to crystallographs is an intended feature.
	
	\begin{theo}
		\label{thm:crystallographs_classify_root_systems}
		
		The association $(\mc G,c) \mapsto \Phi^{\mc G,c}$ restricts to a bijection between crystallographs on $n \geq 1$ nodes and root subsystems of $BC_n$.
	\end{theo}

	\begin{proof}
		The conditions in Def.~\ref{def:crystallographs} are the translation of the relevant identities for the root-hyperplane reflections.
		
		Namely denote
		\begin{equation*}
			\sigma^{\pm}_{ij} \ceqq \sigma_{e^{\pm}_{ij}} \, , \, \sigma_i \ceqq \sigma_{e_i} \in W \bigl( BC_n \bigr) \sse O \bigl( V, ( \cdot \mid \cdot ) \bigr) \, ,
		\end{equation*}
		for $i \neq j \in \ul n$, using the standard notation for Weyl group elements.
		Then the nontrivial identities are	
		\begin{equation}
			\label{eq:simply_laced_identities}
			\sigma^-_{ij} \bigl( e^-_{jk} \bigr) = e^-_{ik} \, , \quad \sigma^-_{ij} \bigl( e^+_{jk} \bigr) = e^+_{ik} \, , \quad \sigma^+_{ij} \bigl( e^-_{jk} \bigr) = -e^+_{ik} \, , \quad \sigma^+_{ij}\bigl( e^+_{jk} \bigr) = -e^-_{ik} \, ,
		\end{equation}
		and 
		\begin{equation}
			\label{eq:loop_identities}
			\sigma_i \bigl( e^{\pm}_{ij} \bigr) = -e^{\mp}_{ij} \, , \qquad \sigma^{\pm}_{ij}(e_j) = \mp e_i \, ,
		\end{equation}
		for distinct indices $i,j,k \in \ul n$.
		
		Thus indeed a bichromatic subgraph $(\mc G,c)$ is a crystallograph if and only if the associated (symmetric) subset of roots is closed under mutual root-hyperplane reflections: the row of identities~\eqref{eq:simply_laced_identities} corresponds to the former ``simply-laced'' condition of Def.~\ref{def:crystallographs} (with $i,j,k \in \mc G_0$ being the vertices of the triangle), and~\eqref{eq:loop_identities} to the condition involving loop edges (with $i,j \in \mc G_0$ being the ends of the straight edge).
	\end{proof}

	\begin{rema}
		\label{rk:graph_root_system_correspondence}
		
		By construction a crystallograph is connected if and only if the associated root system is irreducible; it is simply-laced if and only if the associated root system is a direct sum of simply-laced irreducible reduced root systems; and it has a double loop edge if and only if the associated root system is nonreduced.
		
		Moreover the bichromatic graph $(\mc G^{\Phi^{\dual}},c^{\dual})$ associated with the dual/inverse system is given by $\mc G^{\Phi^{\dual}} = \mc G^{\Phi}$, and then swapping the colour of all loop edges of $(\mc G^{\Phi},c)$ (so there is a natural inherited notion of dual/inverse crystallograph).
		
		Finally the \emph{rank} of a crystallograph can be defined by $\rk(\mc G,c) \ceqq rk(\Phi^{\mc G,c})$, and by construction $\abs{ \Phi^{\mc G}} = 2 \abs{\mc G_1}$---as every edge corresponds to a pair of opposite roots.
	\end{rema}

	\section{Crystallograph classifications}
	\label{sec:crystallograph_classification}
	
	\subsection{}
	
	The point of introducing crystallographs is that one can directly prove the following classification statement.

	\begin{theo}
		\label{thm:crystallograph_classification}
		
		Let $(\mc G,c)$ be a crystallograph on $n \geq 1$ nodes.
		Then $\mc G$ is a disjoint union of the following types of (crystallo)graphs:
		\begin{itemize}
			\item one of the ``classical'' graphs $\mc G^{A_{m-1}}$, $\mc G^{D_m}$, $\mc G^{B_m}$, $\mc G^{C_m}$, $\mc G^{BC_m}$, with $m \leq n$;
			
			\item a bichromatic graph $\mc G^{d_1,d_2}$ on a bipartite set of nodes $\mc G_0^{d_1,d_2} = I_1 \mathsmaller{\coprod} I_2$, with $\abs{I_i} = d_i$ and $d_1 + d_2 \leq n$, such that $\mc G^{d_1,d_2}_R$ is a disjoint union of two simply-laced complete graphs on each part, while $\mc G^{d_1,d_2}_G$ has one straight edge between every pair of nodes lying in two different parts---and none other.
		\end{itemize}
	\end{theo}
	
	\subsubsection{} 
	
	The green part of $\mc G^{d_1,d_2}$ is thus a (monochromatic, simply-laced) complete bipartite graph on the set of nodes $I_1 \mathsmaller{\coprod} I_2$; see an example below:
	
	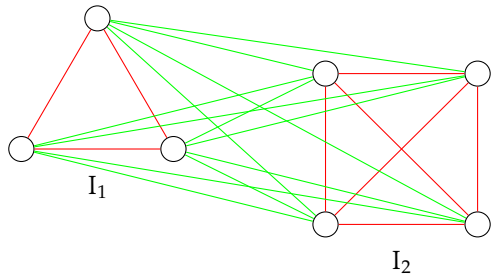
\begin{figure}[H]
		\centering
		\begin{tikzpicture}
			\node (a) at (0,1) {};
			\node (b) at (2,1) {};
			\node (c) at (1,sqrt 3+1) {};
			\node (d) [draw=none] at (1,.5) {$I_1$};
			\node (e) at (4,0) {};
			\node (f) at (6,0) {};
			\node (g) at (6,2) {};
			\node (h) at (4,2) {};
			\node (i) [draw=none] at (5,-.5) {$I_2$};
			\graph{
				(a) --[red] (b) --[red] (c) --[red] (a),
				(e) --[red] (f) --[red] (g) --[red] (h) --[red] (e),
				(e) --[red] (g),
				(f) --[red] (h),
				(a) --[green] (e),
				(a) --[green] (f),
				(a) --[green] (g),
				(a) --[green] (h),
				(b) --[green] (e),
				(b) --[green] (f),
				(b) --[green] (g),
				(b) --[green] (h),
				(c) --[green] (e),
				(c) --[green] (f),
				(c) --[green] (g),
				(c) --[green] (h)
			};
		\end{tikzpicture}
		\caption{Example: a graph of type $\mc G^{d_1,d_2}$ (for $(d_1,d_2) = (3,4)$).}
	\end{figure}

	\begin{proof}
		Let us first suppose $\mc G$ is simply-laced, which involves the first condition of Def.~\ref{def:crystallographs}.
		
		That condition implies that if two nodes are the endpoints of a path of red edges, then they are the endopoints of a red edge, upon completing a sequence of red triangles:
		
		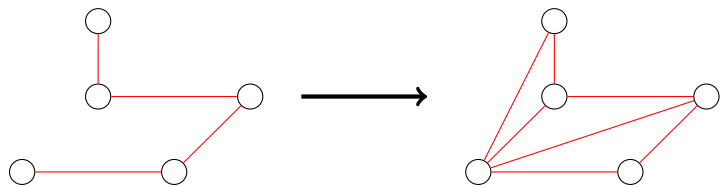
\begin{figure}[H]
			\centering
			\begin{tikzpicture}
				\node (a) at (0,0) {};
				\node (b) at (2,0) {};
				\node (c) at (3,1) {};
				\node (d) at (1,1) {};
				\node (e) at (1,2) {};
				\node (f) [draw=none] at (3.5,1) {};
				\node (g) [draw=none] at (5.5,1) {};
				\node (h) at (6,0) {};
				\node (i) at (8,0) {};
				\node (j) at (9,1) {};
				\node (k) at (7,1) {};
				\node (l) at (7,2) {};
				\graph{
					(a) --[red] (b) --[red] (c) --[red] (d) --[red] (e),
					(f) ->[ultra thick] (g),
					(h) --[red] (i) --[red] (j) --[red] (k) --[red] (l),
					(h) --[red] (j),
					(h) --[red] (k),
					(h) --[red] (l)
				};
			\end{tikzpicture}
			\caption{Example: completing red triangles along a red path.}
		\end{figure}
		
		Thus $\mc G_R$ splits into a disjoint union of simply-laced complete graphs.
		Further adding a green edge within such a component turns this latter into a simply-laced complete bichromatic graph, upon completion of some triangles---both red and green:
		
		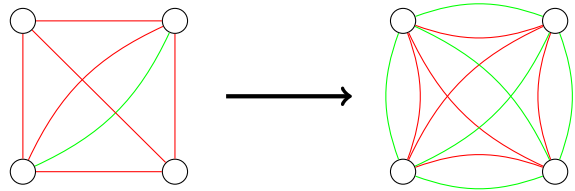
\begin{figure}[H]
			\centering
			\begin{tikzpicture}
				\node (a) at (0,0) {};
				\node (b) at (2,0) {};
				\node (c) at (2,2) {};
				\node (d) at (0,2) {};
				\node (e) [draw=none] at (2.5,1) {};
				\node (f) [draw=none] at (4.5,1) {};
				\node (g) at (5,0) {};
				\node (h) at (7,0) {};
				\node (i) at (7,2) {};
				\node (j) at (5,2) {};
				\graph{
					(a) --[red] (b) --[red] (c) --[red] (d) --[red] (a),
					(a) --[red,bend left = 20] (c),
					(a) --[green,bend right = 20] (c),
					(b) --[red] (d),
					(e) ->[ultra thick] (f),
					(g) --[red,bend left = 20] (h) --[red,bend left = 20] (i) --[red,bend left = 20] (j) --[red,bend left = 20] (g),
					(g) --[green,bend right = 20] (h) --[green,bend right = 20] (i) --[green,bend right = 20] (j) --[green,bend right = 20] (g),
					(g) --[red,bend left=20] (i),
					(h) --[red,bend left=20] (j),
					(g) --[green,bend right=20] (i),
					(h) --[green,bend right=20] (j),
				};
			\end{tikzpicture}
			\caption{Example: completing triangles within a red component.}
		\end{figure}
		
		Thus at this stage $\mc G$ splits as a disjoint union of type-$A$ (= red) and type-$D$ (= bichromatic) components, and we can in principle further add green edges between two such.		
		Reasoning as above shows that if there is a green edge between two nodes lying in different red components, then all pairs of nodes lying in the two different parts will be the endpoints of a green edge:
		
		\begin{figure}[H]
			\centering
			\begin{tikzpicture}[scale=1.5]
			\node (a) at (0,.5) {};
			\node (b) at (1,.5) {};
			\node (c) at (.5,sqrt 3/2+.5) {};
			\node (e) at (2,0) {};
			\node (f) at (3,0) {};
			\node (g) at (3,1) {};
			\node (h) at (2,1) {};
			\node (d) [draw=none]  at (3.5,.5) {};
			\node (i) [draw=none] at (4.5,.5) {};
			\node (j) at (5,.5) {};
			\node (k) at (6,.5) {};
			\node (l) at (5.5,sqrt 3/2+.5) {};
			\node (m) at (7,0) {};
			\node (n) at (8,0) {};
			\node (o) at (8,1) {};
			\node (p) at (7,1) {};
			\graph{
				(a) --[red] (b) --[red] (c) --[red] (a),
				(e) --[red] (f) --[red] (g) --[red] (h) --[red] (e),
				(e) --[red] (g),
				(f) --[red] (h),
				(b) --[green] (h),
				(d) ->[ultra thick] (i),
				(j) --[red] (k) --[red] (l) --[red] (j),
				(m) --[red] (n) --[red] (o) --[red] (p) --[red] (m),
				(m) --[red] (o),
				(n) --[red] (p),
				(j) --[green] (m),
				(j) --[green] (n),
				(j) --[green] (o),
				(j) --[green] (p),
				(k) --[green] (m),
				(k) --[green] (n),
				(k) --[green] (o),
				(k) --[green] (p),
				(l) --[green] (m),
				(l) --[green] (n),
				(l) --[green] (o),
				(l) --[green] (p)
			};
			\end{tikzpicture}
			\caption{Example: completing triangles between two red components.}
		\end{figure}
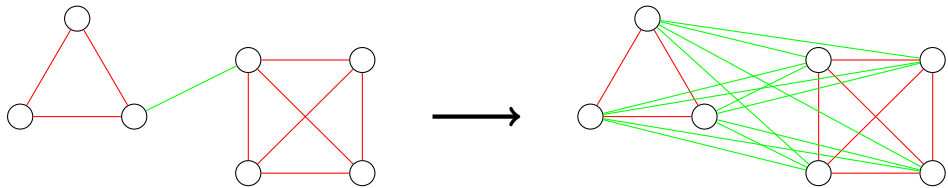
		
		Analogously if a green edge is drawn from a bichromatic component to either a red or a bichromatic one, then both will be ``fused'' into a single bichromatic component---as all green edges will propagate.
				
		Finally suppose there is a node of a red component, which is linked to two different red components via green edges.
		A final application of the ``simply-laced'' condition of Def.~\ref{def:crystallographs} shows the two latter red components are then fused into a single red one, creating a graph of type $\mc G^{d_1,d_2}$:
		
		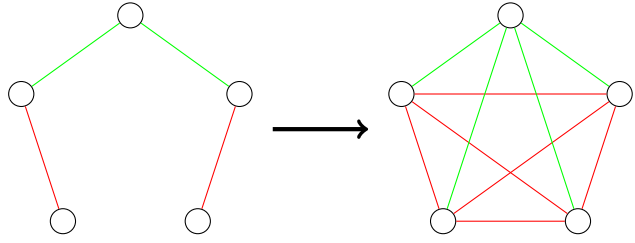
\begin{figure}[H]
			\centering
			\begin{tikzpicture}
				\node [draw=none,regular polygon, regular polygon sides=5, minimum size=3cm] (A) at (0,0) {};
				\foreach \i in {1,...,5}{
					\node (\i) at (A.corner \i) {};
				}
				\node (a) [draw=none] at (1.7,0) {};
				\node (b) [draw=none] at (3.3,0) {};
				\graph{
					(2) --[green] (1) --[green] (5), 
					(2) --[red] (3),
					(4) --[red] (5),
					(a) ->[ultra thick] (b),
				};
				\node [draw=none,regular polygon, regular polygon sides=5, minimum size=3cm] (B) at (5,0) {};
				\foreach \i in {1,...,5}{
					\node (\i) at (B.corner \i) {};
				}
				\graph{
					(2) --[green] (1) --[green] (5), 
					(2) --[red] (3) --[red] (4) --[red] (5) --[red] (2),
					(2) --[red] (4),
					(3) --[red] (5),
					(1) --[green] (3),
					(1) --[green] (4)
				};
				\end{tikzpicture}
			\caption{Example: completing triangles among triples of red components.}
		\end{figure}

		We conclude that $\mc G$ splits into a disjoint union of the simply-laced components of the statements.
		
		To conclude the proof for a generic graph, note gluing a loop edge to a node in any component of $\ol{\mc G}$ turns this latter into a simply-laced complete bichromatic graph, having a loop edge of the same colour at each node---by the second condition of Def.~\ref{def:crystallographs}.
	\end{proof}
	
	\begin{rema}
		\label{rk:red_components}
		
		In view of Thm.~\ref{thm:crystallograph_classification} we can talk of the ``red components'' of a crystallograph: these are the type-$A$ subgraphs that appear in his decomposition, i.e. the simply-laced complete red graphs it contains.
		By convention a trivial component, consisting of a disconnected node, is also considered to be red.\fn{
			It corresponds to the ``rank-zero'' root system $A_0 = \vn$.}
		
		These red components play a central role in the computation of root-kernel intersection, see \S~\ref{sec:kernels}; and the presence/absence of trivial components in turn controls quotients of root systems, see \S~\ref{sec:quotients}.
		
		By the same token, a ``bichromatic'' component is a connected component of any other type $B$, $C$, $D$ or $BC$.
	\end{rema}
	
	\subsection{}
	
	Now we consider the simply-laced bichromatic graphs $\mc G^{d_1,d_2}$, for $d_i \geq 1$; let simply $\Phi^{d_1,d_2} \ceqq \Phi^{\mc G^{d_1,d_2}} \sse BC_n$.

	\begin{lemm}
		\label{lem:bipartite_graph_is_type_A}
		
		The root system $\Phi^{d_1,d_2}$ is equivalent to $A_{d_1+d_2-1}$---inside $BC_n$.
	\end{lemm}

	\begin{proof}
		For $i \in I_1$ consider the element $\sigma_i = \sigma_{e_i} \in W \bigl( BC_{d_1+d_2} \bigr)$.
		This is an orthogonal transformation of $\mb{C}^{d_1+d_2}$, and it will turn $\Phi^{d_1,d_2}$ into a different root subsystem of $BC_{d_1+d_2}$ as follows (using~\eqref{eq:loop_identities}): the crystallograph associated with $\sigma_i(\Phi^{d_1,d_2})$ is obtained by changing the colour of all edges incident at the node $i \in I_1$.
		
		The result is a graph of type $\mc G^{d_1-1,d_2+1}$, since $i$ will now be part of a bigger red component on nodes $I_2 \mathsmaller{\coprod} \set{i}$.
		(There are $d_1 - 1$ red edges incident at $i \in I_1$, and $d_2$ incident green edges, so this operation is \emph{not} symmetric with respect to the two parts.)
		
		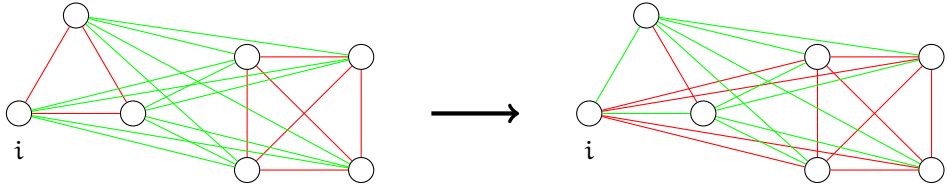
\begin{figure}[H]
			\centering
			\begin{tikzpicture}[scale=1.5]
			\node (a) [label=below:$i$] at (0,.5) {};
			\node (b) at (1,.5) {};
			\node (c) at (.5,sqrt 3/2+.5) {};
			\node (e) at (2,0) {};
			\node (f) at (3,0) {};
			\node (g) at (3,1) {};
			\node (h) at (2,1) {};
			\node (d) [draw=none]  at (3.5,.5) {};
			\node (i) [draw=none] at (4.5,.5) {};
			\node (j) [label=below:$i$] at (5,.5) {};
			\node (k) at (6,.5) {};
			\node (l) at (5.5,sqrt 3/2+.5) {};
			\node (m) at (7,0) {};
			\node (n) at (8,0) {};
			\node (o) at (8,1) {};
			\node (p) at (7,1) {};
			\graph{
				(a) --[red] (b) --[red] (c) --[red] (a),
				(e) --[red] (f) --[red] (g) --[red] (h) --[red] (e),
				(e) --[red] (g),
				(f) --[red] (h),
				(a) --[green] (e),
				(a) --[green] (f),
				(a) --[green] (g),
				(a) --[green] (h),
				(b) --[green] (e),
				(b) --[green] (f),
				(b) --[green] (g),
				(b) --[green] (h),
				(c) --[green] (e),
				(c) --[green] (f),
				(c) --[green] (g),
				(c) --[green] (h),
				(d) ->[ultra thick] (i),
				(j) --[green] (k) --[red] (l) --[green] (j),
				(m) --[red] (n) --[red] (o) --[red] (p) --[red] (m),
				(m) --[red] (o),
				(n) --[red] (p),
				(j) --[red] (m),
				(j) --[red] (n),
				(j) --[red] (o),
				(j) --[red] (p),
				(k) --[green] (m),
				(k) --[green] (n),
				(k) --[green] (o),
				(k) --[green] (p),
				(l) --[green] (m),
				(l) --[green] (n),
				(l) --[green] (o),
				(l) --[green] (p)
			};
			\end{tikzpicture}
			\caption{Example: Weyl reflection at a node $i$ of $\mc G^{3,4}$ (turning it into $\mc G^{2,5}$).}
		\end{figure}

		Then repeating at all nodes of $I_1 \sse \mc G^{d_1,d_2}_0$ yields a transformation
		\begin{equation}
			\label{eq:weyl_element}
			w = \prod_{I_1} \sigma_i \in W(BC_{d_1 + d_2}) \, ,            
		\end{equation}
		after choosing a total order for $I_1$.
		By construction the Weyl element~\eqref{eq:weyl_element} turns $\mc G^{d_1,d_2}$ into the simply-laced complete red graph on nodes $I_1 \mathsmaller{\coprod} I_2$, which is associated with the root subsystem $A_{d_1 + d_2 - 1}$ by the main correspondence.
	\end{proof}

	\begin{rema}
		One can also show (in nonconstructive fashion) that $\Phi^{d_1,d_2}$ is isomorphic to $A_{d_1 + d_2 - 1}$. 
		
		First, by Thm.~\ref{thm:crystallographs_classify_root_systems}, $\Phi^{d_1,d_2}$ is an irreducible reduced simply-laced root system, so it must be of type $A,D$ or $E$.
		Now the number of edges of $\mc G^{d_1,d_2}$ equals 
		\begin{equation*}
			\abs{\mc G^{d_1,d_2}_1} = \binom{d_1}{2} + \binom{d_2}{2} + d_1d_2 = \binom{d_1 + d_2}{2} \, .
		\end{equation*}
		Then one can show that $(\Phi^{d_1,d_2})^{\perp} = \mb{C}v$, as a subspace of 
		\begin{equation*}
			\mb{C}^{d_1 + d_2} \simeq \bigoplus_{i \in \mc G^{d_1,d_2}_0} \mb{C} e_i \sse \mb{C}^n \, ,
		\end{equation*}
		with the restriction of the standard scalar product of $V = \mb C^n$, where  
		\begin{equation*}
			v \ceqq \sum_{I_1} e_i - \sum_{I_2} e_i \in \mb{C}^{d_1+d_2} \, .
		\end{equation*}
		Hence $\Phi^{d_1,d_2}$ has rank equal to $\dim (\mb{C} \Phi^{d_1,d_2}) = d_1+d_2-1$, and it contains $2\abs{ \mc G^{d_1,d_2}_1} = (d_1+d_2)(d_1+d_2-1)$ roots: it must thus be isomorphic to $A_{d_1+d_2-1}$.
	\end{rema}

	\subsubsection{}
	
	Note the Weyl element~\eqref{eq:weyl_element} preserves the orthogonal complement of $\mb C^{d_1 + d_2} \sse \mb C^n$, so the other components of the graph are unaffected.	
	If one only considers root subsystems up to equivalence, it is then enough to work with disjoint unions of ``classical'' crystallographs: indeed the ``inner'' automorphisms suffice to turn any crystallograph into a union of those---by Thm.~\ref{thm:crystallograph_classification} and Lem.~\ref{lem:bipartite_graph_is_type_A}.
	
	\subsection{}
	
	Putting all together we have thus in particular proved the following statement.
	
	\begin{coro}
		\label{cor:root_system_classification}
		
		Suppose $\Phi' \sse \Phi \sse BC_n$ are nested root systems. 
		Then:
		\begin{itemize}
			\item if $\Phi = A_{n-1}$, $\Phi'$ is equivalent to a direct sum of type-$A$ irreducible root systems;
			
			\item if $\Phi = D_n$, $\Phi'$ is equivalent to a direct sum of type-$A$ and type-$D$ irreducible root systems;
			
			\item if $\Phi = B_n$, $\Phi'$ is a equivalent to a direct sum of type-$A$, type-$D$ and type-$B$ irreducible root systems;
			
			\item if $\Phi = C_n$, $\Phi'$ is equivalent to a direct sum of type-$A$, type-$D$ and type-$C$ irreducible root systems.
		\end{itemize}
	\end{coro}
	
% 	\subsubsection{}
% 	
% 	This is coherent with the tables of~\cite[\S~10]{oshima_2006_a_classification_of_subsystems_of_a_root_system}---with a different, more elementary proof.
		
	\section{Root-kernel intersections}
	\label{sec:kernels}

	\subsection{}
	
	Crystallographs can also be used to describe kernel intersections of arbitrary root subsystems $\Phi \sse BC_n$, generalising the prototype of \S~\ref{sec:type_A}.
	
	To this end let us work in the dual viewpoint where $\Phi^{\mc G,c} \sse V^{\dual}$ for any bichromatic graph $(\mc G,c)$.
	Let then $J$ be the set of red components of $\mc G$, as in Rk.~\ref{rk:red_components}, and denote $I_j \ceqq \mc G(j)_0 \sse \mc G_0$ the set of nodes of the component $\mc G(j) \sse \mc G$---for $j \in J$.
	Analogously, let $J'$ be the set of ``bipartite'' components of $\mc G$ of the form $\mc G^{d_1,d_2}$, and denote $I'_j = I'_{j,1} \coprod I'_{j,2} \sse \mc G_0$ the associated set of nodes---for $j \in J'$.
	(So we have tacitly fixed an order for the two parts.)
	
	\begin{prop}[Cf.~\cite{doucot_rembado_tamiozzo_2022_local_wild_mapping_class_groups_and_cabled_braids}]
		\label{prop:kernels}
		
		There is a vector space isomorphism $\mb{C}^{\abs{J}} \oplus \mb C^{\abs{J'}} \simeq \Ker \bigl( \Phi^{\mc G,c} \bigr)$, given by
		\begin{equation}
			\label{eq:kernel_basis_1}
			\ol{e}_j \longmapsto v_{I_j} = e_{I_j} \ceqq \sum_{i \in I_j} e_i \in V \, ,
		\end{equation}
		where $j \in J$, and $(\ol{e}_j)_j$ is the canonical basis of $\mb{C}^{\abs{J}}$; and 
		\begin{equation}
			\label{eq:kernel_basis_2}
			\ol{e}'_j \longmapsto v_{I'_j} = e_{I'_j} \ceqq \sum_{i \in I'_{j,1}} e_i - \sum_{i \in I'_{j,2}} e_i \in V,
		\end{equation} 
		where $j \in J'$, and $(\ol{e}'_j)_j$ is the canonical basis of $\mb C^{\abs{J'}}$.
	\end{prop}
	
	\begin{proof}
		The kernel is obtained by intersecting:
		\begin{itemize}
			\item $\Ker(\alpha^-_{ij})$, if $\set{i,j} \in \mc G_1$ is red;
			
			\item $\Ker(\alpha^+_{ij})$, if $\set{i,j} \in \mc G_1$ is green;
			
			\item $\Ker(\alpha_i)$, if $\set{i} \in \mc G_1$ is a loop edge.
		\end{itemize}
		
		Hence any loop annihilates the generator $e_i \in V$, and the same is true of $e_i,e_j \in V$ if $(\mc G,c)$ contains both a red and a green straight edge $\set{i,j}$.
		Thus only generators attached to nodes in red or ``bipartite'' components will not vanish in the kernel intersection.
		
		Conversely, if there is a red component on nodes $I_k \sse \mc G_0$, then (in the kernel) the coordinates along the vectors $e_i$ and $e_j$ are equal for all $i,j \in I_k$, which yields the line generated by the vector $e_{I_k}$ of~\eqref{eq:kernel_basis_1}.
		Further, if two red components on nodes $I_k$ and $I_l$ are linked by (all possible) green straight edges, then the corresponding generators sum to zero in the kernel intersection, which finally yields the line generated by the vector $e_{I'_{kl}}$ of~\eqref{eq:kernel_basis_2}, with $I'_{kl} = I_k \coprod I_l$.
	\end{proof}
	
	\section{Classification of quotients}
	\label{sec:quotients}
	
	\subsection{}
	
	Finally we can make further use of crystallographs to describe quotients of root systems, which is one of the main motivations.
	
	Consider thus the situation where $\mc G' \sse \mc G$ is an inclusion of crystallographs on $n \geq 1$ nodes.
	By Thm.~\ref{thm:crystallograph_classification} and Lem.~\ref{lem:bipartite_graph_is_type_A} we can assume both are disjoint union of graphs of type $A$, $B$, $C$, $D$ and $BC$, up to acting via the Weyl group.\fn{
		Hence we actually consider a subclass of ``dominant'' graphs, since we are interested in a classification up to equivalence.}

	By an ``edge between two components'' we mean a straight edge drawn between a pair of nodes lying in two distinct connected components; analogously a ``straight edge within a component'' is a straight edge between a pair of nodes lying within one component, and a ``loop edge in a component'' is a loop edge at any node within a component.
	
	\begin{defi}
		\label{def:quotient_graph}
		
		The \emph{quotient} graph $\mc G \bs \mc G'$ is the bichromatic graph whose nodes are given by the set $\Set{ \mc G'(j) }_{j \in J}$ of red components of $\mc G'$, and whose adjacency/colouring is defined as follows:
		\begin{itemize}
			\item if there is a straight edge $e \in \mc G_1$ between two red components $\mc G'(i)$ and $\mc G'(j)$ of $\mc G'$, put a straight edge $\Set{ \mc G'(i),\mc G'(j) }$ of the same colour;
			
			\item if there is a straight green edge $e \in \mc G_1$ within one red component $\mc G'(i)$ of $\mc G'$, put a green loop edge $\Set{ \mc G'(i) }$;
			
			\item if there is a straight edge $e \in \mc G_1$ from one red component $\mc G'(i)$ of $\mc G'$ to a bichromatic component of $\mc G$', put a red loop edge $\Set{ \mc G'(i) }$;
			
			\item if there is a loop edge $l \in \mc G_1$ in a red component $\mc G'(i)$ of $\mc G'$, put a loop edge $\Set{ \mc G'(i) }$ of the same colour.
		\end{itemize}
	\end{defi}
	
	\begin{exem}[Adjacency of a quotient of crystallographs]
	\label{ex:quotient_graph}
		
		In the following pictures, which exemplify all cases of Def.~\ref{def:quotient_graph}, we represent by dashed lines the edges lying in $\mc G_1 \sm \mc G'_1 \sse \mc G_1$.
		
		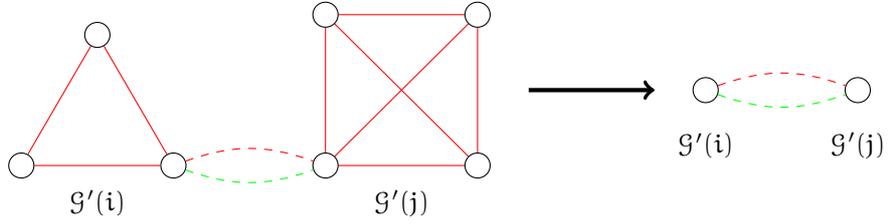
\begin{figure}[H]
			\centering
			\begin{tikzpicture}
				\node (a) at (0,0) {};
				\node (b) at (2,0) {};
				\node (c) at (1,sqrt 3) {};
				\node (d) [draw=none] at (1,-.5) {$\mc G'(i)$};
				\node (e) at (4,0) {};
				\node (f) at (6,0) {};
				\node (g) at (6,2) {};
				\node (h) at (4,2) {};
				\node (i) [draw=none] at (5,-.5) {$\mc G'(j)$};
				\node (j) [draw=none] at (6.5,1) {};
				\node (k) [draw=none] at (8.5,1) {};
				\node (l) [label=below:$\mc G'(i)$] at (9,1) {};
				\node (m) [label=below:$\mc G'(j)$] at (11,1) {};
				\graph{
					(a) --[red] (b) --[red] (c) --[red] (a),
					(e) --[red] (f) --[red] (g) --[red] (h) --[red] (e),
					(e) --[red] (g),
					(f) --[red] (h),
					(b) --[red,bend left=20, dashed] (e),
					(b) --[green,bend right=20, dashed] (e),
					(j) ->[ultra thick] (k),
					(l) --[red,bend left=20,dashed] (m),
					(l) --[green,bend right=20,dashed] (m),
				};
			\end{tikzpicture}
			\caption{Example: straight edge in $\mc G \bs \mc G'$, arising from a straight edge between two red components of $\mc G'$.}
		\end{figure}
		
		\begin{figure}[H]
				\centering
				\begin{tikzpicture}
					\node (a) at (0,0) {};
					\node (b) at (2,0) {};
					\node (c) at (2,2) {};
					\node (d) at (0,2) {};
					\node (e) [draw=none] at (2.5,1) {};
					\node (f) [draw=none] at (4.5,1) {};
					\node (g) [label=below:$\mc G'(i)$] at (5,1) {};
					\node (h) [draw=none] at (1,-.5) {$\mc G'(i)$};
					\graph{
						(a) --[red] (b) --[red] (c) --[red] (d) --[red] (a),
						(a) --[red,bend left = 20] (c),
						(a) --[green,bend right = 20,dashed] (c),
						(b) --[red] (d),
						(e) ->[ultra thick] (f),
					};
					\path
						(g) edge [loop above,green,dashed] (g);
				\end{tikzpicture}
				\caption{Example: green loop edge in $\mc G \bs \mc G'$, arising from a green straight edge within a red component of $\mc G'$.}
			\end{figure}
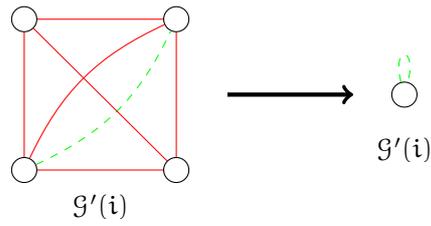
			
			\begin{figure}[H]
			\centering
			\begin{tikzpicture}
				\node (a) at (0,0) {};
				\node (b) at (2,0) {};
				\node (c) at (1,sqrt 3) {};
				\node (d) [draw=none] at (1,-.5) {$\mc G'(i)$};
				\node (e) at (4,0) {};
				\node (f) at (6,0) {};
				\node (g) at (6,2) {};
				\node (h) at (4,2) {};
				\node (i) [draw=none] at (5,-.5) {$\mc G'(j)$};
				\node (j) [draw=none] at (6.5,1) {};
				\node (k) [draw=none] at (8.5,1) {};
				\node (l) [label=below:$\mc G'(i)$] at (9,1) {};
				\graph{
					(a) --[red] (b) --[red] (c) --[red] (a),
					(e) --[red,bend left=20] (f) --[red,bend left=20] (g) --[red,bend left=20] (h) --[red,bend left=20] (e),
					(e) --[green,bend right=20] (f) --[green,bend right=20] (g) --[green,bend right=20] (h) --[green,bend right=20] (e),
					(e) --[red,bend left=20] (g),
					(f) --[red,bend left=20] (h),
					(e) --[green,bend right=20] (g),
					(f) --[green,bend right=20] (h),
					(b) --[red,bend left=20, dashed] (e),
					(b) --[green,bend right=20, dashed] (e),
					(j) ->[ultra thick] (k),
				};
				\path
					(l) edge [loop above,red,dashed] (l);
			\end{tikzpicture}
			\caption{Example: red loop edge in $\mc G \bs \mc G'$, arising from a straight edge from a red component to a bichromatic component of $\mc G'$.}
		\end{figure}
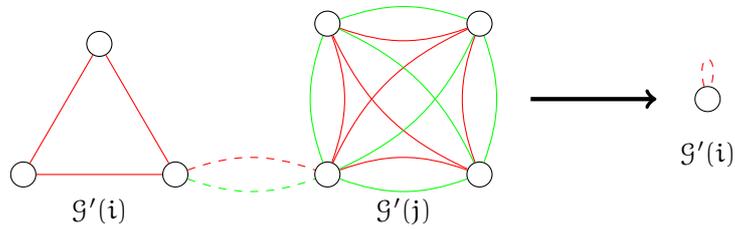

		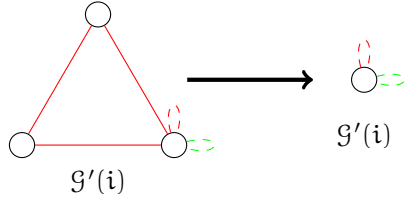
\begin{figure}[H]
			\centering
			\begin{tikzpicture}
				\node (a) at (0,0) {};
				\node (b) at (2,0) {};
				\node (c) at (1,sqrt 3) {};
				\node (d) [draw=none] at (1,-.5) {$\mc G'(i)$};
				\node (e) [draw=none] at (2,sqrt 3/2) {};
				\node (f) [draw=none] at (4,sqrt 3/2) {};
				\node (g) [label=below:$\mc G'(i)$] at (4.5,sqrt 3/2) {};
				\graph{
					(a) --[red] (b) --[red] (c) --[red] (a),
					(e) ->[ultra thick] (f)
				};
				\path
					(b) edge [loop above,red,dashed] (b)
					(b) edge [loop right,green,dashed] (b)
					(g) edge [loop above,red,dashed] (g)
					(g) edge [loop right,green,dashed] (g);
			\end{tikzpicture}
			\caption{Loop edge in $\mc G \bs \mc G'$, arising from a loop edge within a red component of $\mc G'$.}
		\end{figure}
	\end{exem}

	\subsection{}
	
	Now let $U \ceqq \Ker \bigl( \Phi^{\mc G'} \bigr) \sse V$, and consider the restriction of all other (co)roots $\alpha \in \Phi^{\mc G} \sm \Phi^{\mc G'}$ to $U$: this yields the set~\eqref{eq:restricted_system}.
	Clearly this set of restricted functionals is symmetric, so by the main correspondence it is associated with a bichromatic graph on as many nodes as $\dim(U) \leq n$; but by Prop.~\ref{prop:kernels} this dimension is the number of red components of $\mc G'$, which in turn is the set of nodes of the quotient graph in Def.~\ref{def:quotient_graph}: the next statement thus motivates that definition.

	\begin{theo}
		\label{thm:quotient_graphs_classify_restricted_arrangements}
		
		One has $\Phi^{\mc G \bs \mc G'} = \eval[1]{\Phi^{\mc G}}_U \sse U^{\vee}$.
	\end{theo}

	\begin{proof}
		This follows directly by evaluating all covectors $\alpha^{\pm}_{ij}, \alpha_i \in V^{\vee}$ on the vectors $e_{I_k} \in V$ in~\eqref{eq:kernel_basis_1}, which provide a basis of $U$, and then expressing the restrictions $\eval[1]{\alpha^{\pm}_{ij}}_U, \eval[1]{\alpha_i}_U \in U^{\vee}$ in terms of the dual basis $e_{I_k}^{\dual}$.
		
		(Recall we assume there are no bipartite components, so the vectors~\eqref{eq:kernel_basis_2} are no longer relevant.)		
	\end{proof}
	
	\subsubsection{}
	
	As a corollary we not only derive the type-$A$ case of \S~\ref{sec:type_A}, where the quotient root systems are still root systems, but also verify that this is false in general: the quotient graph is \emph{not} always a crystallograph, e.g. in the following situation---corresponding to the root-system inclusion $A_1 \sse D_4$:
	
	
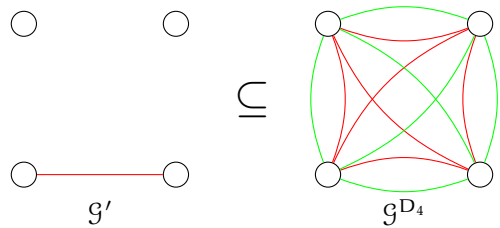
\begin{figure}[H]
		\centering
		\begin{tikzpicture}
			\node (a) at (0,0) {};
			\node (b) at (2,0) {};
			\node (c) at (2,2) {};
			\node (d) at (0,2) {};
			\node (e) [draw=none] at (1,-.5) {$\mc G'$};
			\node (f) [draw=none] at (3,1) {$\mathlarger{\mathlarger{\mathlarger{\mathlarger{\sse}}}}$};
			\node (g) at (4,0) {};
			\node (h) at (6,0) {};
			\node (i) at (6,2) {};
			\node (j) at (4,2) {};
			\node (k) [draw=none] at (5,-.5) {$\mc G^{D_4}$};
			\graph{
				(a) --[red] (b),
				(g) --[red,bend left=20] (h) --[red,bend left=20] (i) --[red,bend left=20] (j) --[red,bend left=20] (g),
				(g) --[red,bend left=20] (i),
				(h) --[red,bend left=20] (j),
				(g) --[green,bend right=20] (h) --[green,bend right=20] (i) --[green,bend right=20] (j) --[green,bend right=20] (g),
				(g) --[green,bend right=20] (i),
				(h) --[green,bend right=20] (j)
			};
		\end{tikzpicture}
	\end{figure}

	In this case the resulting arrangement consists of seven hyperplanes in $\mb C^3$, so it is \emph{not} crystallographic (cf. \S~\ref{sec:arrangements}).
	
	Nonetheless the algorithm of Def.~\ref{def:quotient_graph} yields the general way to encode all (sub)quotients of classical irreducible root systems, and moreover it is possible to explicitly describe the class of bichromatic graphs one obtains.
	
	\subsection{Crystallographic obstruction}
	\label{sec:quotient_crystallographic_obstruction}
	
	Here we study the obstruction for a quotient graph $\mc G \bs \mc G'$ to be a crystallograph.
	The result is that \emph{all but one} condition of Def.~\ref{def:crystallographs} apply.
	
	\begin{defi}
		\label{def:quasi_crystallographs}
		
		A \emph{quasi-crystallograph} is a bichromatic graph $(\mc G,c)$ satisfying all conditions of Def.~\ref{def:crystallographs}, with the following exception: if a straight edge and a green loop edge are incident at a common node, then the green loop edge need \emph{not} propagate at the other end of the straight edge---but the straight edge will double.
	\end{defi}
	
	Compare the next two figures with Figg.~\ref{fig:green_loop_edge_closure_crystallograph_I} and~\ref{fig:green_loop_edge_closure_crystallograph_II}.
	
	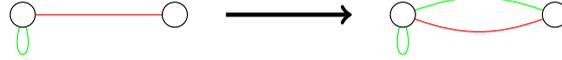
\begin{figure}[H]
		\centering
		\begin{tikzpicture}
			\node (a) at (0,0) {};
			\node (b) at (2,0) {};
			\node (c) [draw=none] at (2.5,0) {};
			\node (d) [draw=none] at (4.5,0) {};
			\node (e) at (5,0) {};
			\node (f) at (7,0) {};
			\graph{
				(a) --[red] (b),
				(c) ->[ultra thick] (d),
				(e) --[red,bend right=20] (f),
				(e) --[green,bend left=20] (f)
			};
			\path
				(a) edge [green,loop below] (a)
				(e) edge [green,loop below] (e);
		\end{tikzpicture}
		\caption{Green loop edge closure condition for quasi-crystallographs (I).} 
	\end{figure}
	
	\begin{figure}[H]
		\centering
		\begin{tikzpicture}
			\node (a) at (0,0) {};
			\node (b) at (2,0) {};
			\node (c) [draw=none] at (2.5,0) {};
			\node (d) [draw=none] at (4.5,0) {};
			\node (e) at (5,0) {};
			\node (f) at (7,0) {};
			\graph{
				(a) --[green] (b),
				(c) ->[ultra thick] (d),
				(e) --[red,bend right=20] (f),
				(e) --[green,bend left=20] (f)
			};
			\path
				(a) edge [green,loop below] (a)
				(e) edge [green,loop below] (e);
		\end{tikzpicture}
		\caption{Green loop edge closure condition for quasi-crystallographs (II).} 
	\end{figure}
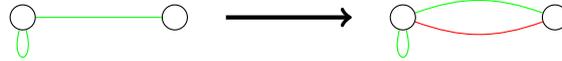
	
	\begin{theo}
		\label{thm:quasi_crystallographs_classify_quotients}
		
		All quotient graphs are quasi-crystallographs.
	\end{theo}
	
	\subsubsection{}
	
	We will prove this statement with a case-by-case analysis of the crystallograph conditions.
	In all the following local pictures we depict edges lying in $\mc G_1 \sm \mc G'_1$ with dashed lines---as above.

	\begin{prop}
		The simply-laced part of any quotient graph is a crystallograph.	
	\end{prop}
	
	\subsubsection{}
	
	Hence by Thm.~\ref{thm:crystallograph_classification} it will be a disjoint union of graphs of type $A$ and $D$.

	\begin{proof}
		We must prove the first condition of Def.~\ref{def:crystallographs} is satisfied for any pair of crystallographs $\mc G' \sse \mc G$, i.e. that suitable triangles close in the quotient.
		
		Suppose then we have two consecutive straight red edges $\Set{\mc G'(i),\mc G'(j)}$ and $\Set{ \mc G'(j),\mc G'(k) }$ in $\mc G \bs \mc G'$, for $i,j,k \in J$.
		This means there are red straight edges $\set{i_0,j_0}, \set{\wt j_0,k_0} \in \mc G_1 \sm \mc G'_1$ (in $\mc G$) such that $i_0 \in \mc G'(i)_0$, $j_0,\wt j_0 \in \mc G'(j)_0$ and $k_0 \in \mc G'(k)_0$:
		
		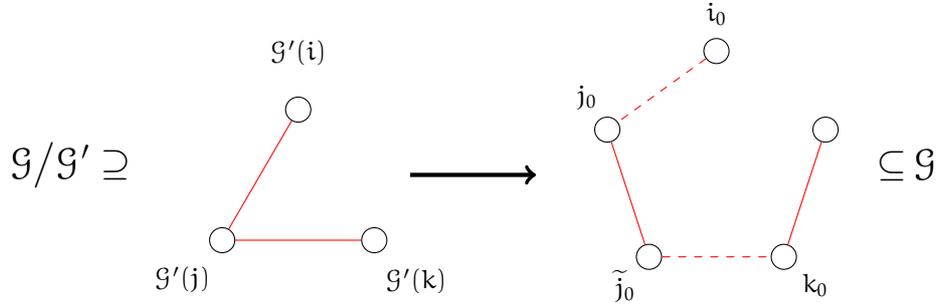
\begin{figure}[H]
			\centering
			\begin{tikzpicture}
				\node (a) [label=below left:$\mc G'(j)$] at (0,0) {};
				\node (b) [label=below right:$\mc G'(k)$] at (2,0) {};
				\node (c) [label=above:$\mc G'(i)$] at (1,sqrt 3) {};
				\node (d) [draw=none] at (2.3,sqrt 3/2) {};
				\node (e) [draw=none] at (4.3,sqrt 3/2) {};
				\node [draw=none,regular polygon, regular polygon sides=5, minimum size=3cm] (A) at (6.5,1) {};
				\foreach \i in {1,...,5}{
					\node (\i) at (A.corner \i) {};
				};
				\node (f) [draw=none] at (6.5,3) {$i_0$};
				\node (g) [draw=none] at (4.8,1.9) {$j_0$};
				\node (h) [draw=none] at (5.3,-.6) {$\wt j_0$};
				\node (i) [draw=none] at (7.8,-.6) {$k_0$};
				\node (j) [draw=none] at (-2,1) {$\mathlarger{\mathlarger{\mathlarger{\mc G \bs \mc G' \supseteq}}}$};
				\node (k) [draw=none] at (9,1) {$\mathlarger{\mathlarger{\mathlarger{\sse \mc G}}}$}; 
				\graph{
					(c) --[red] (a) --[red] (b),
					(d) ->[ultra thick] (e),
					(1) --[red,dashed] (2) --[red] (3) --[red,dashed] (4) --[red] (5),
				};
			\end{tikzpicture}
			\caption{Example: completing a red triangle in the quotient (I).}
		\end{figure}
		
		But since $\mc G$ is a crystallograph these three red components of $\mc G'$ are then contained in a single component of $\mc G$:
		
		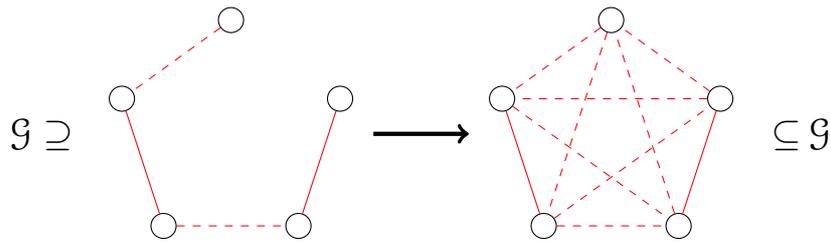
\begin{figure}[H]
			\centering
			\begin{tikzpicture}
				\node [draw=none,regular polygon, regular polygon sides=5, minimum size=3cm] (A) at (0,0) {};
				\foreach \i in {1,...,5}{
					\node (\i) at (A.corner \i) {};
				}
				\node (a) [draw=none] at (1.7,0) {};
				\node (b) [draw=none] at (3.3,0) {};
				\graph{
					(1) --[red,dashed] (2) --[red] (3) --[red,dashed] (4) --[red] (5),
					(a) ->[ultra thick] (b),
				};
				\node [draw=none,regular polygon, regular polygon sides=5, minimum size=3cm] (B) at (5,0) {};
				\foreach \i in {1,...,5}{
					\node (\i) at (B.corner \i) {};
				}
				\graph{
					(1) --[red,dashed] (2) --[red] (3) --[red,dashed] (4) --[red] (5) --[red,dashed] (1),
					(1) --[red,dashed] (3) --[red,dashed] (5) --[red,dashed] (2) --[red,dashed] (4) --[red,dashed] (1)
				};
				\node (j) [draw=none] at (-2.5,0) {$\mathlarger{\mathlarger{\mathlarger{\mc G \supseteq}}}$};
				\node (k) [draw=none] at (7.5,0) {$\mathlarger{\mathlarger{\mathlarger{\sse \mc G}}}$};
				\end{tikzpicture}
			\caption{Example: completing a red triangle in the quotient (II).}
		\end{figure}
		
		In particular $\set{i_0,k_0} \in \mc G_1 \sm \mc G'_1$ and the red triangle closes in the quotient:
		
		\begin{figure}[H]
			\centering
			\begin{tikzpicture}
				\node [draw=none,regular polygon, regular polygon sides=5, minimum size=3cm] (A) at (0,1) {};
				\foreach \i in {1,...,5}{
					\node (\i) at (A.corner \i) {};
				}
				\node (a) [draw=none] at (1.7,1) {};
				\node (b) [draw=none] at (3.3,1) {};
				\graph{
					(1) --[red,dashed] (2) --[red] (3) --[red,dashed] (4) --[red] (5) --[red,dashed] (1),
					(1) --[red,dashed] (3) --[red,dashed] (5) --[red,dashed] (2) --[red,dashed] (4) --[red,dashed] (1),
					(a) ->[ultra thick] (b)
				};
				\node (a) [label=below left:$\mc G'(j)$] at (4,0) {};
				\node (b) [label=below right:$\mc G'(k)$] at (6,0) {};
				\node (c) [label=above:$\mc G'(i)$] at (5,sqrt 3) {};
				\graph{
					(a) --[red] (b) --[red] (c) --[red](a)
				};
				\node (j) [draw=none] at (-2.5,1) {$\mathlarger{\mathlarger{\mathlarger{\mc G \supseteq}}}$};
				\node (k) [draw=none] at (7.5,1) {$\mathlarger{\mathlarger{\mathlarger{\sse \mc G \bs \mc G'}}}$};
				\end{tikzpicture}
			\caption{Example: completing a red triangle in the quotient (III).}
		\end{figure}
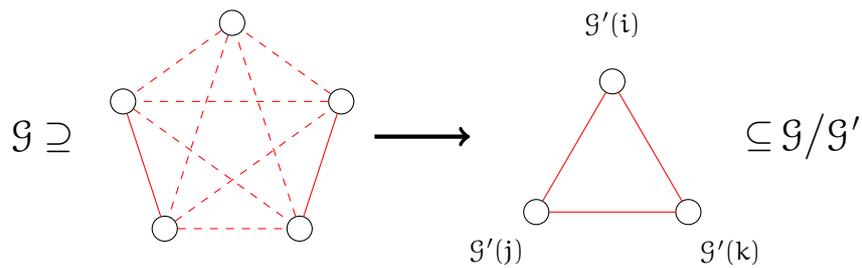
		
		Analogous arguments work in the cases where there are two consecutive straight edges of green/different colour.
	\end{proof}

	\subsubsection{}
	
	This is coherent with the fact that all quotients of type-$A$ crystallographs are crystallographs (of type-$A$)---and strengthens it.

	\begin{lemm}
		\label{lem:red_loop_closure_quotient}
		
		The second condition of Def.~\ref{def:crystallographs} is verified for red loop edges in $\mc G \bs \mc G'$.
	\end{lemm}

	\subsubsection{}
	
	In particular if the quotient has no green loop edges then it is a crystallograph: by Thm.~\ref{thm:crystallograph_classification} it will be a disjoint union of graphs of type $A,D$ and $B$.
	
	Hereafter a yellow straight edge (in a picture) stands for a straight edge of either colour.\fn{
		Recall $Y = R+G$ in the additive RGB colour model.}
	
	\begin{proof}
		Suppose there is a red loop edge $\Set{ \mc G'(i) }$ at a node of the quotient, and a straight edge $\Set{ \mc G'(i),\mc G'(j)}$ incident at the same node, for $i,j \in J$.
		Then we have a straight edge $\set{i_0,j_0} \in \mc G_1 \sm \mc G'_1$ (in $\mc G$) with $i_0 \in \mc G'(i)_0$ and $j_0 \in \mc G'(j)_0$, and one of the following two cases:
		\begin{enumerate}
			\item the component $\mc G'(i)$ has red loop edges in $\mc G$;
			
		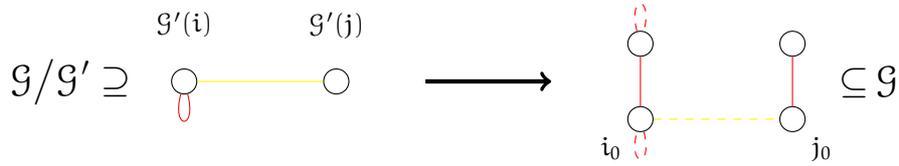
\begin{figure}[H]
			\centering 
			\begin{tikzpicture}
				\node (a) [label=above:$\mc G'(i)$] at (0,0) {};
				\node (b) [label=above:$\mc G'(j)$] at (2,0) {};
				\node (c) [draw=none] at (3,0) {};
				\node (d) [draw=none] at (5,0) {};
				\node (e) [label=below left:$i_0$] at (6,-.5) {};
				\node (f) at (6,.5) {};
				\node (g) [label=below right:$j_0$] at (8,-.5) {};
				\node (h) at (8,.5) {};
				\node(i) [draw=none] at (-1.5,0) {$\mathlarger{\mathlarger{\mathlarger{\mc G \bs \mc G' \supseteq}}}$};
				\node (j) [draw=none] at (9,0) {$\mathlarger{\mathlarger{\mathlarger{\sse \mc G}}}$};
				\graph{
					(a) --[yellow] (b),
					(c) ->[ultra thick] (d),
					(f) --[red] (e) --[yellow,dashed] (g) --[red] (h)
				};
				\path
					(a) edge [loop below,red] (a)
					(e) edge [loop below,red,dashed] (e)
					(f) edge [loop above,red,dashed] (f);
			\end{tikzpicture}
			\caption{Example: red loop edge closure condition for quotients, first case (I).}
			\label{fig:red_loop_closure_quotient_I}
		\end{figure}
			
		\item there is a straight edge $\set{\wt i_0,k_0} \in \mc G_1 \sm \mc G'_1$ (in $\mc G$) such that $k_0 \in \mc G_0$ lies within a bichromatic component of $\mc G'$.
		\end{enumerate}
		
		\begin{figure}[H]
			\centering 
			\begin{tikzpicture}
				\node (a) [label=above:$\mc G'(i)$] at (0,0) {};
				\node (b) [label=above:$\mc G'(j)$] at (2,0) {};
				\node (c) [draw=none] at (3,0) {};
				\node (d) [draw=none] at (5,0) {};
				\node [draw=none,regular polygon, regular polygon sides=5, minimum size=3cm] (B) at (7,0) {};
				\foreach \i in {1,...,5}{
					\node (\i) at (B.corner \i) {};
				}
				\graph{
					(1) --[yellow,dashed] (2) --[red] (3) --[yellow,dashed] (4) --[red,bend left =20] (5),
					(4) --[green,bend right=20] (5)
				};
				\node (e) [draw=none] at (7,2) {$j_0$};
				\node (f) [draw=none] at (5,.7) {$i_0$};
				\node (g) [draw=none] at (5.5,-1.2) {$\wt i_0$};
				\node (h) [draw=none] at (8.5,-1.2) {$k_0$};
				\node (i) [draw=none] at (-1.5,0) {$\mathlarger{\mathlarger{\mathlarger{\mc G \bs \mc G' \supseteq}}}$};
				\node (j) [draw=none] at (9.5,0) {$\mathlarger{\mathlarger{\mathlarger{\sse \mc G}}}$};
				\graph{
					(a) --[yellow] (b),
					(c) ->[ultra thick] (d),
				};
				\path
					(a) edge [loop below,red] (a);
			\end{tikzpicture}
			\caption{Example: red loop edge closure condition for quotients, second case (I).}
		\end{figure}
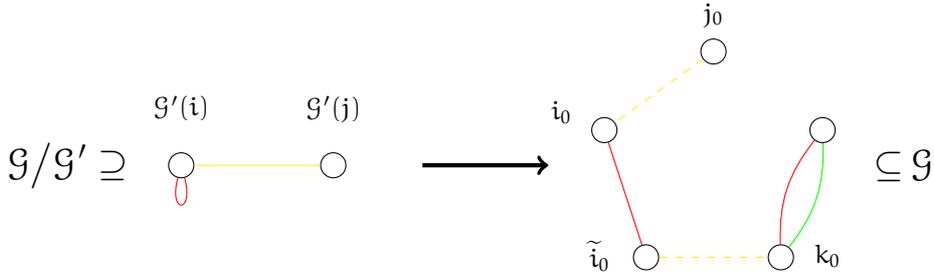
		
		In the first case the red loop edges propagate to the nodes of $\mc G'(j) \sse \mc G'$, since $\mc G$ is a crystallograph, and moreover the straight edge $\set{i_0,j_0}$ will ``double'' (in $\mc G$):
		
		\begin{figure}[H]
			\centering 
			\begin{tikzpicture}
				
				\node (a) at (0,-.5) {};
				\node (b) at (0,.5) {};
				\node(k) at (2,-.5) {};
				\node(l) at (2,.5) {};
				\node (c) [draw=none] at (3,0) {};
				\node (d) [draw=none] at (5,0) {};
				\node (e) at (6,-.5) {};
				\node (f) at (6,.5) {};
				\node (g) at (8,-.5) {};
				\node (h) at (8,.5) {};
				\node(i) [draw=none] at (-1,0) {$\mathlarger{\mathlarger{\mathlarger{\mc G \supseteq}}}$};
				\node (j) [draw=none] at (9,0) {$\mathlarger{\mathlarger{\mathlarger{\sse \mc G}}}$};
				\graph{
					(b) --[red] (a) --[yellow,dashed] (k) --[red] (l);
					(c) ->[ultra thick] (d),
					(f) --[red] (e) --[red,bend left=20,dashed] (g) --[red] (h),
					(e) --[green,bend right=20,dashed] (g)
				};	
				\path
					(a) edge [loop below,red,dashed] (a)
					(b) edge [loop above,red,dashed] (b)
					(e) edge [loop below,red,dashed] (e)
					(f) edge [loop above,red,dashed] (f)
					(g) edge [loop below,red,dashed] (g)
					(h) edge [loop above,red,dashed] (h);
			\end{tikzpicture}
			\caption{Example: red loop edge closure condition for quotients, first case (II).}
			\label{fig:red_loop_closure_quotient_II}
		\end{figure}
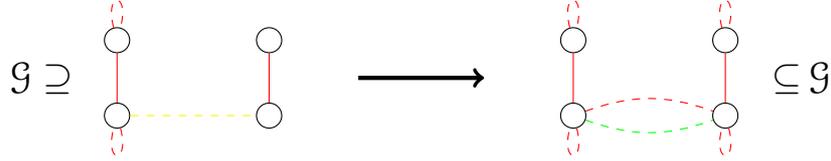
		
		Hence the same will happen in the quotient:
		
		\begin{figure}[H]
			\centering 
			\begin{tikzpicture}
				\node (a) at (0,-.5) {};
				\node (b) at (0,.5) {};
				\node(k) at (2,-.5) {};
				\node(l) at (2,.5) {};
				\node (c) [draw=none] at (3,0) {};
				\node (d) [draw=none] at (5,0) {};
				\node (e) [label=above:$\mc G'(i)$] at (6,0) {};
				\node (f) [label=above:$\mc G'(j)$] at (8,0) {};
				\node(i) [draw=none] at (-1,0) {$\mathlarger{\mathlarger{\mathlarger{\mc G \supseteq}}}$};
				\node (j) [draw=none] at (9.5,0) {$\mathlarger{\mathlarger{\mathlarger{\sse \mc G \bs \mc G'}}}$};
				\graph{
					(b) --[red] (a) --[red,bend left=20,dashed] (k) --[red] (l),
					(a) --[green,bend right=20,dashed] (k),
					(c) ->[ultra thick] (d),
					(e) --[red,bend left=20] (f),
					(e) --[green,bend right=20] (f)
				};	
				\path
					(a) edge [loop below,red,dashed] (a)
					(b) edge [loop above,red,dashed] (b)
					(k) edge [loop below,red,dashed] (k)
					(l) edge [loop above,red,dashed] (l)
					(e) edge [loop below,red] (e)
					(f) edge [loop below,red] (f);
			\end{tikzpicture}
			\caption{Example: red loop edge closure condition for quotients, first case (III).}
			\label{fig:red_loop_closure_quotient_III}
		\end{figure}
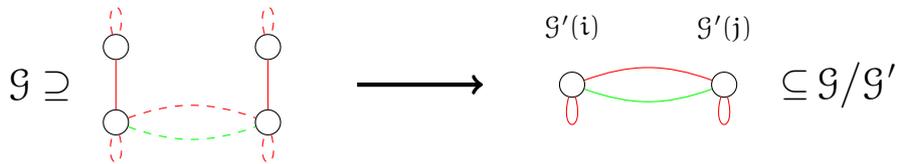
		
		In the second case the red component $\mc G'(i) \sse \mc G'$ is contained within a bichromatic component of $\mc G$ (since it is linked to a bichromatic component of $\mc G'$ by a straight edge); then analogously the subgraph $\mc G'(j) \sse \mc G'$ is also contained in that bichromatic component of $\mc G$:
		
		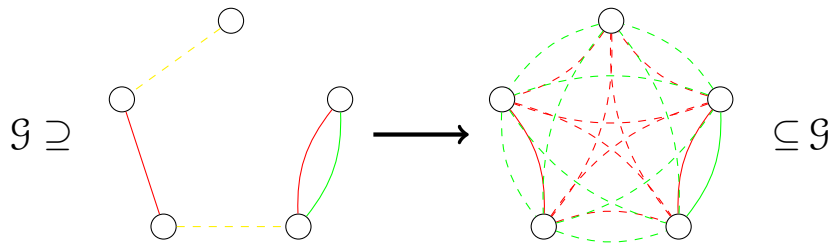
\begin{figure}[H]
			\centering
			\begin{tikzpicture}
				\node [draw=none,regular polygon, regular polygon sides=5, minimum size=3cm] (A) at (0,0) {};
				\foreach \i in {1,...,5}{
					\node (\i) at (A.corner \i) {};
				}
				\node (a) [draw=none] at (1.7,0) {};
				\node (b) [draw=none] at (3.3,0) {};
				\graph{
					(1) --[yellow,dashed] (2) --[red] (3) --[yellow,dashed] (4) --[red,bend left =20] (5),
					(4) --[green,bend right=20] (5),
					(a) ->[ultra thick] (b)
				};
				\node [draw=none,regular polygon, regular polygon sides=5, minimum size=3cm] (B) at (5,0) {};
				\foreach \i in {1,...,5}{
					\node (\i) at (B.corner \i) {};
				}
				\graph{
					(1) --[red,bend left=20,dashed] (2) --[red,bend left=20] (3) --[red,bend left=20,dashed] (4) --[red,bend left=20] (5) --[red,bend left=20,dashed] (1),
					(1) --[red,bend left=20,dashed] (3) --[red,bend left=20,dashed] (5) --[red,bend left=20,dashed] (2) --[red,bend left=20,dashed] (4) --[red,bend left=20,dashed] (1),
					(1) --[green,bend right=20,dashed] (2) --[green,bend right=20,dashed] (3) --[green,bend right=20,dashed] (4) --[green,bend right=20] (5) --[green,bend right=20,dashed] (1),
					(1) --[green,bend right=20,dashed] (3) --[green,bend right=20,dashed] (5) --[green,bend right=20,dashed] (2) --[green,bend right=20,dashed] (4) --[green,bend right=20,dashed] (1)
				};
				\node (j) [draw=none] at (-2.5,0) {$\mathlarger{\mathlarger{\mathlarger{\mc G \supseteq}}}$};
				\node (k) [draw=none] at (7.5,0) {$\mathlarger{\mathlarger{\mathlarger{\sse \mc G}}}$};
				\end{tikzpicture}
			\caption{Example: red loop edges closure condition for quotients, second case (II).}
			\label{fig:second_red_loop_closure_quotient_I}
		\end{figure}
		
		In particular there is a straight edge $\set{j_0,k_0} \in \mc G_1 \sm \mc G'_1$ that yields a red loop $\Set{ \mc G'(j) }$ in the quotient, as well as a straight edge of opposite colour $\set{i_0,j_0}$, which lead to the required crystallographic configuration in the quotient:
			
		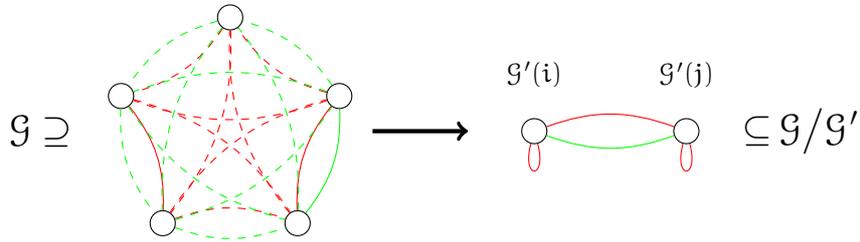
\begin{figure}[H]
			\centering
			\begin{tikzpicture}
				\node [draw=none,regular polygon, regular polygon sides=5, minimum size=3cm] (B) at (0,0) {};
				\foreach \i in {1,...,5}{
					\node (\i) at (B.corner \i) {};
				}
				\node (a) [draw=none] at (1.7,0) {};
				\node (b) [draw=none] at (3.3,0) {};
				\node (e) [label=above:$\mc G'(i)$] at (4,0) {};
				\node (f) [label=above:$\mc G'(j)$] at (6,0) {};
				\node (j) [draw=none] at (-2.5,0) {$\mathlarger{\mathlarger{\mathlarger{\mc G \supseteq}}}$};
				\node (k) [draw=none] at (7.5,0) {$\mathlarger{\mathlarger{\mathlarger{\sse \mc G \bs \mc G'}}}$};
				\graph{
					(1) --[red,bend left=20,dashed] (2) --[red,bend left=20] (3) --[red,bend left=20,dashed] (4) --[red,bend left=20] (5) --[red,bend left=20,dashed] (1),
					(1) --[red,bend left=20,dashed] (3) --[red,bend left=20,dashed] (5) --[red,bend left=20,dashed] (2) --[red,bend left=20,dashed] (4) --[red,bend left=20,dashed] (1),
					(1) --[green,bend right=20,dashed] (2) --[green,bend right=20,dashed] (3) --[green,bend right=20,dashed] (4) --[green,bend right=20] (5) --[green,bend right=20,dashed] (1),
					(1) --[green,bend right=20,dashed] (3) --[green,bend right=20,dashed] (5) --[green,bend right=20,dashed] (2) --[green,bend right=20,dashed] (4) --[green,bend right=20,dashed] (1),
					(a) ->[ultra thick] (b),
					(e) --[red,bend left=20] (f),
					(e) --[green,bend right=20] (f)
				};
				\path
					(e) edge [loop below,red] (e)
					(f) edge [loop below,red] (f);
				\end{tikzpicture}
			\caption{Example: red loop edges closure condition for quotients, second case (III). \qedhere}
			\label{fig:second_red_loop_closure_quotient_II}
		\end{figure}
	\end{proof}
	
	\subsubsection{}
	
	Let us finally consider \emph{green} loop edges: suppose there is a green loop edge at a node of the quotient $\mc G \bs \mc G'$, and a straight edge incident to it.

	\begin{lemm}
		If a green loop edge $\Set{ \mc G'(i) }$ in the quotient originates from green loop edges at the nodes of the red component $\mc G'(i) \sse \mc G'$, then the second condition of Def.~\ref{def:crystallographs} is verified for $\mc G \bs \mc G'$.
	\end{lemm}

	\begin{proof}
		This follows as in the proof of the first case in the previous Lem.~\ref{lem:red_loop_closure_quotient}; compare the following sequence of pictures with Figg.~\ref{fig:red_loop_closure_quotient_I},~\ref{fig:red_loop_closure_quotient_II} and~\ref{fig:red_loop_closure_quotient_III}:
		
		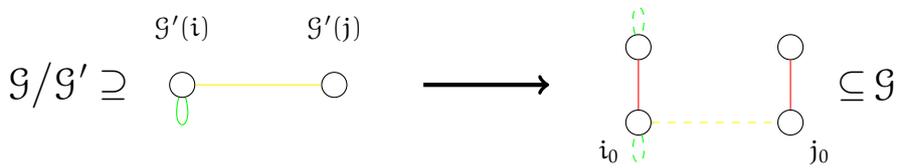
\begin{figure}[H]
			\centering 
			\begin{tikzpicture}
				\node (a) [label=above:$\mc G'(i)$] at (0,0) {};
				\node (b) [label=above:$\mc G'(j)$] at (2,0) {};
				\node (c) [draw=none] at (3,0) {};
				\node (d) [draw=none] at (5,0) {};
				\node (e) [label=below left:$i_0$] at (6,-.5) {};
				\node (f) at (6,.5) {};
				\node (g) [label=below right:$j_0$] at (8,-.5) {};
				\node (h) at (8,.5) {};
				\node(i) [draw=none] at (-1.5,0) {$\mathlarger{\mathlarger{\mathlarger{\mc G \bs \mc G' \supseteq}}}$};
				\node (j) [draw=none] at (9,0) {$\mathlarger{\mathlarger{\mathlarger{\sse \mc G}}}$};
				\graph{
					(a) --[yellow] (b),
					(c) ->[ultra thick] (d),
					(f) --[red] (e) --[yellow,dashed] (g) --[red] (h)
				};
				\path
					(a) edge [loop below,green] (a)
					(e) edge [loop below,green,dashed] (e)
					(f) edge [loop above,green,dashed] (f);
			\end{tikzpicture}
			\caption{Example: green loop edge closure condition for quotients, first case (I).}
		\end{figure}
		
		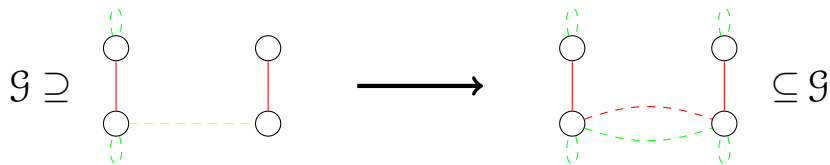
\begin{figure}[H]
			\centering 
			\begin{tikzpicture}
				
				\node (a) at (0,-.5) {};
				\node (b) at (0,.5) {};
				\node(k) at (2,-.5) {};
				\node(l) at (2,.5) {};
				\node (c) [draw=none] at (3,0) {};
				\node (d) [draw=none] at (5,0) {};
				\node (e) at (6,-.5) {};
				\node (f) at (6,.5) {};
				\node (g) at (8,-.5) {};
				\node (h) at (8,.5) {};
				\node(i) [draw=none] at (-1,0) {$\mathlarger{\mathlarger{\mathlarger{\mc G \supseteq}}}$};
				\node (j) [draw=none] at (9,0) {$\mathlarger{\mathlarger{\mathlarger{\sse \mc G}}}$};
				\graph{
					(b) --[red] (a) --[yellow,dashed] (k) --[red] (l);
					(c) ->[ultra thick] (d),
					(f) --[red] (e) --[red,bend left=20,dashed] (g) --[red] (h),
					(e) --[green,bend right=20,dashed] (g)
				};	
				\path
					(a) edge [loop below,green,dashed] (a)
					(b) edge [loop above,green,dashed] (b)
					(e) edge [loop below,green,dashed] (e)
					(f) edge [loop above,green,dashed] (f)
					(g) edge [loop below,green,dashed] (g)
					(h) edge [loop above,green,dashed] (h);
			\end{tikzpicture}
			\caption{Example: green loop edge closure condition for quotients, first case (II).}
		\end{figure}
		
		\begin{figure}[H]
			\centering 
			\begin{tikzpicture}
				\node (a) at (0,-.5) {};
				\node (b) at (0,.5) {};
				\node(k) at (2,-.5) {};
				\node(l) at (2,.5) {};
				\node (c) [draw=none] at (3,0) {};
				\node (d) [draw=none] at (5,0) {};
				\node (e) [label=above:$\mc G'(i)$] at (6,0) {};
				\node (f) [label=above:$\mc G'(j)$] at (8,0) {};
				\node(i) [draw=none] at (-1,0) {$\mathlarger{\mathlarger{\mathlarger{\mc G \supseteq}}}$};
				\node (j) [draw=none] at (9.5,0) {$\mathlarger{\mathlarger{\mathlarger{\sse \mc G \bs \mc G'}}}$};
				\graph{
					(b) --[red] (a) --[red,bend left=20,dashed] (k) --[red] (l),
					(a) --[green,bend right=20,dashed] (k),
					(c) ->[ultra thick] (d),
					(e) --[red,bend left=20] (f),
					(e) --[green,bend right=20] (f)
				};	
				\path
					(a) edge [loop below,green,dashed] (a)
					(b) edge [loop above,green,dashed] (b)
					(k) edge [loop below,green,dashed] (k)
					(l) edge [loop above,green,dashed] (l)
					(e) edge [loop below,green] (e)
					(f) edge [loop below,green] (f);
			\end{tikzpicture}
			\caption{Example: green loop edge closure condition for quotients, first case (III). \qedhere}
		\end{figure}
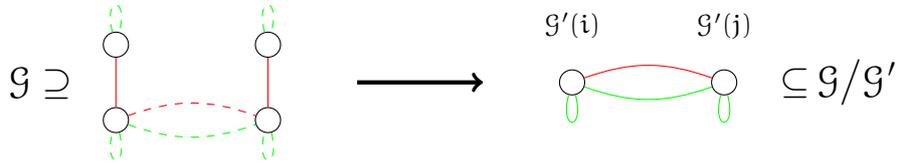
	\end{proof}

	\subsubsection{}
	
	The final case is the following: there is a red component $\mc G'(i) \sse \mc G'$ such that there exists a green straight edge $\set{i_0,\wt i_0} \in \mc G_1 \sm \mc G'_1$ for (distinct) nodes $i_0,\wt i_0 \in \mc G'(i)_0$, and $\mc G'(i)$ is connected by a straight edge $\set{\ol i_0,j_0} \in \mc G_1 \sm \mc G'_1$ to another red component $\mc G'(j) \sse \mc G'$.
	In this case we must distinguish the case where $\mc G'(j)$ is trivial or not:
	
	\begin{figure}[H]
		\centering 
		\begin{tikzpicture}
			\node (a) [label=above:$\mc G'(i)$] at (0,0) {};
			\node (b) [label=above:$\mc G'(j)$] at (2,0) {};
			\node (c) [draw=none] at (3,0) {};
			\node (d) [draw=none] at (5,0) {};
			\node [draw=none,regular polygon, regular polygon sides=5, minimum size=3cm] (B) at (7,0) {};
			\foreach \i in {1,...,5}{
				\node (\i) at (B.corner \i) {};
			}
			\graph{
				(2) --[red,bend left=20] (3) --[red] (4) -- [red] (2),
				(2) --[green,bend right=20,dashed] (3),
				(4) --[yellow,dashed] (1) --[red] (5)
			};
			\node (e) [draw=none] at (7,2) {$j_0$};
			\node (f) [draw=none] at (5,.7) {$i_0$};
			\node (g) [draw=none] at (5.5,-1.2) {$\wt i_0$};
			\node (h) [draw=none] at (8.5,-1.2) {$\ol i_0$};
			\node (i) [draw=none] at (-1.5,0) {$\mathlarger{\mathlarger{\mathlarger{\mc G \bs \mc G' \supseteq}}}$};
			\node (j) [draw=none] at (9.5,0) {$\mathlarger{\mathlarger{\mathlarger{\sse \mc G}}}$};
			\graph{
				(a) --[yellow] (b),
				(c) ->[ultra thick] (d),
			};
			\path
				(a) edge [loop below,green] (a);
		\end{tikzpicture}
		\caption{Example: final configuration for green loop edges in the quotient, first case (I).}
	\end{figure}
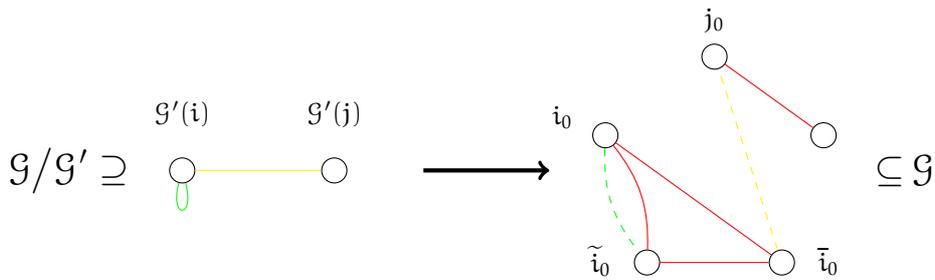
	
	\begin{figure}[H]
		\centering 
		\begin{tikzpicture}
			\node (a) [label=above:$\mc G'(i)$] at (0,0) {};
			\node (b) [label=above:$\mc G'(j)$] at (2,0) {};
			\node (c) [draw=none] at (3,0) {};
			\node (d) [draw=none] at (5,0) {};
			\node [draw=none,regular polygon, regular polygon sides=5, minimum size=3cm] (B) at (7,0) {};
			\foreach \i in {1,...,4}{
				\node (\i) at (B.corner \i) {};
			}
			\graph{
				(2) --[red,bend left=20] (3) --[red] (4) -- [red] (2),
				(2) --[green,bend right=20,dashed] (3),
				(4) --[yellow,dashed] (1)
			};
			\node (e) [draw=none] at (7,2) {$j_0$};
			\node (f) [draw=none] at (5,.7) {$i_0$};
			\node (g) [draw=none] at (5.5,-1.2) {$\wt i_0$};
			\node (h) [draw=none] at (8.5,-1.2) {$\ol i_0$};
			\node (i) [draw=none] at (-1.5,0) {$\mathlarger{\mathlarger{\mathlarger{\mc G \bs \mc G' \supseteq}}}$};
			\node (j) [draw=none] at (9.5,0) {$\mathlarger{\mathlarger{\mathlarger{\sse \mc G}}}$};
			\graph{
				(a) --[yellow] (b),
				(c) ->[ultra thick] (d),
			};
			\path
				(a) edge [loop below,green] (a);
		\end{tikzpicture}
		\caption{Example: final configuration for green loop edges in the quotient, second case (I).}
	\end{figure}
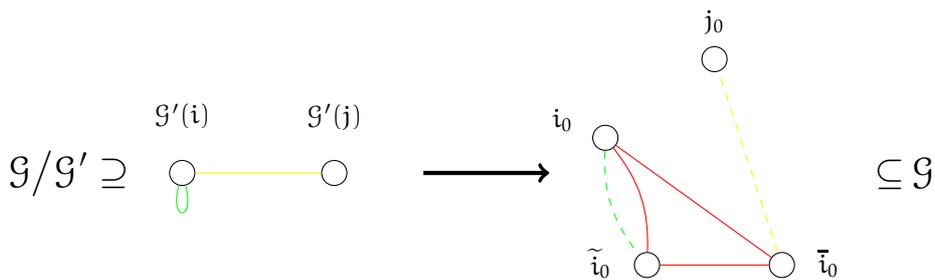

	In both cases, reasoning as in the proof of the second case of Lem.~\ref{lem:red_loop_closure_quotient}, the subgraphs $\mc G'(i), \mc G'(j) \sse \mc G'$ are then contained in a bichromatic component of $\mc G$.
	Hence the straight edge $\Set{ \mc G'(i),\mc G'(j)}$ will double in the quotient, but if $\mc G'(j)$ is trivial the green loop edge $\Set{ \mc G'(j) }$ need \emph{not} arise in $\mc G \bs \mc G'$ (compare the following with Figg.~\ref{fig:second_red_loop_closure_quotient_I} and~\ref{fig:second_red_loop_closure_quotient_II}):
	
	\begin{figure}[H]
		\centering
		\begin{tikzpicture}
			\node [draw=none,regular polygon, regular polygon sides=5, minimum size=3cm] (A) at (0,0) {};
			\foreach \i in {1,...,5}{
				\node (\i) at (A.corner \i) {};
			}
			\node (a) [draw=none] at (1.7,0) {};
			\node (b) [draw=none] at (3.3,0) {};
			\graph{
				(2) --[red,bend left=20] (3) --[red] (4) -- [red] (2),
				(2) --[green,bend right=20,dashed] (3),
				(4) --[yellow,dashed] (1) --[red] (5),
				(a) ->[ultra thick] (b)
			};
			\node [draw=none,regular polygon, regular polygon sides=5, minimum size=3cm] (B) at (5,0) {};
			\foreach \i in {1,...,5}{
				\node (\i) at (B.corner \i) {};
			}
			\graph{
				(1) --[red,bend left=20,dashed] (2) --[red,bend left=20] (3) --[red,bend left=20] (4) --[red,bend left=20,dashed] (5) --[red,bend left=20] (1),
				(1) --[red,bend left=20,dashed] (3) --[red,bend left=20,dashed] (5) --[red,bend left=20,dashed] (2) --[red,bend left=20] (4) --[red,bend left=20,dashed] (1),
				(1) --[green,bend right=20,dashed] (2) --[green,bend right=20,dashed] (3) --[green,bend right=20,dashed] (4) --[green,bend right=20,dashed] (5) --[green,bend right=20,dashed] (1),
				(1) --[green,bend right=20,dashed] (3) --[green,bend right=20,dashed] (5) --[green,bend right=20,dashed] (2) --[green,bend right=20,dashed] (4) --[green,bend right=20,dashed] (1)
			};
			\node (j) [draw=none] at (-2.5,0) {$\mathlarger{\mathlarger{\mathlarger{\mc G \supseteq}}}$};
			\node (k) [draw=none] at (7.5,0) {$\mathlarger{\mathlarger{\mathlarger{\sse \mc G}}}$};
			\end{tikzpicture}
		\caption{Example: final configuration for green loop edges in the quotient, first case (II).}
	\end{figure}
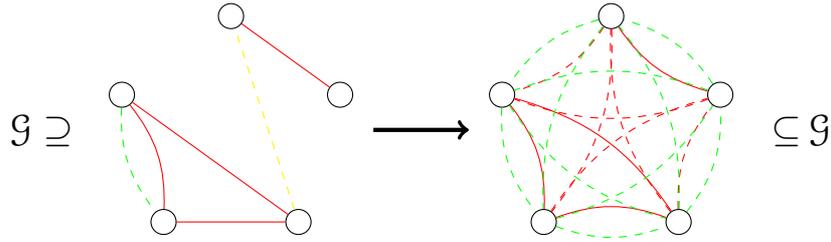
	
	\begin{figure}[H]
		\centering
		\begin{tikzpicture}
			\node [draw=none,regular polygon, regular polygon sides=5, minimum size=3cm] (B) at (0,0) {};
			\foreach \i in {1,...,5}{
				\node (\i) at (B.corner \i) {};
			}
			\node (a) [draw=none] at (1.7,0) {};
			\node (b) [draw=none] at (3.3,0) {};
			\node (e) [label=above:$\mc G'(i)$] at (4,0) {};
			\node (f) [label=above:$\mc G'(j)$] at (6,0) {};
			\node (j) [draw=none] at (-2.5,0) {$\mathlarger{\mathlarger{\mathlarger{\mc G \supseteq}}}$};
			\node (k) [draw=none] at (7.5,0) {$\mathlarger{\mathlarger{\mathlarger{\sse \mc G \bs \mc G'}}}$};
			\graph{
				(1) --[red,bend left=20,dashed] (2) --[red,bend left=20] (3) --[red,bend left=20] (4) --[red,bend left=20,dashed] (5) --[red,bend left=20] (1),
				(1) --[red,bend left=20,dashed] (3) --[red,bend left=20,dashed] (5) --[red,bend left=20,dashed] (2) --[red,bend left=20] (4) --[red,bend left=20,dashed] (1),
				(1) --[green,bend right=20,dashed] (2) --[green,bend right=20,dashed] (3) --[green,bend right=20,dashed] (4) --[green,bend right=20,dashed] (5) --[green,bend right=20,dashed] (1),
				(1) --[green,bend right=20,dashed] (3) --[green,bend right=20,dashed] (5) --[green,bend right=20,dashed] (2) --[green,bend right=20,dashed] (4) --[green,bend right=20,dashed] (1),
				(a) ->[ultra thick] (b),
				(e) --[red,bend left=20] (f),
				(e) --[green,bend right=20] (f)
			};
			\path
				(e) edge [loop below,green] (e)
				(f) edge [loop below,green] (f);
		\end{tikzpicture}
		\caption{Example: final configuration for green loop edges in the quotient, first case (III).}
	\end{figure}
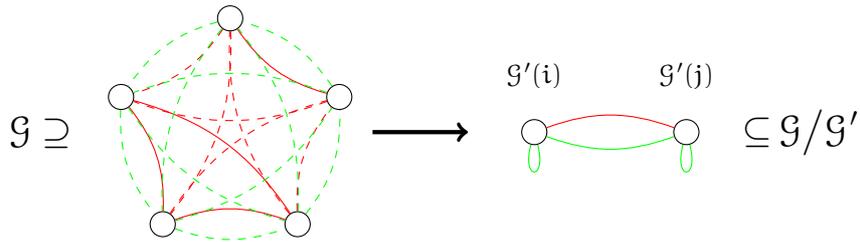
	
	\begin{figure}[H]
		\centering
		\begin{tikzpicture}
			\node [draw=none,regular polygon, regular polygon sides=5, minimum size=3cm] (A) at (0,0) {};
			\foreach \i in {1,...,4}{
				\node (\i) at (A.corner \i) {};
			}
			\node (a) [draw=none] at (1.7,0) {};
			\node (b) [draw=none] at (3.3,0) {};
			\graph{
				(2) --[red,bend left=20] (3) --[red] (4) -- [red] (2),
				(2) --[green,bend right=20,dashed] (3),
				(4) --[yellow,dashed] (1),
				(a) ->[ultra thick] (b)
			};
			\node [draw=none,regular polygon, regular polygon sides=5, minimum size=3cm] (B) at (5,0) {};
			\foreach \i in {1,...,4}{
				\node (\i) at (B.corner \i) {};
			}
			\graph{
				(1) --[red,bend left=20,dashed] (2) --[red,bend left=20] (3) --[red,bend left=20] (4) --[red,bend left=20,dashed] (1),
				(1) --[green,bend right=20,dashed] (2) --[green,bend right=20,dashed] (3) --[green,bend right=20,dashed] (4) --[green,bend right=20,dashed] (1),
				(1) --[green,bend right=20,dashed] (3),
				(1) --[red,bend left=20,dashed] (3),
				(2) --[red,bend left=20] (4),
				(2) --[green,bend right=20,dashed] (4)
			};
			\node (j) [draw=none] at (-2.5,0) {$\mathlarger{\mathlarger{\mathlarger{\mc G \supseteq}}}$};
			\node (k) [draw=none] at (7.5,0) {$\mathlarger{\mathlarger{\mathlarger{\sse \mc G}}}$};
			\end{tikzpicture}
		\caption{Example: final configuration for green loop edges in the quotient, second case (II).}
	\end{figure}
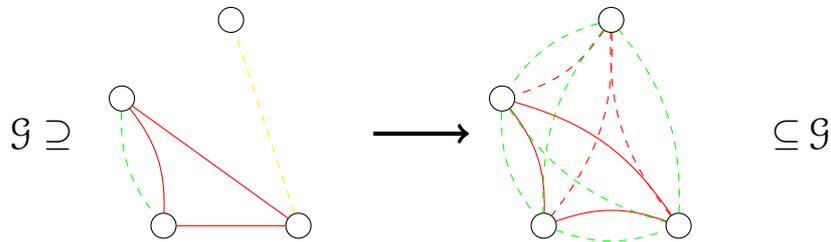
	
	\begin{figure}[H]
		\centering
		\begin{tikzpicture}
			\node [draw=none,regular polygon, regular polygon sides=5, minimum size=3cm] (B) at (0,0) {};
			\foreach \i in {1,...,4}{
				\node (\i) at (B.corner \i) {};
			}
			\node (a) [draw=none] at (1.7,0) {};
			\node (b) [draw=none] at (3.3,0) {};
			\node (e) [label=above:$\mc G'(i)$] at (4,0) {};
			\node (f) [label=above:$\mc G'(j)$] at (6,0) {};
			\node (j) [draw=none] at (-2.5,0) {$\mathlarger{\mathlarger{\mathlarger{\mc G \supseteq}}}$};
			\node (k) [draw=none] at (7.5,0) {$\mathlarger{\mathlarger{\mathlarger{\sse \mc G \bs \mc G'}}}$};
			\graph{
				(1) --[red,bend left=20,dashed] (2) --[red,bend left=20] (3) --[red,bend left=20] (4) --[red,bend left=20,dashed] (1),
				(1) --[green,bend right=20,dashed] (2) --[green,bend right=20,dashed] (3) --[green,bend right=20,dashed] (4) --[green,bend right=20,dashed] (1),
				(1) --[green,bend right=20,dashed] (3),
				(1) --[red,bend left=20,dashed] (3),
				(2) --[red,bend left=20] (4),
				(2) --[green,bend right=20,dashed] (4),
				(a) ->[ultra thick] (b),
				(e) --[red,bend left=20] (f),
				(e) --[green,bend right=20] (f)
			};
			\path
				(e) edge [loop below,green] (e);
		\end{tikzpicture}
		\caption{Example: final configuration for green loop edges in the quotient, second case (III).}
	\end{figure}
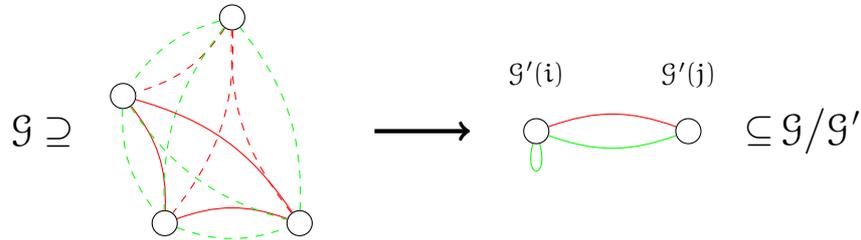
		
	This completes the proof of Thm.~\ref{thm:quasi_crystallographs_classify_quotients}.
	
	\subsection{Exotic components}
	
	It follows that there are two new connected components of quasi-crystallographs: starting from a component of type $B$ or $D$, we can glue green loop edges to any subset of nodes.
	For integers $r,s \geq 0$ we will say these are bichromatic graphs of type $B_{r+s}C_r$ and $C_rD_{r+s}$, respectively, where $r$ is the number of nodes with green loop edges, and $s$ the number of nodes without---so they have $n = r+s$ nodes, and we have the tautological identities $B_rC_r = BC_r$, $B_sC_0 = B_s$, $C_0D_s = D_s$ and $C_rD_r = C_r$.
	
	\begin{figure}[H]
		\centering
		\begin{tikzpicture}
			\node (a) at (0,0) {};
			\node (b) at (2,0) {};
			\node (c) at (2,2) {};
			\node (d) at (0,2) {};
			\node(e) [draw=none] at (1,-.5) {$\mc G^{B_4C_2}$};
			\node (f) at (4,0) {};
			\node (g) at (6,0) {};
			\node (h) at (6,2) {};
			\node (i) at (4,2) {};
			\node (j) [draw=none] at (5,-.5) {$\mc G^{C_2D_4}$};
			\graph{
				(a) --[red,bend left=20] (b) --[red,bend left=20] (c) --[red,bend left=20] (d) --[red,bend left=20] (a),
				(a) --[green,bend right=20] (b) -- [green,bend right=20] (c) --[green,bend right=20] (d) --[green,bend right=20] (a),
				(a) --[red,bend left=20] (c),
				(a) --[green,bend right=20] (c),
				(b) --[red,bend left=20] (d),
				(b) --[green,bend right=20] (d),
				(f) --[red,bend left=20] (g) -- [red,bend left=20] (h) --[red,bend left=20] (i) --[red,bend left=20] (f),
				(f) --[green,bend right=20] (g) -- [green,bend right=20] (h) --[green,bend right=20] (i) --[green,bend right=20] (f),
				(f) --[red,bend left=20] (h),
				(f) --[green,bend right=20] (h),
				(g) --[red,bend left=20] (i),
				(g) --[green,bend right=20] (i)
			};
			\path
				(a) edge [red,loop below] (a)
				(b) edge [red,loop below] (b)
				(c) edge [red,loop above] (c)
				(d) edge [red,loop above] (d)
				(b) edge [green,loop right] (b)
				(d) edge [green,loop left] (d)
				(f) edge [green,loop below] (f)
				(h) edge [green,loop above] (h);
		\end{tikzpicture}
		\caption{Examples of ''exotic`` quasi-crystallographs.}
	\end{figure}
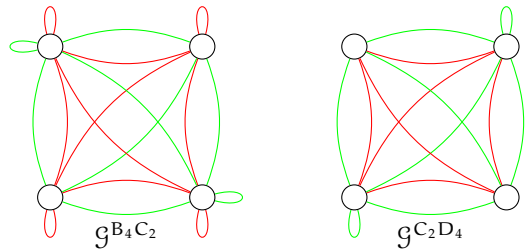
	
	\begin{rema}
		\label{rk:exotic_arrangement}
		
		The hyperplane arrangement associated with the symmetric subset $\Phi^{B_{r+s}C_r} \sse BC_n$ is of type $B_n/C_n$, where $n = r+s$, so this is always crystallographic.
		The point is dilating any linear functional will not change the arrangement, so there is no difference in having one or two loop edges at a node.
		
		On the contrary there \emph{is} a difference between having zero or one loop edges, so $\Phi^{C_rD_{r+s}} \sse BC_n$ leads to an ''exotic`` hyperplane arrangement, which is not always crystallographic (as mentioned above).
		We say it is of type $(B_r/C_r)D_{r+s}$---slightly changing the notation of~\cite{doucot_rembado_tamiozzo_2022_local_wild_mapping_class_groups_and_cabled_braids}.
	\end{rema}
	
	\subsubsection{}
	
	Putting all together we have proven the following classification of quasi-crystallographs, in view of Thm.~\ref{thm:crystallograph_classification}.
	
	\begin{coro}
		\label{cor:quasi_crystallographs_classification}
		
		All quasi-crystallographs are disjoint unions of bichromatic graphs of type $A_n$, $B_n$, $C_n$, $D_n$, $BC_n$, $B_{r+s}C_r$ or $C_rD_{r+s}$, for integers $n,r,s \geq 0$.	
	\end{coro}
	
	\subsection{}
	
	Finally we can prove the converse of Thm.~\ref{thm:quasi_crystallographs_classify_quotients}, so that in the end this is a classification of all possible restrictions~\eqref{eq:restricted_system} of root systems onto the kernels of their subsystems, up to equivalence---which was one of the main goals.
	
	\begin{prop}
		\label{prop:quasi_crystallographs_are_quotients}
		
		All quasi-crystallographs are quotient graphs.
	\end{prop}
	
	\begin{proof}
		By Cor.~\ref{cor:quasi_crystallographs_classification}, it is enough to prove the statement for all components, since taking quotients and disjoint unions of crystallographs are commutative operations: if $\mc G' = \coprod_I \mc G'_i$ is a subgraph of $\mc G = \coprod_I \mc G_i$, so that $\mc G'_i \sse \mc G_i$ for $i \in I$, then by definition
		\begin{equation*}
			\mc G \bs \mc G' = \coprod_I \mc G_i \bs \mc G'_i \, .
		\end{equation*}
		Moreover if $\mc G'$ is a totally disconnected graph on $n \geq 1$ nodes then $\mc G \bs \mc G' = \mc G$, so certainly all crystallographic components arise.
		
		Finally consider the crystallograph $\mc G_{r,s}$ obtained as the disjoint union of $r \geq 0$ graphs of type $A_1$, and $s \geq 0$ graphs of type $A_0$:
		
		\begin{figure}[H]
			\centering
			\begin{tikzpicture}
				\node (a) at (0,0) {};
				\node (b) [draw=none] at (1,0) {$\dm$};
				\node (c) at (2,0) {};
				\node (d) at (0,1) {};
				\node (e) at (0,2) {};
				\node (f) [draw=none] at (1,1.5) {$\dm$};
				\node (g) at (2,1) {};
				\node (h) at (2,2) {};
				\node (i) [draw=none] at (-1.3,1) {$\mathlarger{\mathlarger{\mathlarger{\mc G_{r,s}=}}}$};
				\node (j) [draw=none] at (-.3,2.5) {};
				\node (k) [draw=none] at (2.3,2.5) {};
				\node (l) [draw=none] at (1,3) {$s \text{ times}$};
				\node (m) [draw=none] at (-.3,-.5) {};
				\node (n) [draw=none] at (2.3,-.5) {};
				\node (o) [draw=none] at (1,-1) {$r \text{ times}$};
				\graph{
					(d) --[red] (e),
					(g) --[red] (h)
				};
				\draw[ultra thick,decoration={calligraphic brace},decorate] (j)--(k);
				\draw[ultra thick,decoration={calligraphic brace,mirror},decorate] (m)--(n);
			\end{tikzpicture}
			\caption{The graph $\mc G_{r,s}$.}
		\end{figure}
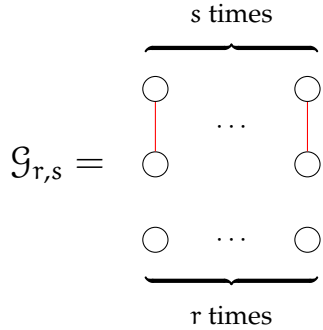

		This is naturally embedded as a subgraph of both $B_{r+2s}$ and $D_{r+2s}$, and one readily verifies that 
		\begin{equation*}
			\mc G^{B_{r+2s}} \bs \mc G_{r,s} = \mc G^{B_{r+s}C_r} \, , \qquad \mc G^{D_{r+2s}} \bs \mc G_{r,s} = \mc G^{C_rD_{r+s}} \, .
		\end{equation*}
		Indeed in both cases the quotient has $r+s$ nodes and all possible straight edges, and in the former case it also has all red loop edges; finally there are $s$ green loop edges at the nodes corresponding to the nontrivial red components of $\mc G_{r,s}$---and no other.
	\end{proof}

	\subsubsection{}
	
	Hence in conclusion all the information about restricted hyperplane arrangements is encoded in a quasi-crystallograph, but there is some clear redundancy that will be taken care of in the next section.
	
	\section{Hyperplane arrangements}
	\label{sec:arrangements}
	
	\subsection{}
	
	Of course $\Ker(\Phi) \sse V = \mb C^n$ only depends on $\Phi \sse V^{\dual} \sm \set{0}$ up to dilation of each covector, i.e. on the subset $\mb P(\Phi) \ceqq \pi(\Phi) \sse \mb P \bigl( V^{\dual} \bigr)$ of the projective dual space, using the canonical projection 
	\begin{equation*}
		\pi \cl V^{\dual} \sm \set{0} \lra \mb P \bigl( V^{\dual} \bigr) \eqqcolon \bigl( V^{\dual} \sm \set{0} \bigr) \bs \mb C^{\times} \, .
	\end{equation*}
	In particular in our situation we may restrict to positive systems $\Phi_+ \sse \Phi$ of roots, and further to reduced ones.
	
	Now the bichromatic graph of a simply-laced root system is already naturally associated with a choice of a positive system, since it has one (unoriented) edge for each opposite pair $\pm \alpha \in \Phi$; on the contrary we kept track of short/long roots for non-simply-laced types, by having loop edges of two different colours, which now ought to be fused into a single loop edge (of a different colour).
	
	Introduce thus the extended set of (primary) colours $\set{R,G,B}$, where now ''blue`` is allowed.

	\begin{defi}
		A \emph{trichromatic} graph is a graph $\mc G = (\mc G_0,\mc G_1,m)$ equipped with a colour function $c \cl \mc G_1 \to \set{R,G,B}$, such that $c(e,m_e) \neq c(e',m_{e'})$ if $e = e' \sse \mc G_0$.
		
		A trichromatic \emph{subgraph} of $(\mc G,c)$ is a subgraph $\mc G' \sse \mc G$ equipped with the restricted colour function $c' = \eval[1]{c}_{\mc G'_1}$.
	\end{defi}
	
	\subsubsection{}
	
	By definition bichromatic graphs are (all the) trichromatic graphs $(\mc G,c)$ with $c(\mc G_1) \sse \set{R,G}$.
	This way the relevant terminology about bichromatic graphs can be naturally translated here---e.g. being simply-laced, etc.
	
	Now our viewpoint is that the adjacency of a trichromatic graph on $n \geq 1$ nodes can be used to encode certain hyperplane arrangements in $V = \mb C^n$, viz. certain lines in $V^{\dual}$, as follows.
	For $i \neq j \in \ul n$ denote 
	\begin{equation*}
		H^{\pm}_{ij} \ceqq \Ker \bigl( \alpha^{\pm}_{ij} \bigr) = \Set{ \sum_k \alpha_k e_k \in V | \alpha_i \pm \alpha_j = 0 } \sse V \, ,
	\end{equation*}
	and 
	\begin{equation*}
		H_i \ceqq \Ker(\alpha_i) = \Set{ \sum_k \alpha_k e_k | \alpha_i = 0 } \sse V \, ,
	\end{equation*}
	identifying the covectors $\alpha_i = e_i^{\dual} \in V^{\dual}$ with the linear coordinates $V \to \mb C$ in the basis $(e_1,\dc,e_n)$.
	Then the root-hyperplane arrangements of the classical root systems are 
	\begin{equation*}
		\mc H^{A_{n-1}} = \Set{ H^-_{ij} | i \neq j \in \ul n } \, , \quad \mc H^{D_n} = \mc H^{A_{n-1}} \cup \Set{ H^+_{ij} | i \neq j \in \ul n } \sse \mb P \bigl( V^{\dual} \bigr) \, ,
	\end{equation*}
	and 
	\begin{equation*}
		\mc H^{B_n} = \mc H^{C_n} = \mc H^{BC_n} = \mc H^{D_n} \cup \Set{ H_i | i \in \ul n} \sse \mb P \bigl( V^{\dual} \bigr) \, .
	\end{equation*}
	
	Now if $\mc H \sse \mc H^{BC_n}$ we associate a trichromatic graph $\mc G^{\mc H}$ with it.
	It has nodes $\mc G^{\mc H}_0 = \ul n$, and its adjacency/colouring are prescribed as follows:
	\begin{enumerate}
		\item if $H^-_{ij} \in \mc H$, put a red straight edge $e = \set{i,j} \in \mc G^{\mc H}_1$;
			
		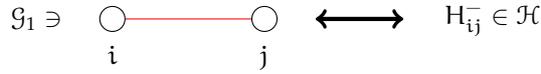
\begin{figure}[H]
			\centering
			\begin{tikzpicture}
				\node (a) [label=below:$i$] at (0,0) {};
				\node (b) [label=below:$j$] at (2,0) {};
				\node (c) [draw=none] at (2.5,0) {};
				\node (d) [draw=none] at (4,0) {};
				\node (e) [draw=none] at (5,0) {$H^-_{ij} \in \mc H$};
				\node (f) [draw=none] at (-1,0) {$\mc G_1 \ni$};
				\graph{
					(a) --[red] (b),
					(c) --[to-to,ultra thick] (d)
				};
			\end{tikzpicture}
			\caption{Projective correspondence, straight edges (I).}
		\end{figure}
		
		\item if $H^+_{ij} \in \mc H$, put a green straight edge $e = \set{i,j} \in \mc G^{\mc H}_1$;
		
		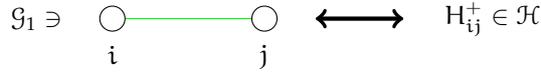
\begin{figure}[H]
			\centering
			\begin{tikzpicture}
				\node (a) [label=below:$i$] at (0,0) {};
				\node (b) [label=below:$j$] at (2,0) {};
				\node (c) [draw=none] at (2.5,0) {};
				\node (d) [draw=none] at (4,0) {};
				\node (e) [draw=none] at (5,0) {$H^+_{ij} \in \mc H$};
				\node (f) [draw=none] at (-1,0) {$\mc G_1 \ni$};
				\graph{
					(a) --[green] (b),
					(c) --[to-to,ultra thick] (d)
				};
			\end{tikzpicture}
			\caption{Projective correspondence, straight edges (II).}
		\end{figure}
		
		\item if $H_k \in \mc H$, put a blue loop edge $l = \set{k} \in \mc G^{\mc H}_1$.
			
		\begin{figure}[H]
			\centering
			\begin{tikzpicture}
				\node (a) [label=below:$k$] at (0,0) {};
				\node (c) [draw=none] at (.5,0) {};
				\node (d) [draw=none] at (2,0) {};
				\node (e) [draw=none] at (3,0) {$H_k \in \mc H$};
				\node (f) [draw=none] at (-1,0) {$\mc G_1 \ni$};
				\graph{
					(c) --[to-to,ultra thick] (d)
				};
				\path
				(a) edge [loop above,blue] (a);
			\end{tikzpicture}
			\caption{Projective correspondence, loop edges.}
		\end{figure}
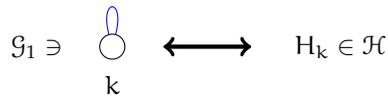
	\end{enumerate}
	
	The conditions for the straight edges are completely analogous to the bichromatic case of \S~\ref{sec:crystallographs}, while the case of loop edges should be compared with Figg.~\ref{fig:main_correspondence_loop_edges_I} and~\ref{fig:main_correspondence_loop_edges_II} (which are now ''fused``).
	
	\begin{rema}
		Just as in the case of root systems, a disjoint-union decomposition of $\mc G^{\mc H}$ corresponds to a direct sum decomposition of $\mc H$, and the number of hyperplanes is $\abs{\mc H} = \abs{ \mc G^{\mc H}_1}$ (cf. Rk.~\ref{rk:graph_root_system_correspondence}).
	\end{rema}
	
	\subsubsection{}
	
	The above prescription singles out a special class of trichromatic graphs, which are in natural inclusion-preserving bijection with sub-arrangements of the root-hyperplane arrangement of type $B/C$: in particular $\ol{\mc G^{\mc H}} \sse \mc G^{\mc H}$ is a (simply-laced) bichromatic graph as above, and the full graph is obtained by gluing blue loop edges at some nodes.
	
	\begin{exem}
		The classical root-hyperplane arrangements yield the following trichromatic graphs: the same bichromatic graphs as in Ex.~\ref{ex:classical_crystallographs} for the simply-laced cases, and a simply-laced complete bichromatic graph with blue loop edges glued at each node for types $B/C/BC$.
		Hence all non-simply-laced cases are ''fused`` into a single one; see below the case on $n = 4$ nodes:
		\begin{figure}[H]
			\centering
			\begin{tikzpicture}
				\node (a) at (0,0) {};
				\node (b) at (2,0) {};
				\node (c) at (2,2) {};
				\node (d) at (0,2) {};
				\node(e) [draw=none] at (1,-1) {$\mc G^{\mc H_{B_4}} = \mc G^{\mc H_{C_4}} = \mc G^{\mc H_{BC_4}}$};
				\graph{
					(a) --[red,bend left=20] (b) --[red,bend left=20] (c) --[red,bend left=20] (d) --[red,bend left=20] (a),
					(a) --[green,bend right=20] (b) -- [green,bend right=20] (c) --[green,bend right=20] (d) --[green,bend right=20] (a),
					(a) --[red,bend left=20] (c),
					(a) --[green,bend right=20] (c),
					(b) --[red,bend left=20] (d),
					(b) --[green,bend right=20] (d)
				};
				\path
					(a) edge [blue,loop below] (a)
					(b) edge [blue,loop below] (b)
					(c) edge [blue,loop above] (c)
					(d) edge [blue,loop above] (d);
			\end{tikzpicture}
			\caption{Example: the rank-four non-simply-laced ''classical`` trichromatic graph. \qedhere}
		\end{figure}
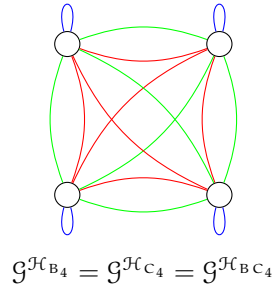
	\end{exem}
	
	\subsection{}
	
	Of course not all subsets $\mc H \sse \mc H^{BC_n}$ correspond to \emph{root}-hyperplane arrangements, i.e. they are not the projectification of root subsystems $\Phi \sse BC_n$ (and so are \emph{not} closed under mutual reflections): but those which do happen to be crystallographic are captured by a variation of the main definition.
	
	\begin{defi}
		\label{def:projective_crystallographs}
		
		A \emph{projective} crystallograph is a trichromatic graph $(\mc G,c)$ such that $c(\ol{\mc G}) \sse \set{R,G}$ and $c(l) = B$ for all loop edges $l \in \mc G_1$, and satisfying the two conditions of Def.~\ref{def:crystallographs}. 
	\end{defi}
	
	\subsubsection{}
	
	By definition the ''simply-laced`` condition yields the same local pictures as in the bichromatic case, viz. Figg.~\ref{fig:straight_edge_closure_condition_crystallographs_I}~\ref{fig:straight_edge_closure_condition_crystallographs_II} and~\ref{fig:straight_edge_closure_condition_crystallographs_III}.
	The conditions involving loop edges instead become the following (compare with Figg.~\ref{fig:red_loop_edge_closure_crystallograph_I},~\ref{fig:red_loop_edge_closure_crystallograph_II},~\ref{fig:green_loop_edge_closure_crystallograph_I} and~\ref{fig:green_loop_edge_closure_crystallograph_II}):	
	
	\begin{figure}[H]
		\centering
		\begin{tikzpicture}
			\node (a) at (0,0) {};
			\node (b) at (2,0) {};
			\node (c) [draw=none] at (2.5,0) {};
			\node (d) [draw=none] at (4.5,0) {};
			\node (e) at (5,0) {};
			\node (f) at (7,0) {};
			\graph{
				(a) --[red] (b),
				(c) ->[ultra thick] (d),
				(e) --[red,bend right=20] (f),
				(e) --[green,bend left=20] (f)
			};
			\path
				(a) edge [blue,loop below] (a)
				(e) edge [blue,loop below] (e)
				(f) edge [blue,loop below] (f);
		\end{tikzpicture}
		\caption{Loop edge closure condition for projective crystallographs (I).} 
	\end{figure}
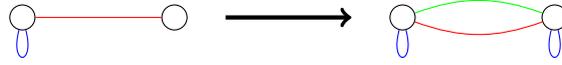
	
	\begin{figure}[H]
		\centering
		\begin{tikzpicture}
			\node (a) at (0,0) {};
			\node (b) at (2,0) {};
			\node (c) [draw=none] at (2.5,0) {};
			\node (d) [draw=none] at (4.5,0) {};
			\node (e) at (5,0) {};
			\node (f) at (7,0) {};
			\graph{
				(a) --[green] (b),
				(c) ->[ultra thick] (d),
				(e) --[red,bend right=20] (f),
				(e) --[green,bend left=20] (f)
			};
			\path
				(a) edge [blue,loop below] (a)
				(e) edge [blue,loop below] (e)
				(f) edge [blue,loop below] (f);
		\end{tikzpicture}
		\caption{Loop edge closure condition for projective crystallographs (II).} 
	\end{figure}
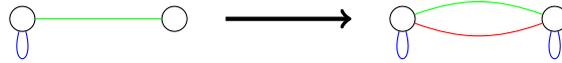
	
	Any such graph is associated with a root-hyperplane arrangement $\mc H = \mc H^{\mc G,c}$, inverting the above prescription.
	
	Then running the same argument in the proof of Thm.~\ref{thm:crystallograph_classification} yields the following classification.
	
	\begin{coro}
		\label{cor:projective_crystallograph_classification}
		
		Let $(\mc G,c)$ be a projective crystallograph on $n \geq 1$ nodes.
		Then $\mc G$ is a disjoint union of the ''classical`` (projective crystallo)graphs $\mc G^{\mc H_{A_{m-1}}}$, $\mc G^{\mc H_{D_m}}$ and $\mc G^{\mc H_{B_m}} = \mc G^{\mc H_{C_m}} = \mc G^{\mc H_{BC_m}}$ ($m \leq n$), and of the simply-laced bichromatic graphs $\mc G^{d_1,d_2}$ ($d_1 + d_2 \leq n$).
	\end{coro}
	
	\subsubsection{}
	
	Again up to acting via the Weyl group we can restrict to ''classical`` components only (cf. Lem.~\ref{lem:bipartite_graph_is_type_A}).
	
	\subsection{}
	
	The last idea is that there is a natural operation on graphs which mimics taking the projection $\pi \cl V^{\dual} \sm \set{0} \to \mb P \bigl( V^{\dual} \bigr)$.
	Namely, if $(\mc G,c)$ is \emph{bi}chromatic, we can associate to it a \emph{tri}chromatic graph by replacing all red/green loop edges by a (single) blue loop edge.
	The resulting trichromatic graph is denoted $\mb P(\mc G,c)$, and called its \emph{projectification}; see an example below:
	
	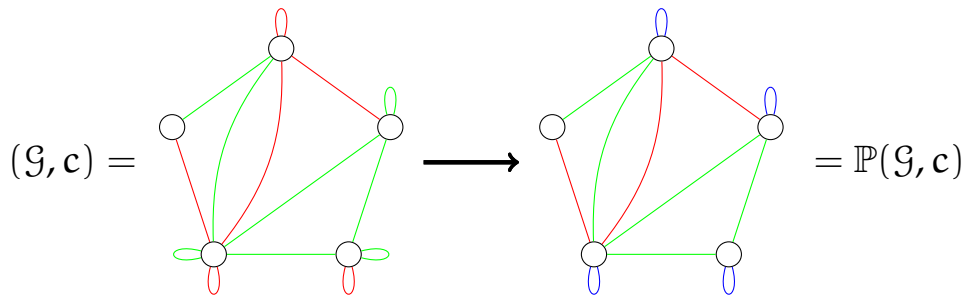
\begin{figure}[H]
		\centering
		\begin{tikzpicture}
			\node [draw=none,regular polygon, regular polygon sides=5, minimum size=3cm] (A) at (0,0) {};
			\foreach \i in {1,...,5}{
				\node (\i) at (A.corner \i) {};
			}
			\node (a) [draw=none] at (1.7,0) {};
			\node (b) [draw=none] at (3.3,0) {};
			\node (c) [draw=none] at (-2.75,0) {$\mathlarger{\mathlarger{\mathlarger{(\mc G,c) =}}}$};
			\graph{
				(1) --[green] (2) --[red] (3) --[green] (4) --[green] (5) --[red] (1), 
				(1) --[red,bend left = 20] (3) --[green] (5),
				(1) --[green,bend right=20] (3),
				(a) ->[ultra thick] (b)
			};
			\path
				(1) edge [loop above,red] (1)
				(3) edge [loop below,red] (3)
				(3) edge [loop left,green] (3)
				(4) edge [loop below,red] (4)
				(4) edge [loop right,green] (4)
				(5) edge [loop above,green] (5);
			\node [draw=none,regular polygon, regular polygon sides=5, minimum size=3cm] (A) at (5,0) {};
			\foreach \i in {1,...,5}{
				\node (\i) at (A.corner \i) {};
			}
			\node (a) [draw=none] at (1.7,0) {};
			\node (b) [draw=none] at (3.3,0) {};
			\node (c) [draw=none] at (8,0) {$\mathlarger{\mathlarger{\mathlarger{= \mb P(\mc G,c)}}}$};
			\graph{
				(1) --[green] (2) --[red] (3) --[green] (4) --[green] (5) --[red] (1), 
				(1) --[red,bend left = 20] (3) --[green] (5),
				(1) --[green,bend right=20] (3),
				(a) ->[ultra thick] (b)
			};
			\path
				(1) edge [loop above,blue] (1)
				(3) edge [loop below,blue] (3)
				(4) edge [loop below,blue] (4)
				(5) edge [loop above,blue] (5);
		\end{tikzpicture}
		\caption{Example: projectification of a bichromatic graph.}
	\end{figure}
	
	\subsubsection{}
	
	Now if $\mc G = \mc G^{\Phi}$ for some symmetric subset $\Phi \sse V^{\dual}$ (in the dual reading), then $\mb P(\mc G)$ encodes the projectification $\mb P(\Phi) \sse \mb P \bigl( V^{\dual} \bigr)$: indeed this latter is obtained by replacing each pair $\pm \alpha \in \Phi$ of opposite roots with the line it generates, i.e. with the hyperplane $\Ker(\pm \alpha) \sse V$, so in particular the overall operation ''fuses`` the short/long roots $\alpha_i, 2\alpha_i \in V^{\dual}$ into a single element.
	
	\begin{theo}
		\label{thm:projective_crystallograph_classificy_hyperplane_arrangements}
			
		A sub-arrangement $\mc H \sse BC_n$ is crystallographic if and only if $\mc G^{\mc H}$ is a projective crystallograph, or equivalently if and only if it is the projectification of a crystallograph.
	\end{theo}
	
	\begin{proof}
		The first statement is a corollary of Thm.~\ref{thm:crystallographs_classify_root_systems}, since the Weyl element associated to $\alpha_i$ and $2\alpha_i$ are the same for $i \in \ul n$.
		
		The second statement follows from the classifications of Thm.~\ref{thm:crystallograph_classification} and Cor.~\ref{cor:projective_crystallograph_classification}.
		For the less trivial implication suppose $(\mc G,c)$ is a projective crystallograph: then repainting its (blue) loop edges---if any---in red provides a ''lifted`` bichromatic graph $(\wt{\mc G},\wt c)$ such that $\mb P(\wt{\mc G},\wt c) = (\mc G,c)$, and this latter is a crystallograph.\fn{
			Of course other lifts are possible, e.g. putting green loop edges, which corresponds to the fact that the (dual) root systems of type $B$ and $C$ have the same root-hyperplane arrangement.}
	\end{proof}
	
	\subsubsection{}
	
	Hence up to isomorphism we obtain a classification of the hyperplane arrangements of root subsystems of all classical root systems, reading from Cor.~\ref{cor:root_system_classification}: in particular they are direct sums of classical root-hyperplane arrangements.
	
	\subsection{About quotients}
	
	Finally an analogous variation of the main definition can be given for quasi-crystallographs, with a view towards encoding the hyperplane arrangements of quotients of root systems.
			
	Namely we can define the quotient of two nested projective crystallographs using the same algorithm of Def.~\ref{def:quotient_graph}, but replacing red/green loop edges with blue loop edges throughout.
	
	\begin{exem}[Adjacency of a quotient of projective crystallographs]
		Compare the following pictures with those of Ex.~\ref{ex:quotient_graph}.
		
		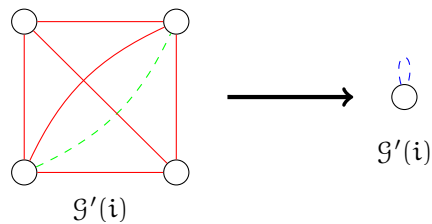
\begin{figure}[H]
				\centering
				\begin{tikzpicture}
					\node (a) at (0,0) {};
					\node (b) at (2,0) {};
					\node (c) at (2,2) {};
					\node (d) at (0,2) {};
					\node (e) [draw=none] at (2.5,1) {};
					\node (f) [draw=none] at (4.5,1) {};
					\node (g) [label=below:$\mc G'(i)$] at (5,1) {};
					\node (h) [draw=none] at (1,-.5) {$\mc G'(i)$};
					\graph{
						(a) --[red] (b) --[red] (c) --[red] (d) --[red] (a),
						(a) --[red,bend left = 20] (c),
						(a) --[green,bend right = 20,dashed] (c),
						(b) --[red] (d),
						(e) ->[ultra thick] (f),
					};
					\path
						(g) edge [loop above,blue,dashed] (g);
				\end{tikzpicture}
				\caption{Example: loop edge in $\mc G \bs \mc G'$, arising from a green straight edge within a red component of $\mc G'$.}
			\end{figure}
			
			\begin{figure}[H]
			\centering
			\begin{tikzpicture}
				\node (a) at (0,0) {};
				\node (b) at (2,0) {};
				\node (c) at (1,sqrt 3) {};
				\node (d) [draw=none] at (1,-.5) {$\mc G'(i)$};
				\node (e) at (4,0) {};
				\node (f) at (6,0) {};
				\node (g) at (6,2) {};
				\node (h) at (4,2) {};
				\node (i) [draw=none] at (5,-.5) {$\mc G'(j)$};
				\node (j) [draw=none] at (6.5,1) {};
				\node (k) [draw=none] at (8.5,1) {};
				\node (l) [label=below:$\mc G'(i)$] at (9,1) {};
				\graph{
					(a) --[red] (b) --[red] (c) --[red] (a),
					(e) --[red,bend left=20] (f) --[red,bend left=20] (g) --[red,bend left=20] (h) --[red,bend left=20] (e),
					(e) --[green,bend right=20] (f) --[green,bend right=20] (g) --[green,bend right=20] (h) --[green,bend right=20] (e),
					(e) --[red,bend left=20] (g),
					(f) --[red,bend left=20] (h),
					(e) --[green,bend right=20] (g),
					(f) --[green,bend right=20] (h),
					(b) --[red,bend left=20, dashed] (e),
					(b) --[green,bend right=20, dashed] (e),
					(j) ->[ultra thick] (k),
				};
				\path
					(l) edge [loop above,blue,dashed] (l);
			\end{tikzpicture}
			\caption{Example: loop edge in $\mc G \bs \mc G'$, arising from a straight edge from a red component to a bichromatic component of $\mc G'$.}
		\end{figure}
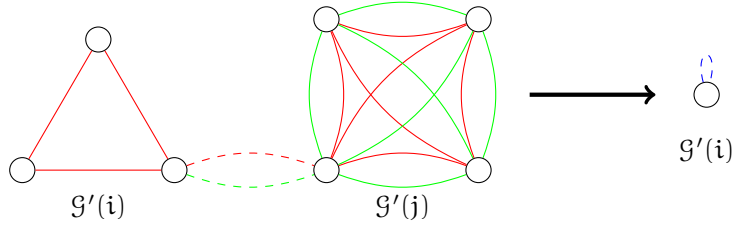

		\begin{figure}[H]
			\centering
			\begin{tikzpicture}
				\node (a) at (0,0) {};
				\node (b) at (2,0) {};
				\node (c) at (1,sqrt 3) {};
				\node (d) [draw=none] at (1,-.5) {$\mc G'(i)$};
				\node (e) [draw=none] at (2,sqrt 3/2) {};
				\node (f) [draw=none] at (4,sqrt 3/2) {};
				\node (g) [label=below:$\mc G'(i)$] at (4.5,sqrt 3/2) {};
				\graph{
					(a) --[red] (b) --[red] (c) --[red] (a),
					(e) ->[ultra thick] (f)
				};
				\path
					(b) edge [loop right,blue,dashed] (b)
					(g) edge [loop right,blue,dashed] (g);
			\end{tikzpicture}
			\caption{Loop edge in $\mc G \bs \mc G'$, arising from a loop edge within a red component of $\mc G'$. \qedhere}
		\end{figure}
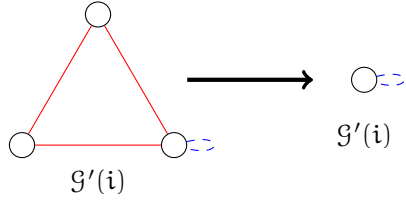
	\end{exem}
	
	\subsubsection{}
	
	By construction the resulting graph has no blue straight edges, and no red/green loop edges, but it need not be a projective crystallograph.
	Nonetheless this definition leads to the following compatibility.
	
	\begin{lemm}
		\label{lem:compatibility_quotients_projectifications}
		
		Projectifications and quotients of crystallographs are commutative operations: if $\mc G' \sse \mc G$ are nested crystallographs then there is an equality
		\begin{equation}
			\label{eq:graphs_comparison_quotients_projectifications}
			\mb P (\mc G) \bs \mb P (\mc G') = \mb P \bigl( \mc G \bs \mc G' \bigr) \, ,
		\end{equation}
		of trichromatic graphs.
	\end{lemm}

	\begin{proof}
		The set of nodes are the same, since (simply-laced) red components of $\mc G'$ are unaffected by operation $\mc G \mapsto \mb P(\mc G)$.		
		As for the adjacency/colouring, by the same token the simply-laced parts of both sides of~\eqref{eq:graphs_comparison_quotients_projectifications} coincide, so we need only show that the subsets of nodes with blue loop edges are the same.
		
		Now let $\mc G'(i) \sse \mc G'$ be a red component.
		By definition the left-hand side of~\eqref{eq:graphs_comparison_quotients_projectifications} has a blue loop edge $l = \Set{ \mc G'(i) }$ if and only if:
		\begin{itemize}
			\item there is a straight green edge $e \in \mc G_1$ within $\mc G'(i)$;
			
			\item there is a straight edge $e \in \mc G_1$ from $\mc G'(i)$ to a bichromatic component of $\mc G'$;
			
			\item there is a loop at a node of $\mc G'(i)$.
		\end{itemize}
		
		In all these cases there will be a loop $l = \Set{ \mc G'(i) }$ (of the same colour) in the quotient $\mc G \bs \mc G'$, reasoning as in the proofs of the lemmata in \S~\ref{sec:quotients}: this becomes blue after projectification---i.e. on the right-hand side of~\eqref{eq:graphs_comparison_quotients_projectifications}.
	\end{proof}
	
	\subsubsection{}

	Hence the quotients of projective crystallographs are exactly the projectifications of quasi-crystallographs, in view of Thm.~\ref{thm:quasi_crystallographs_classify_quotients} and Prop.~\ref{prop:quasi_crystallographs_are_quotients}.
	In turn the two ''exotic`` components of Cor.~\ref{cor:quasi_crystallographs_classification} only yield one ''exotic`` hyperplane arrangement, corresponding to a simply-laced complete bichromatic graph with blue loop edges glued at any \emph{proper} subset of nodes (cf. Rk.~\ref{rk:exotic_arrangement}).

	This leads to the final statement, classifying \emph{all} restricted root-hyperplane arrangements.
	
	\begin{theo}
		\label{thm:projectified_quotients_classify_restricted_arrangements}
		Let $\Phi' \sse \Phi \sse BC_n$ be nested root systems.
		Then the hyperplane arrangement of the restricted system~\eqref{eq:restricted_system} is isomorphic to a direct sum of hyperplane arrangements of classical type, or of (unique) ''exotic`` type $(B_r/C_r)D_{r+s}$, for integers $r,s \geq 0$---with $r + s \leq n$.
	\end{theo}
	
	\section{About closed subsystems}
	\label{sec:closed_subystems}
	
	\subsection{}
	
	In these two final sections we show how the main Def.~\ref{def:crystallographs} can be refined to classify distinguished root subsystems of classical type.
	
	First, it is possible to encode symmetric \emph{closed} subsets of roots of $BC_n$ via bichromatic graphs, with the colouring convention of \S~\ref{sec:quotient_crystallographic_obstruction} (cf.~\eqref{eq:closed_subset}):
	
	\begin{defi}
		\label{def:closed_bichromatic_graphs}
		
		A bichromatic graph $(\mc G,c)$ is \emph{closed} if it satisfies the ``simply-laced'' conditions of Figg.~\ref{fig:straight_edge_closure_condition_crystallographs_I},~\ref{fig:green_loop_edge_closure_crystallograph_II} and~\ref{fig:straight_edge_closure_condition_crystallographs_III}, and moreover:
		\begin{enumerate}
			\item if a node is incident to a green loop edge and a straight edge, then $\mc G$ contains the straight edge of opposite colour; 

			\begin{figure}[H]
				\centering
				\begin{tikzpicture}
					\node (a) at (0,0) {};
					\node (b) at (2,0) {};
					\node (c) [draw=none] at (2.5,0) {};
					\node (d) [draw=none] at (4.5,0) {};
					\node (e) at (5,0) {};
					\node (f) at (7,0) {};
					\graph{
						(a) --[yellow] (b),
						(c) ->[ultra thick] (d),
						(e) --[red,bend right=20] (f),
						(e) --[green,bend left=20] (f)
					};
					\path
					(a) edge [green,loop below] (a)
					(e) edge [green,loop below] (e);
				\end{tikzpicture}
				\caption{Green loop closure condition for closed bichromatic graphs (I).}
				\label{fig:green_loop_closure_condition_closed_graphs_I}
			\end{figure}
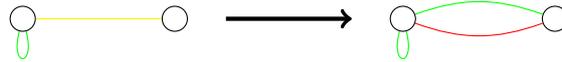
			
% 			\begin{figure}[H]
% 				\centering
% 				\begin{tikzpicture}
% 					\node (a) at (0,0) {};
% 					\node (b) at (2,0) {};
% 					\node (c) [draw=none] at (2.5,0) {};
% 					\node (d) [draw=none] at (4.5,0) {};
% 					\node (e) at (5,0) {};
% 					\node (f) at (7,0) {};
% 					\graph{
% 						(a) --[green] (b),
% 						(c) ->[ultra thick] (d),
% 						(e) --[red,bend right=20] (f),
% 						(e) --[green,bend left=20] (f)
% 					};
% 					\path
% 					(a) edge [green,loop below] (a)
% 					(e) edge [green,loop below] (e);
% 				\end{tikzpicture}
% 				\caption{Green loop closure condition for closed bichromatic graphs (II).}
% 				\label{fig:green_loop_closure_condition_closed_graphs_II}
% 			\end{figure}
			
			\item if $\mc G$ contains a double straight edge, then it contains the green loop edges at both ends;
			
			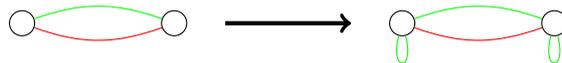
\begin{figure}[H]
				\centering
				\begin{tikzpicture}
					\node (a) at (0,0) {};
					\node (b) at (2,0) {};
					\node (c) [draw=none] at (2.5,0) {};
					\node (d) [draw=none] at (4.5,0) {};
					\node (e) at (5,0) {};
					\node (f) at (7,0) {};
					\graph{
						(a) --[red,bend right=20] (b),
						(a) --[green,bend left=20] (b),
						(c) ->[ultra thick] (d),
						(e) --[red,bend right=20] (f),
						(e) --[green,bend left=20] (f)
					};
					\path
					(e) edge [green,loop below] (e)
					(f) edge [green,loop below] (f);
				\end{tikzpicture}
				\caption{Green loop closure condition for closed bichromatic graphs (II).}
				\label{fig:green_loop_closure_condition_closed_graphs_II}
			\end{figure}
			
			\item if a node is incident at a red loop and a straight edge, then $\mc G$ contains the red loop at the opposite side;
			
			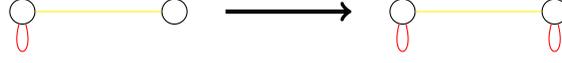
\begin{figure}[H]
				\centering
				\begin{tikzpicture}
					\node (a) at (0,0) {};
					\node (b) at (2,0) {};
					\node (c) [draw=none] at (2.5,0) {};
					\node (d) [draw=none] at (4.5,0) {};
					\node (e) at (5,0) {};
					\node (f) at (7,0) {};
					\graph{
						(a) --[yellow] (b),
						(c) ->[ultra thick] (d),
						(e) --[yellow] (f)
					};
					\path
					(a) edge [red,loop below] (a)
					(e) edge [red,loop below] (e)
					(f) edge [red,loop below] (f);
				\end{tikzpicture}
				\caption{Red loop closure condition for closed bichromatic graphs (I).}
				\label{fig:red_loop_closure_condition_closed_graphs_I}
			\end{figure}
			
% 			\begin{figure}[H]
% 				\centering
% 				\begin{tikzpicture}
% 					\node (a) at (0,0) {};
% 					\node (b) at (2,0) {};
% 					\node (c) [draw=none] at (2.5,0) {};
% 					\node (d) [draw=none] at (4.5,0) {};
% 					\node (e) at (5,0) {};
% 					\node (f) at (7,0) {};
% 					\graph{
% 						(a) --[green] (b),
% 						(c) ->[ultra thick] (d),
% 						(e) --[green] (f)
% 					};
% 					\path
% 					(a) edge [red,loop below] (a)
% 					(e) edge [red,loop below] (e)
% 					(f) edge [red,loop below] (f);
% 				\end{tikzpicture}
% 				\caption{Red loop closure condition for closed bichromatic graphs (II).}
% 				\label{fig:red_loop_closure_condition_closed_graphs_II}
% 			\end{figure}
			
			\item if two distinct nodes are incident at red loop edges, then $\mc G$ contains a double straight edge incident at those nodes;
			
			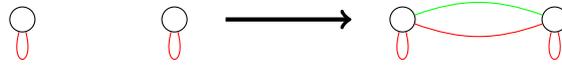
\begin{figure}[H]
				\centering
				\begin{tikzpicture}
					\node (a) at (0,0) {};
					\node (b) at (2,0) {};
					\node (c) [draw=none] at (2.5,0) {};
					\node (d) [draw=none] at (4.5,0) {};
					\node (e) at (5,0) {};
					\node (f) at (7,0) {};
					\graph{
						(c) ->[ultra thick] (d),
						(e) --[red,bend right=20] (f),
						(e) --[green,bend left=20] (f)
					};
					\path
					(a) edge [red,loop below] (a)
					(b) edge [red,loop below] (b)
					(e) edge [red,loop below] (e)
					(f) edge [red,loop below] (f);
				\end{tikzpicture}
				\caption{Red loop closure condition for closed bichromatic graphs (II).}
				\label{fig:red_loop_closure_condition_closed_graphs_II}
			\end{figure}
			
			\item if a node is incident at a red loop edge, then it is also incident at a green loop edge.
		\end{enumerate}
	\end{defi}
	
	\subsubsection{}
	
	Then the main correspondence yields the following.
	
	\begin{theo}[Compare with Thm.~\ref{thm:crystallographs_classify_root_systems}]
		\label{thm:closed_bichromatic_graphs_classify_closed_symmetric_subsets}
		
		The association $(\mc G,c) \mapsto \Phi^{\mc G,c}$ restricts to a bijection between closed bichromatic graphs on $n \geq 1$ nodes and symmetric closed subsets of $BC_n$.
	\end{theo}
	
	\begin{proof}
		The conditions in Def.~\ref{def:closed_bichromatic_graphs} are the translation of the constraint~\eqref{eq:closed_subset}.
		Indeed one has
		\begin{equation}
			e^-_{ij} + e^-_{jk} = e^-_{jk}, \quad e^-_{ij} + e^+_{jk} = e^+_{ik}, \quad e^+_{ij} - e^+_{jk} = e^-_{ik},
		\end{equation} 
		and 
		\begin{equation}
			e^-_{ij} + 2e_i = e^+_{ij}, \quad e^+_{ij} - 2e_j = e^-_{ij},
		\end{equation} 
		and 
		\begin{equation}
			e^-_{ij} + e^+_{ij} = 2e_i, \quad e^-_{ji} + e^+_{ij} = 2e_j, 
		\end{equation} 
		and 
		\begin{equation}
			e^-_{ij} + e_j = e_i, \quad e^-_{ij} - e_i = -e_j,
		\end{equation} 
		and 
		\begin{equation}
			e^+_{ij} - e_i = e_j, \quad e^+_{ij} - e_j = e_i,
		\end{equation} 
		and the rest follows by definition.
		
		To check no other sum yields an element of $BC_n$ one can (repeatedly) use the following elementary fact: let $S \sse V$ be a set of linearly independent vectors, $v_1,\dc,v_k \in S$ \emph{distinct} elements, and $\lambda_1,\dc,\lambda_k \in \mb C$ \emph{nonvanishing} numbers; then the vector $v = \sum_i \lambda_iv_i \in V$ is \emph{not} a linear combination of a smaller number of elements of $S$.
		Thus e.g. if $\lambda,\mu \in \set{\pm 1}$ then the vector $v = \lambda e^-_{ij} + \mu e^-_{kl}$ cannot be an element of $BC_n$ if $\set{i,j} \cap \set{k,l} = \vn$---and it is easy to see $\set{i,j} = \set{k,l}$ cannot work.
		Suppose instead the intersection is a singleton, say $j = k$ and $i \neq l$, so that $v = \lambda e_i + (\mu - \lambda) e_j - \mu e_l$: by the above we must have $\lambda = \mu$, and the only roots one finds are $v = \pm e^-_{il}$; etc.
	\end{proof}
	
	\subsection{}
	
	As expected there is a partial superposition with the crystallographic conditions, since few Weyl reflections act via sums, notably $\sigma^-_{ij}(e^-_{jk}) = e^-_{ik}$ in type $A$.
	Then one can define a \emph{closed crystallograph} by intersecting the conditions of Def.~\ref{def:crystallographs} and~\ref{def:closed_bichromatic_graphs}, finding:
	
	\begin{coro}
		The association $(\mc G,c) \mapsto \Phi^{\mc G,c}$ restricts to a bijection between closed crystallographs on $n \geq 1$ nodes and closed root subsystems of $BC_n$.
	\end{coro}
	
	\begin{proof}
		This follows from Thm.~\ref{thm:crystallographs_classify_root_systems} and~\ref{thm:closed_bichromatic_graphs_classify_closed_symmetric_subsets}.
	\end{proof}

	\subsubsection{}
	
	This builds towards a graph-theoretic classification of all closed root subsystems of the classical root systems, along the lines of Cor.~\ref{cor:root_system_classification}.
	Of course one must take care that e.g. $B_n \sse BC_n$ is \emph{not} a closed subset, so in general a variation of Def.~\ref{def:closed_bichromatic_graphs} is needed.
	
	Rather than pursuing this we shall consider a stronger (and more symmetric) condition which is much relevant to hyperplane complements and wild isomonodromic deformations.
	
	\section{About Levi subsystems}
	\label{sec:levi_subsystems}
	
	\subsection{}
	
	We finally wish to encode the algebraic condition of \emph{Levi} subsystems into a bichromatic graph, cf.~\eqref{eq:levi_subsystem}.
	
	\begin{defi}
		\label{def:levigraph}
		
		A \emph{Levigraph} (\emph{of type} $BC$) is a crystallograph $(\mc G,c)$ satisfying the following conditions:
		\begin{enumerate}			
			\item if a node is incident to a loop edge, then it is incident to two loop edges;
			
			\item if $\mc G$ contains a double straight edge, then it contains the two loop edges at both ends;
			
			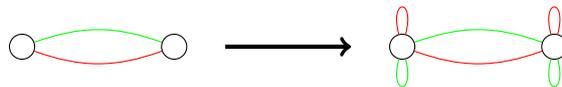
\begin{figure}[H]
				\centering
				\begin{tikzpicture}
					\node (a) at (0,0) {};
					\node (b) at (2,0) {};
					\node (c) [draw=none] at (2.5,0) {};
					\node (d) [draw=none] at (4.5,0) {};
					\node (e) at (5,0) {};
					\node (f) at (7,0) {};
					\graph{
						(a) --[red,bend right=20] (b),
						(a) --[green,bend left=20] (b),
						(c) ->[ultra thick] (d),
						(e) --[red,bend right=20] (f),
						(e) --[green,bend left=20] (f)
					};
					\path
					(e) edge [green,loop below] (e)
					(e) edge [red,loop above] (e)
					(f) edge [green,loop below] (f)
					(f) edge [red,loop above] (f);
				\end{tikzpicture}
				\caption{Loop closure condition for Levigraphs (I).}
				\label{fig:loop_closure_condition_levigraphs_I}
			\end{figure}
			
			\item if two distinct nodes are incident at loop edges, then $\mc G$ contains a double straight edge incident at those nodes.
			
			\begin{figure}[H]
				\centering
				\begin{tikzpicture}
					\node (a) at (0,0) {};
					\node (b) at (2,0) {};
					\node (c) [draw=none] at (2.5,0) {};
					\node (d) [draw=none] at (4.5,0) {};
					\node (e) at (5,0) {};
					\node (f) at (7,0) {};
					\graph{
						(c) ->[ultra thick] (d),
						(e) --[red,bend right=20] (f),
						(e) --[green,bend left=20] (f)
					};
					\path
					(a) edge [yellow,loop below] (a)
					(b) edge [yellow,loop below] (b)
					(e) edge [yellow,loop below] (e)
					(f) edge [yellow,loop below] (f);
				\end{tikzpicture}
				\caption{Loop closure condition for Levigraphs (II).}
				\label{fig:loop_closure_condition_levigraphs_II}
			\end{figure}
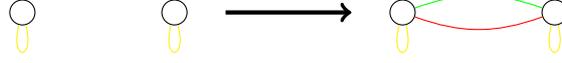
		\end{enumerate}
	\end{defi}
	
	\begin{rema}
		The conditions of Def.~\ref{def:levigraph} are ``minimal'', and some more follow by iteration.
		In particular the following one is important in type $D$ (cf. Def.~\ref{def:levigraph_type_D}): if two distinct pairs of nodes are incident at a double straight edge, then $\mc G$ contains the complete simply-laced bichromatic graph on the corresponding four nodes.
		
		\begin{figure}[H]
			\centering
			\begin{tikzpicture}
				\node (a) at (0,0) {};
				\node (b) at (0,2) {};
				\node (c) at (2,0) {};
				\node (d) at (2,2) {};				
				\node (e) [draw=none] at (2.5,1) {};
				\node (f) [draw=none] at (4.5,1) {};
				\node (g) at (5,0) {};
				\node (h) at (5,2) {};
				\node (i) at (7,0) {};
				\node (j) at (7,2) {};
				
				\graph{
					(a) --[red,bend right=20] (b),
					(a) --[green,bend left=20] (b),
					
					(c) --[green,bend right=20] (d),
					(c) --[red,bend left=20] (d),

					(e) ->[ultra thick] (f),
					
					(g) --[red,bend right=20] (h),
					(g) --[green,bend left=20] (h),
					(g) --[red,bend right=20] (i),
					(g) --[green,bend left=20] (i),
					(i) --[green,bend right=20] (j),
					(i) --[red,bend left=20] (j),
					(h) --[green,bend right=20] (j),
					(h) --[red,bend left=20] (j),
				};
			\end{tikzpicture}
			\caption{Double straight edge closure condition for Levigraphs.}
			\label{fig:double_straight_edge_closure_condition_levigraphs}
		\end{figure}
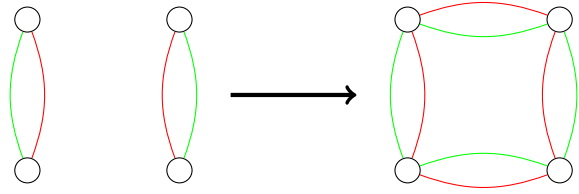
		
		(The point is this does \emph{not} appeal to loop edges.)
	\end{rema}
	
	\subsubsection{}
	
	By definition all Levigraphs are (closed) crystallographs, but not conversely.
	Indeed, the main correspondence yields:
	
	\begin{theo}[Compare with Thm.~\ref{thm:crystallographs_classify_root_systems} and~\ref{thm:closed_bichromatic_graphs_classify_closed_symmetric_subsets}]
		\label{thm:levigraphs_classify_levi_subsystems}
		
		The association $(\mc G,c) \mapsto \Phi^{\mc G,c}$ restricts to a bijection between Levigraphs on $n \geq 1$ nodes and Levi subsystems of $BC_n$.		
	\end{theo}
	
	\begin{proof}
		The conditions of Def.~\ref{def:levigraph} translate all possible ways in which a root $\gamma \in BC_n$ can be expressed as 
		\begin{equation}
			\gamma = \alpha + \lambda \beta, \qquad \alpha,\beta \in BC_n, \quad \lambda \in \mb Q.		                                                                                      
		\end{equation} 
		(To show no other root arises one can again appeal to the elementary fact stated in the proof of Thm.~\ref{thm:closed_bichromatic_graphs_classify_closed_symmetric_subsets}.)
		Generic linear combinations iteratively follow, as prescribed by~\eqref{eq:levi_subsystem}; e.g. 
		\begin{equation}
			e_i = \frac 1 2 (e^-_{ij} + e^+_{ij}) = \Bigl( e^-_{ij} - \frac 1 2 e^-_{ij} \Bigr) + \frac 1 2 e^+_{ij},
		\end{equation} 
		in addition to~\eqref{eq:closed_subset}, and one can freely pass from $e_i$ to $2e_i$ up to multiplying by $\lambda \in \set{ 2^{\pm 1} }$, etc.
	\end{proof}

	\subsection{}
	
	Now a variation of Thm.~\ref{thm:crystallograph_classification} yields the following characterisation:
	
	\begin{coro}
		\label{cor:levigraph_classification}
		
		Let $(\mc G,c)$ be a Levigraph on $n \geq 1$ nodes. 
		Then $\mc G$ is a disjoint union of arbitrarily many (Levi)graphs of type $\mc G^{A_{m-1}}$ and $\mc G^{d_1,d_2}$, and at most a \emph{single} (Levi)graph of type $\mc G^{BC_m}$---where $m,d_1+d_2 \leq n$.
	\end{coro}
	
	\begin{proof}
		By Thm.~\ref{thm:crystallograph_classification}, the crystallographic conditions yield connected components of classical/bipartite type.
		
		First, in view of Fig.~\ref{fig:loop_closure_condition_levigraphs_I}, all components of type $B$, $C$, or $D$ are actually of type $BC$; second, by Fig.~\ref{fig:loop_closure_condition_levigraphs_II} they are all fused into a single component---of type $BC$.
		
		Nothing more happens to any type-$A$/bipartite component, as no further double straight edges or loop edges can be created.
	\end{proof}

	\subsubsection{}
	
	In view of Thm.~\ref{thm:levigraphs_classify_levi_subsystems} and Lem.~\ref{lem:bipartite_graph_is_type_A}, this immediately yields a proof of the following statement:
	
	\begin{coro}
		\label{cor:classification_levi_subsystems_type_BC}
		
		Suppose $\Phi \sse BC_n$ is a Levi subsystem.
		Then $\Phi$ is equivalent to a direct sum of arbitrarily many type-$A$ irreducible root systems, plus at most a single irreducible root system of type $BC$.
	\end{coro}
	
	\subsection{Classical Levi subsystems}
	
	We now wish to classify the Levi subsystems for all classical root systems, refining the statement of Cor.~\ref{cor:classification_levi_subsystems_type_BC}.
	This means varying the closure condition~\eqref{eq:levi_subsystem} by changing the ``ambient'' root system, i.e. e.g. imposing that $\spann_{\mb C}(\Phi) \cap A_{n-1} = \Phi$ for a subset $\Phi \sse A_{n-1}$, etc.
	
	The result is the following list of definitions.
	
	\begin{defi}
		A \emph{type}-$A$ Levigraph is just a type-$A$ crystallograph.
	\end{defi}
	
	\begin{defi}
		A \emph{type}-$B$ (resp. \emph{type}-$C$) Levigraph is a crystallograph $\mc G \sse \mc G^B_n$ (resp. $\mc G \sse \mc G^C_n$) satisfying the following conditions:
		\begin{enumerate}
			\item if $\mc G$ contains a double straight edge, then it contains the two red loop edges (resp. green loop edges) at both ends, cf. Fig.~\ref{fig:loop_closure_condition_levigraphs_I};
			
			\item if two distinct nodes are incident at red loop edges (resp. at green loop edges), then $\mc G$ contains a double straight edge incident at those nodes, cf. Fig.~\ref{fig:loop_closure_condition_levigraphs_II}.
		\end{enumerate}
	\end{defi}
	
	\begin{defi}
		\label{def:levigraph_type_D}
		
		A \emph{type}-$D$ Levigraph is a crystallograph $\mc G \sse \mc G^D_n$ satisfying the condition of Fig.~\ref{fig:double_straight_edge_closure_condition_levigraphs}.
	\end{defi}
	
	\subsubsection{}
	
	Again one can prove the analogue of Thm.~\ref{thm:levigraphs_classify_levi_subsystems}, i.e. that type-$A$ (resp. type-$B$, type-$C$, type-$D$) Levigraphs on $n \geq 1$ nodes correspond bijectively to Levi subsystems of $A_{n-1}$ (resp. $B_n$, $C_n$, $D_n$).
	This is simply because the above list of definitions correspond to restricting the closure conditions of Def.~\ref{def:levigraph} to each classical type.
	Then using Cor.~\ref{cor:levigraph_classification} it is now easy to classify such Levigraphs of classical types; e.g. all type-$B$ Levigraphs on $n \geq 1$ nodes are disjoint unions of arbitrarily many (Levi)graphs of type $\mc G^A_{m-1}$ and $\mc G^{d_1,d_2}$, plus at most a \emph{single} (Levi)graph of type $\mc G^B_m$---where $m,d_1+d_2 \leq n$.	
	
	\subsection{}
	
	Finally, yet another variation of Cor.~\ref{cor:root_system_classification} proves the following:
	
	\begin{coro}
		Let $\Phi' \sse \Phi \sse BC_n$ be nested root systems, and suppose $\Phi'$ is a Levi subsystem of $\Phi$.
		Then:
		\begin{itemize}
			\item if $\Phi = A_{n-1}$, $\Phi'$ is equivalent to a direct sum of (arbitrarily many) type-$A$ irreducible root systems;\fn{
				Proving the standard fact that all root subsystems of $A_{n-1}$ are Levi.}
			
			\item if $\Phi = B_n$ (resp. $\Phi = C_n$, $\Phi = D_n$), $\Phi'$ is equivalent to a direct sum of arbitrarily many type-$A$ irreducible root systems, plus at most a \emph{single} irreducible root system of type $B$ (resp. type $C$, type$D$).
		\end{itemize}
	\end{coro}

	This is in accordance with the standard characterisations of Levi subsystems in terms of subdiagrams of the Dynkin diagram of the associated simple Lie algebra, but contrary to the Dynkin diagram it does \emph{not} use the choice of a base of simple roots.

	\section*{Outlook}
	
	The Borel--de Siebenthal theory~\cite{borel_de_siebenthal_1949_les_sous_groupes_fermes_de_rang_maximum_des_groupes_de_lie_clos} classifies connected closed subgroups of maximal rank, inside connected compact Lie groups $G$---with a given maximal torus $T \sse G$.
	This can be expressed in terms of certain operations on the extended Dynkin diagram of $\mf g = \Lie(G)$, and leads in particular to a classification of maximal closed subsystems of $\Phi_{\mf g}$, cf.~\cite{kane_2001_reflection_groups_and_invariant_theory}.
	
	The present graph-theoretic construction instead does \emph{not} rely on (extended) Dynkin diagrams, but rather replaces them with new combinatorial objects that retain more information.
	In particular we thus classify root subsystems of \emph{any} rank, not just maximal ones (cf.~\cite{wallach_1968_on_maximal_subsystems_of_root_systems}); in principle a Lie-group theoretic analogue should exist.
	Moreover in our viewpoint maximal root subsystems simply correspond to maximal sub-crystallographs, so can still be characterised in combinatorial fashion, encompassing e.g. maximal simply-laced complete subgraphs in type $A$.
	
	\section*{Acknowledgements}
	
	We thank D.~Schwein for their interest in this paper, and for helpful discussions.
	
	\appendix 
	
	\section{Notations/conventions}
	\label{sec:appendix}
	
	\subsection*{Graphs}
	
	Here it is convenient to think about graphs as follows.
	
	A \emph{multiset} is a set $S$ with a multiplicity function $m \cl S \to \mb{Z}_{\geq 1}$.
	If the underlying set $S$ is finite then $\abs m = \sum_{i \in S} m(i) \in \mb{Z}_{\geq 1}$ is the \emph{cardinality} of $(S,m)$.

	A \emph{graph} is a triple $\mc G = (\mc G_0,\mc G_1,m)$ consisting of a set of \emph{nodes} $\mc G_0$ and a multiset $(\mc G_1,m)$ of (unoriented) \emph{edges}.
	In turn an edge is a multiset $(e,m_e)$ with $\abs{m_e} = 2$, and such that its underlying set lies in $\mc G_0$: it is a \emph{loop edge} if $e$ is a singleton, else it is a \emph{straight edge}.
	The forgetful function $\mc G_1 \to \mc{P}(\mc G_0)$ defined by $(e,m_e) \mapsto e$ is the \emph{incidence} of the graph: the edge $(e,e_m)$ is incident at the node(s) of $e \sse \mc G_0$, and conversely.

	A straight edge incident at the nodes $i \neq j \in \mc G_0$ is written $e = \set{i,j} \in \mc G_1$---with $m(i) = m(j) = 1$---, while a loop edge incident at $k \in \mc G_0$ is $l = \set{k}$---with $m(k) = 2$.
	These are depicted as customary:
	
	\begin{figure}[H]
		\centering
		$e = $\begin{tikzpicture}[baseline=-1mm]
			\node (i) [label=$i$] {};
			\node (j) [right of=i,label=$j$] {};
			\graph{ 
				(i) -- (j)
			};
		\end{tikzpicture} \qquad \qquad 
		$l =$\begin{tikzpicture}[baseline=-1mm]
			\node (k) [label=$k$] {};
			\path
				(k) edge [loop below] (k);
		\end{tikzpicture}
		\caption{Straight and loop edges}
	\end{figure}
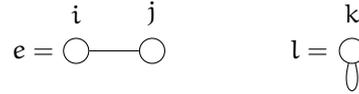
	
	A \emph{subgraph} $\mc G'$ of $\mc G$ is obtained by removing some node or edge: if a node is removed all its adjacent edges are also dropped, and the resulting partial order is denoted $\mc G' \sse \mc G$.
	
	The disconnected/disjoint union of two graphs $\mc G$ and $\mc G'$ is denoted $\mc G \mathsmaller{\coprod} \mc G'$, and underlies the disjoint union $\mc G_0 \mathsmaller{\coprod} \mc G'_0$ of the sets of nodes.
	
	\subsection*{Root systems and Weyl groups}
	
	For an integer $n \geq 1$ let $V = \mb{C}^n$ with canonical basis $(e_1,\dotsc,e_n)$---and scalar product $(e_i \mid e_j) = \delta_{ij}$, whenever helpful.
	There are finite subsets 
	\begin{equation*}
		A_{n-1} \sse D_n \sse B_n,C_n \sse BC_n \sse V \setminus \set{0} \, ,
	\end{equation*}
	of nested (crystallographic, irreducible) root systems, defined as follows.
	If $e^{\pm}_{ij} \ceqq e_i \pm e_j \in V$ for $i,j \in \ul n = \set{1,\dc,n}$, then
	\begin{equation*}
		A_{n-1} = \Set{ e^-_{ij} | i \neq j \in \ul n } \, , \qquad D_n = A_{n-1} \cup \Set{ \pm e^+_{ij} | i \neq j \in \ul n } \, ,
	\end{equation*}
	and 
	\begin{equation*}
		B_n = D_n \cup \Set{ e_i | i \in \ul n } \, , \qquad C_n = D_n \cup \Set{ 2e_i | i \in \ul n } \, ,
	\end{equation*}
	and finally $BC_n = B_n \cup C_n$.

	The single-letter families are the classical \emph{reduced} irreducible root systems, while $BC_n$ is a \emph{nonreduced} one---the unique one of rank $n$, up to isomorphism, cf.~\cite[Ch.~VI]{bourbaki_1968_groupes_et_algebres_de_lie_chapitres_4_6} and~\cite{humphreys_1972_introduction_to_lie_algebras_and_representation_theory}.
	
	There is a dual viewpoint where $BC_n \sse V^{\dual}$, consisting of covectors naturally associated with the dual basis $\alpha_i = e_i^{\dual} \in V^{\dual}$---i.e. taking dual/inverse root systems.
	Whenever helpful we equip $V^{\dual}$ with the push-forward scalar product along the vector space isomorphism $( \cdot \mid \cdot )^{\flat} \cl V \to V^{\dual}$, and we identify (canonically) the Weyl groups $W(\Phi) \sse \GL(V)$ and $W(\Phi^{\dual}) \sse \GL(V^{\dual})$ by $w \mapsto \prescript{t}{}w^{-1}$.

	A root \emph{subsystem} of a root system $\Phi \sse V$ is a subset $\Phi' \sse \Phi$ such that $\sigma_{\alpha}(\Phi') \sse \Phi'$ for $\alpha \in \Phi'$, where $\sigma_{\alpha} \in W$ is the correspond Weyl reflection.
	The Weyl group permutes root subsystems: two such are \emph{equivalent} if they lie in the same Weyl-group orbit.
	This is a stricter relation than that of isomorphism, since there are in general ''outer`` automorphisms of $\Phi$---the isomorphisms of the Dynkin diagram---which also preserve root subsystems.
	
	A subset $\Phi' \sse \Phi$ is \emph{closed} if 
	\begin{equation}
		\label{eq:closed_subset}
		(\Phi' + \Phi') \cap \Phi = \Phi',
	\end{equation} 
	i.e. if it contains the sum of all its elements, provided the result is still a root.
	E.g. the subsystem of short roots of $B_2$, which is a copy of $D_2 = A_1 \oplus A_1$, is \emph{not} closed (in $B_2$); but the subsystem of long roots is.
	
	A root subsystem $\Phi' \sse \Phi$ is a \emph{Levi} subsystem if it is obtained by intersecting $\Phi$ with a vector subspace of $\spann_{\mb C}(\Phi)$, or equivalently if 
	\begin{equation}
		\label{eq:levi_subsystem}
		\spann_{\mb C}(\Phi') \cap \Phi = \Phi',
	\end{equation} 
	i.e. it contains the linear combinations of all its elements, provided the result is still a root.
	(It is enough to take $\mb Q$-linear combinations.)
	This implies $\Phi'$ is a closed root subsystem of $\Phi$, but not conversely: e.g. the closed subsystem of long roots of $B_2$, which is another copy of $D_2$, is \emph{not} Levi.
	Further, this is equivalent to asking that the root hyperplane complement~\eqref{eq:restricted_complement_2} is nonempty, so this leads to the topological spaces whose fundamental groups are (pure) local wild mapping class groups.
	
	\subsection*{Hyperplane arrangements}
	
	A (linear) \emph{hyperplane arrangement} $\mc H$ in a (finite-dimensional) complex vector space $V$ is a set of linear hyperplanes $H \sse V$: in this paper we only consider finite such arrangements, and we identify hyperplanes in $V$ with lines in the dual space via $\lambda \mapsto \Ker(\lambda) \sse V$, for $\lambda \in V^{\dual} \sm \set{0}$.
	
	An isomorphism of two hyperplane arrangements $\mc H, \mc H'$ in $V$ is a linear automorphism $w \in \GL(V)$ such that $\Set{ w(H) | H \in \mc H} = \mc H'$.
	
	A \emph{sub-arrangement} of $\mc H \sse \mb P \bigl( V^{\dual} \bigr)$ is a subset $\mc H' \sse \mc H$ of hyperplanes of $\mc H$.
	
	The \emph{direct sum} of hyperplane arrangements $\mc H_i \sse \mb P\bigl( V_i^{\dual} \bigr)$, $i = 1,2$, is the set of hyperplanes 
	\begin{equation*}
		H_1 \oplus V_2, \, V_1 \oplus H_2 \sse V_1 \oplus V_2 \, , \qquad H_i \in \mc H_i \, .
	\end{equation*}
	
	The hyperplane arrangement of a root-system $\Phi \sse \mf t^{\dual}$---viz. a ''root-hyperplane arrangement``, for short---is the set 
	\begin{equation*}
		\mc H^{\Phi} \ceqq \Set{ \Ker(\alpha) | \alpha \in \Phi } \simeq \Phi \bs \mb C^{\times} \sse \mb P \bigl( \mf t^{\dual} \bigr) \, .
	\end{equation*}
	Such hyperplane arrangements are said to be \emph{crystallographic}.
	
	\backmatter 
	
	%Bibliography: BibTeX
	\bibliographystyle{amsalpha}
		\bibliography{/home/gabriele/Desktop/bibliography_macros/bibliography}	
\end{document}